\documentclass[12pt]{article}
\usepackage{amsmath}
\usepackage{amssymb}
\numberwithin{equation}{section}
\usepackage{epic}

\title{
Virasoro action on %Fourier series, Schur functions
 Schur function expansions, skew
Young tableaux and random walks
 }

\author{
M. Adler\thanks{ Department of Mathematics, Brandeis
University, Waltham, Mass 02454, USA. E-mail:
adler@brandeis.edu.  The support of a National Science
Foundation grant \# DMS-01-00782 is gratefully
acknowledged.}~~~~~~ P. van Moerbeke\thanks{ Department
of Mathematics, Universit\'e de Louvain, 1348
Louvain-la-Neuve, Belgium and Clay Mathematics
Institute, One Bow Street, Cambridge, MA 02138, USA.
E-mail: vanmoerbeke@math.ucl.ac.be and @brandeis.edu .
The support of a National Science Foundation grant \#
DMS-01-00782, European Science Foundation, Nato, FNRS
and Francqui Foundation grants is gratefully
acknowledged.}}

\date{}

\catcode`\@=11
\let\c@equation=\relax
\newcounter{equation}[subsection]

\catcode`\@=12

\newcommand{\MAT}[1]{\left(\begin{array}{*#1c}}
\newcommand{\mat}{\end{array}\right)}
\newcommand{\qed}{\leavevmode\unskip\nobreak\penalty200\hskip2pt\null
\nobreak\hfill\rule{1.1ex}{1.1ex}%\parfillskip=0pt
\medbreak }

\newcommand{\pp}{\ldots}

\newcommand{\LR}{{\cal L}}

\newcommand{\VR}{{\cal V}}

\newcommand{\BV}{{\mathbb V}}

\newcommand{\BZ}{{\mathbb Z}}

\newcommand{\iy}{\infty}
\newcommand{\pl}{\partial}
\newcommand{\al}{\alpha}

\newcommand{\gs}{{\bf s}}
\newcommand{\no}{\nonumber}

\newcommand{\tx}{\tilde x}
\newcommand{\ty}{\tilde y}
\newcommand{\ba}{{\backslash}}

\newcommand{\acht}{{\rm ht}}

\newcommand{\shade}{/\!\!/\!\!/\!\!/\!\!%/\!\!
}
\newcommand{\bshade}{\backslash\!\!
\backslash\!\!\backslash\!\!\backslash%\!\!\backslash~~~
}

\newenvironment{proof}{\medskip\noindent{\it Proof:\/} }{\qed}
\newenvironment
        {remark}{\medskip\noindent\underline{\it Remark:\/} }{\medbreak}
\newenvironment
        {example}{\medskip\noindent\underline{\it Example:\/} }{\medbreak}

\newcommand{\la}{\langle}
\newcommand{\ra}{\rangle}

\newcommand{\dt}{\delta}
\newcommand{\Dt}{\Delta}
 \newcommand{\vr}{\varepsilon}

\newcommand{\BR}{{\mathbb R}}
\newcommand{\lb}{\lambda}
\newcommand{\Lb}{\Lambda}

\def\be#1\ee{\begin{equation}#1\end{equation}}
\def\bea#1\eea{\begin{eqnarray}#1\end{eqnarray}}
\def\bean#1\eean{\begin{eqnarray*}#1\end{eqnarray*}}

\newcommand{\Tr}{\operatorname{\rm Tr}}

\newtheorem{definition}{Definition}[section]
\newtheorem{theorem}[definition]{Theorem}
\newtheorem{lemma}[definition]{Lemma}
\newtheorem{corollary}[definition]{Corollary}
\newtheorem{proposition}[definition]{Proposition}

\catcode `!=11

\newdimen\squaresize
\newdimen\thickness
\newdimen\Thickness
\newdimen\ll! \newdimen \uu! \newdimen\dd! \newdimen \rr! \newdimen
\temp!

\def\sq!#1#2#3#4#5{%
\ll!=#1 \uu!=#2 \dd!=#3 \rr!=#4
\setbox0=\hbox{%
 \temp!=\squaresize\advance\temp! by .5\uu!
 \rlap{\kern -.5\ll!
 \vbox{\hrule height \temp! width#1 depth .5\dd!}}%
 \temp!=\squaresize\advance\temp! by -.5\uu!
 \rlap{\raise\temp!
 \vbox{\hrule height #2 width \squaresize}}%
 \rlap{\raise -.5\dd!
 \vbox{\hrule height #3 width \squaresize}}%
 \temp!=\squaresize\advance\temp! by .5\uu!
 \rlap{\kern \squaresize \kern-.5\rr!
 \vbox{\hrule height \temp! width#4 depth .5\dd!}}%
 \rlap{\kern .5\squaresize\raise .5\squaresize
 \vbox to 0pt{\vss\hbox to 0pt{\hss $#5$\hss}\vss}}%
}%end of \hbox
 \ht0=0pt \dp0=0pt \box0
}%end of \sq!

\def\vsq!#1#2#3#4#5\endvsq!{\vbox to \squaresize{\hrule
width\squaresize height 0pt%
\vss\sq!{#1}{#2}{#3}{#4}{#5}}}

\newdimen \LL! \newdimen \UU! \newdimen \DD! \newdimen \RR!

\def\vvsq!{\futurelet\next\vvvsq!}
\def\vvvsq!{\relax
  \ifx     \next l\LL!=\Thickness \let\continue=\skipnexttoken!
  \else\ifx\next u\UU!=\Thickness \let\continue=\skipnexttoken!
  \else\ifx\next d\DD!=\Thickness \let\continue=\skipnexttoken!
  \else\ifx\next r\RR!=\Thickness \let\continue=\skipnexttoken!
  \else\def\continue{\vsq!\LL!\UU!\DD!\RR!}%
  \fi\fi\fi\fi
  \continue}

\def\skipnexttoken!#1{\vvsq!}

\def\place#1#2#3{\vbox to 0pt{\vss
\rlap{\kern#1\squaresize
  \raise#2\squaresize\hbox{$#3$}}
\vss}}

\def\Young#1{\LL!=\thickness \UU!=\thickness \DD! = \thickness \RR! =
\thickness \vbox{\smallskip\offinterlineskip
\halign{&\vvsq! ##
\endvsq!\cr #1}}}

\def\blank{\omit\hskip\squaresize}
\catcode `!=12

\begin{document}
\maketitle

%\subsection{Permutations and integrals over groups}
%\begin{abstract}

%It is known that some matrix integrals over $U(n)$
%satisfy an sl(2,$\BR$)-algebra of Virasoro constraints.
%Acting with these Virasoro generators on 2-dimensional
%Schur function expansion leads to difference relations
%on the coefficients of this expansion. These difference
%relations, set equal to zero, are precisely the backward
%and forward equations for non-intersecting random walks.
%The transition probabilities for these random walks
%appear as the coefficients of an expansion of a
%$U(n)$-matrix integrals (of the type above), by
%inserting in the integral the product of two Schur
%polynomials associated with two partitions; the latter
%are specified by the initial and final positions of the
%non-intersecting random walk. An essential ingredient in
%this work is the generalization of the
%Murnaghan-Nakayama rule to the action of Virasoro on
%Schur polynomials.

%\end{abstract}

\tableofcontents

%\newpage
\section{Introduction}

{\em The aim of this paper is to show that the Virasoro
action on two-dimensional Fourier series in Schur
polynomials and the backward/forward equations for
random walks are close allies!}

%{\em The aim of this paper is to show that the Virasoro
%algebra and the backward/forward equations for random
%walks are close allies, via two-dimensional Fourier
%analysis on symmetric functions!}

Some matrix integrals over $U(n)$ are known to satisfy a
sl(2,$\BR$)-algebra of Virasoro constraints \cite{AvM2}
 ($k=-1,0,1$)
\be \BV_k(t,s)
  \int_{U(n)}
 % e^{\sum_1^{\iy}\Tr(t_j^{(0)}M^j-s_j^{(0)}\bar M^j)}
  e^{\Tr\VR(M,\bar M)}dM=0
 ,~ \mbox{with}~~
   \VR(x,y):=\sum_1^{\iy}(t_j x^j-s_jy^j)
\label{MatrixIntegral} \ee
where

\vspace{-1.3cm}

\begin{eqnarray}
\BV_{-1}(t,s)&=&V_{-1}(t)-V_1(s)+n\left(t_1+\frac{\pl}{\pl
s_{1}}\right)\nonumber\\
 \BV_{0}(t,s)&=&V_0(t)-V_0(s)
 \label{BVirasoroNotation}\\
%{\cal D}_{1}I_n
 \BV_{1}(t,s)&=&-V_{-1}(s)+V_1(t)
  +n\left(s_1+\frac{\pl}{\pl t_1}
\right),\no
\end{eqnarray}
involving standard Virasoro operators\footnote{For
$k=-1,0,1$, they take on the following simple form:
 \be
V_k(t):=\sum_{i\geq \max
(k+1,1)}(i-k)t_{i-k}\frac{\pl}{\pl t_i},
~~~~~~\mbox{for}~ k=-1,0,1.
 \label{VirasoroNotation}\ee} in the variables
$t_1,t_2,\ldots
 $ for all $k \in \BZ$,
 \be
V_k(t)=\frac{1}{2}\sum_{i+j=k}\frac{\pl^2}{\pl t_i\pl
t_j}+\sum_{-i+j=k}it_i\frac{\pl}{\pl
t_j}+\frac{1}{2}\sum_{-i-j=k}(it_i)(jt_j).
\label{VirasoroIntro}\ee
Appropriate shifts of the $t_i$'s and $s_i$'s in the
matrix integral appearing in (\ref{MatrixIntegral}) lead
to the following matrix integrals, already considered by
several authors:
\cite{Gessel1,Borodin-Okounkov,Borodin-Olshanski, TW2}
$$
\begin{array}{ll|l}

&\mbox{shifts}  &  \mbox{matrix integrals}   \\

\hline

1:&\begin{array}{l}
     {it_i\mapsto it_i+z
     \delta_{i1}}  \\
    {is_i\mapsto is_i-z\delta_{i1}}
    \end{array}
 &
 I_1=\displaystyle{\int_{U(n)}}
   e^{z\Tr(M+\bar M)}e^{\Tr\VR(M)}dM  \\ &&\\

2:&\begin{array}{l}
  {it_i\mapsto it_i-q(-1)^i}
       \\
       {is_i\mapsto is_i- z\dt_{i1}}
       \end{array}
   &
  I_2= \displaystyle{\int_{U(n)}}
   \det(I+M)^q e^{z\Tr \bar M}e^{\Tr\VR(M)}dM      \\&& \\

 3:&\begin{array}{l}
    {it_i\mapsto it_i-p(-z)^i}  \\
     {is_i\mapsto is_i+q(-z)^i}
     \end{array}
  &
   I_3=\displaystyle{\int_{U(n)}}
  \det(I\!+\!zM)^p\det(I\!+\!z\bar M)^q e^{\Tr\VR(M)} dM

\end{array}
$$
$$\mbox{table 1}$$

 Clearly, applying the shifts $1,~2$ or $3$ to
$\BV_k$, as in table 1 to (\ref{BVirasoroNotation}),
lead to the Virasoro constraints
 $$
 \left.\BV_k\right|_{\mbox{\tiny shifted}}(I_i)=0
 ~~~~\mbox{for}~ k=-1,0,1$$ for the corresponding matrix integrals
 above.
 %, with the Virasoro operator $\left.\BV_k\right|_{\mbox{\tiny
 %shifted}}$ depending on $z$; the shifts 1,2,or 3 correspond
 %to the matrix integrals $I_1,~I_2$ and $I_3$.
   Consider the generators $\tilde\BV_z$ in the span of
$\BV_{-1},\BV_{0},\BV_{1}|_{\mbox{\tiny shifted}}$,
involving finite sums of $V_k(t),V_k(s), \pl /\pl t_k,
\pl /\pl s_k,\\ t_k, s_k$. They are given by the
$\tilde\BV_z$ in Table 2:

$$
\begin{array}{c|c|c|c}
 \mbox{shifts}&\tilde \BV_z(t,s) &\tilde\BV_{\Lambda} (t,s)
                 & \tilde L_{\Lb} \\
\hline
  & & &\\

  1:&\BV_0,~\BV_{\pm 1}\Bigr|_{{it_i\longmapsto it_i+z
     \delta_{i1}}\atop
    {is_i\longmapsto is_i-z\delta_{i1}}}
  &
  \tilde \BV_z \Bigr|_{z \mapsto k\Lb^{-1}_k}
   &
      \tilde L^{(1)}_0,~
      \tilde L^{(1)}_{\pm}

   \\  & & & \\

  2:& \pm(\BV_0+\BV_{\pm 1})
    \Bigr|_{{it_i\longmapsto it_i-q(-1)^i}
       \atop \!\!\!\!\!{is_i\longmapsto is_i- z\dt_{i1}}}
      &
      \tilde \BV_z \Bigr|_{z \mapsto k\Lb^{-1}_k}
      &
      \tilde L^{(2)}_{\pm}
      \\ & & &\\

   3:  & \left(\mathbb V_{-1}+(z+z^{-1})\mathbb
V_{0}+\mathbb V_{1} \right)
 \Bigr|_{{it_i\longmapsto it_i-p(-z)^i}  \atop
                        {is_i\longmapsto is_i+q(-z)^i}}
     &
     \tilde \BV_z \Bigr|_{z \mapsto \Lb^{-1}_k}
     &
   \tilde L^{(3)}

\end{array}
$$

$$\mbox{table 2}$$
%\vspace{1cm}
%
\noindent It is also useful to replace $z$ in the column
$\tilde\BV_z$ of Table 2 by -roughly speaking- the
operator $\Lb_k^{-1}$, where \be \Lb_k^{-1}
f(k)=f(k-1),~~k\in \BZ .\ee

The operators $\tilde\BV_{\Lb}(t,s)$ act on functions of
$(t,s,k) \in \BR^{\iy}\times \BR^{\iy} \times \BZ$. Let
them act on {\em two-dimensional Fourier series in Schur
polynomials} $\gs_{\lb}(t)$ and $\gs_{\mu}(s)$ with
partitions $\lb$ and $\mu$ having first column bounded
by $n$ and with arbitrary coefficients $\tilde
b^{(k)}_{\lb\mu}$. It turns out that the functions
obtained can be expressed again in terms of double
Fourier series, with coefficients $\tilde L (\tilde
b^{(k)}_{\lb\mu})$, which are linear difference
operators of finite order on the coefficients $\tilde
b^{(k)}_{\lb\mu}$:
\be
  \fbox{~~~$ \tilde\BV_{\Lb}(t,s)
   \displaystyle{
   \sum_{\lb,\mu\atop{\lb^{\top}_1,\mu^{\top}_1\leq
n}}}\tilde b^{(k)}_{\lb\mu}  \gs_{\lb}(t)\gs_{\mu}(-s) =
\displaystyle{
  \sum_{\lb,\mu\atop{\lb^{\top}_1,\mu^{\top}_1\leq n}}}
  \tilde L_{\Lb} (\tilde b^{(k)}_{\lb\mu})~
\gs_{\lb}(t)\gs_{\mu}(-s)$~~~}
\label{box1}\ee
 The {\em surprise} is that the
expressions $\tilde L_{\Lb}$, appearing in the right
hand column of table 2 and given explicitly in section
6,
%thus obtained, as spelled out in section 6,
 are precisely the difference equations for the transition
probabilities for certain {\em random walks, naturally
generated by the matrix integrals} above. This circle of
ideas will now be explained.

As a first ingredient, the three matrix integrals
containing Schur\footnote{
                           Throughout the paper, the
                           $\gs_{\lb}(t)$'s denote Schur
                           polynomials for a
                           partition $\lb$, expressed in
                           terms of the %less customary
                           symmetric
                           functions $kt_k=\sum_i x_i^k$,
                           and not in terms of the $x_k$
                           themselves.
                           The elementary
                  Schur polynomials ${\bf s}_k(t)$ are
                  defined by
                  $e^{\sum^{\infty}_1 t_i z^i }=
                  \sum^{\infty}_0 {\bf
                    s}_k(t)z^k$ and ${\bf s}_k(t)=0$
                    for $k<0$. Given a unitary matrix $M$,
                    we shall also use the notation
                    $\gs_{\lb}(M)$ to denote a
                           symmetric function of the
                           eigenvalues $x_1,\ldots,x_n$
                           of the unitary matrix $M$ and
                           thus in the notation of the
                           present paper
 $$
  {\bf s}_{\lb}(M):={\bf s}_{\lb}\bigl(\Tr M,\frac{1}{2}\Tr
M^2,\frac{1}{3}\Tr M^3,...\bigr).$$ }
polynomials $\gs_{\lb}(t)$ admit the following
expansions in $z$:
$$
\int_{U(n)}\gs_{\lb}(M)
            \gs_{\mu}(\bar M)
            \left\{
              \begin{array}{c}
               e^{z\Tr(M+\bar M)}\\  \\
              \det(I+M)^q e^{z\Tr \bar M} \\  \\
               \det(I+zM)^p\det(I+z\bar M)^q
               \end{array}
              \right\}
            dM
      =\sum_0^{\iy}\frac{z^k}{k!} \tilde b_{\lb\mu}^{(k)}
 $$
whereas the following matrix integrals below admit
(double) Fourier expansions in Schur polynomials:
 \bean
\lefteqn{
   \int_{U(n)}
            \left\{
              \begin{array}{c}
               (\Tr(M+\bar M))^k\\  \\
               (\Tr\bar M)^k \det(I+ M)^q \\  \\
              k! \gs_k (\ldots,
                -\frac{1}{i}
                 \Tr (p(-M)^i+q(-\bar M)^i),\ldots)
               \end{array}
              \right\}
          e^{\Tr \VR(M,\bar M)}  dM
   }\\  \\  \\
      &&\hspace{7cm}=
       \sum_{\lb,\mu\atop{\lb^{\top}_1,\mu^{\top}_1\leq n}}
  \tilde b_{\lb \mu}^{(k)}
            \gs_{\lb}(t)
            \gs_{\mu}(-s),
 \eean
whose coefficients are given by the same $\tilde
b^{(k)}_{\lb\mu}$. This is to say that the $\tilde
b^{(k)}_{\lb\mu}$'s appear in expansions of two
different integrals.

It is useful here to make the change of variables, from
partitions $\lb$ and $\mu$,
  with first columns smaller than $n$,
  %$\lb_1^{\top}\leq n$
  %and $\mu_1^{\top}\leq n$,
  to strictly increasing sets of integers
$x$ (initial position of the random walk) and $y $
(final position of the random walk) in $\BZ_{\geq 0}$,
defined by
 \bea
  x&:=& (x_1<x_2<\ldots<x_n)=
 (0+\lb_n,1+\lb_{n-1},\ldots , n-1+\lb_1)
 %\mbox{(initial position)}
 \no\\
 y&:=& (y_1<y_2<\ldots<y_n)=
 (0+\mu_n,1+\mu_{n-1},\ldots ,n-1+ \mu_1)
 %\mbox{(final position)}\no\\
 \no\\ \label{def(x-lb)}
   \eea
    So, the partitions {\sl
$\lb$ and $\mu$ measure the discrepancy from close
packing} $0,1,\ldots,n-1$ for $x$ and $y$! So the new
expressions
 $$
 b_{xy}^{(k)}=\tilde b^{(k)}_{\lb\mu}
  \Bigr|_{{x= (0+\lb_n,1+\lb_{n-1},\ldots , n-1+\lb_1)}
       \atop{y=(0+\mu_n,1+\mu_{n-1},\ldots ,n-1+ \mu_1)}}
~~\mbox{and}~~
 L_{\Lb}:=\tilde L_{\Lb}
  \Bigr|_{{x= (0+\lb_n,1+\lb_{n-1},\ldots , n-1+\lb_1)}
       \atop{y=(0+\mu_n,1+\mu_{n-1},\ldots ,n-1+ \mu_1)}}
 $$ satisfy the
 difference relations
 \be
  \fbox{ ~~$L_{\Lb}~ ( b_{xy}^{(k)})=0$~~}
 \label{L-relation}
 \ee
where the $b_{xy}^{(k)}$ have the following
interpretation in terms of walks\footnote{The number of
``effective moves"
   counts the actual steps taken by all walkers; i.e.,
   two walkers walking simultaneously counts for two moves,
   a walker not walking contributes nothing!}:
 \bea
 \lefteqn{\mbox{\bf Case 1}}\no\\
 b_{xy}^{(k)}&=& \#\left\{\begin{array}{l}
\mbox{ways that $n$ non-intersecting walkers
  in $\mathbb Z$ move during}\\
\mbox{$k$ instants from $x_1<x_2< ... <x_n$ to $y_1<y_2<
... <y_n$,}\\
 \mbox{where at each instant exactly one walker moves}\\
 \mbox{either one step to the left, or
one step to the right}
    \\
    \mbox{leading to $k$ effective moves  .}
\end{array}\right\}
     \no\\  \label{case 1}\\ %\no\\
    % \newpage
 \lefteqn{\mbox{\bf Case 2}}\no\\
 b_{xy}^{(k)}&=&
  \# \left\{
\begin{array}{l}
  \mbox{ways that n non-intersecting walkers move during $q+k$}\\
  \mbox{instants from $x_1< ... <x_n$ to $y_1< ... <y_n$,
   where at}\\
  \mbox{the instants $1$ to $q$, walkers
  \underline{may} move one step to the}\\
  \mbox{right, or \underline{stay put},
   and at the instants $q+1,...,q+k$
   }\\
  \mbox{\underline{exactly one} walker moves one step to the left, with}\\
  \mbox{total $\#\{\mbox{effective moves}\}=2k+\sum_1^n(y_i-x_i)
  $.}
  \end{array} \right\}
       \no\\  \label{case 2}\\ %\no\\
 \lefteqn{\mbox{\bf Case 3}}\no\\
 b_{xy}^{(k)}&=&  k!
  \# \left\{
\begin{array}{l}
  \mbox{ways that n non-intersecting walkers move during $p+q$}\\
  \mbox{instants from $x_1< ... <x_n$ to $y_1< ... <y_n$,
   where at}\\
  \mbox{the instants $1$ to $p$, walkers
  \underline{may} move one step to the}\\
  \mbox{right, or \underline{stay put},
   and at the instants $p+1,...,p+q$
   }\\
  \mbox{walkers
  \underline{may} move one step to the left or \underline{stay put},}\\
  \mbox{with total $\#\{\mbox{effective moves}\}=k
  $.}
  \end{array} \right\}\no\\\label{case 3}
\eea

\vspace{.0cm}

%\newpage

Case 1 leads, in particular, to a {\bf forward and
backward equation} for the transition probability, as
shown in section 7,

 \bean
 P(k,x,y)&:=&P\left( \begin{array}{l}
\mbox{ways that $n$ non-intersecting walkers
  in $\mathbb Z$ move during}\\
\mbox{$k$ instants from $x_1<x_2< ... <x_n$ to $y_1<y_2<
... <y_n$,}\\
 \mbox{where at each instant exactly one walker moves}\\
 \mbox{either one step to the left, or
one step to the right}
\end{array}
   \right)%=\frac{b^{(k)}_{xy}}{(2n)^{k}}
   \\
  &=&
 \frac{b^{(k)}_{xy}}{(2n)^k},
   \eean
namely,%satisfies the {\bf backward and forward} equations
  \be
 {\cal A}_i
 %\tilde L^{(1)}_{\pm}
 P(k,x,y)=0,
  \label{1.0.11}\ee
  where the %$\tilde L^{(1)}_{\pm}$
 ${\cal A}_i$
  are the same difference operators
  as the operator $L_{\Lb}$ in (\ref{L-relation}), except for the
  division by $2n$, which accounts for considering
  the probability rather than the $b^{(k)}_{xy}$'s\footnote{Define, for $\al\in\BZ$, $\al\neq 0$, the following
difference operators, acting on functions $f(k,x,y),$
with $ k\in \BZ_+, ~x,y\in \BZ$: \bea
\pl^+_{\al x_i}f &:=&f( k, x+\al e_i,y)-f(k,  x,y)  \no\\
\pl^-_{\al x_i}f &:=&f(  k,x,y)-f( k, x-\al e_i,y) \no\\
\Lb^{-1}_k f&:=& f(k-1,x,y) \label{definition}\eea}:
\bea
  {\cal A}_1
  %\tilde L^{(1)}_{+}
  &:=&
    \sum_{i=1}^n
   \left( \frac{k}{2n}\Lb^{-1}_k
     \pl^+_{2y_i}
     +
      x_i\pl^-_{x_i}
      +\pl^+_{y_i} y_i
   -
      (x_i-{y_i})
    \right)\no\\
  {\cal A}_2%\tilde L^{(1)}_{-}
  &:=&
    \sum_{i=1}^n
   \left( \frac{k}{2n}\Lb^{-1}_k
     \pl^+_{2x_i}
     +
      y_i\pl^-_{y_i}
      + \pl^+_{x_i} x_i
   -
      (y_i-{x_i})
    \right)\label{A-operator}
\eea
The operators ${\cal A}_1$ and ${\cal A}_2$ are the
\emph{forward and backward random walk equation},
because ${\cal A}_1$ essentially involves the end points
$y$, whereas ${\cal A}_2$ involves the initial points
$x$.

\remark The transition probabilities for a random walk
in $\BZ^n$ absorbed at the boundary of the Weyl chamber
$z_1<z_2< ... <z_n$, with equally likely steps $\pm
e_1,\ldots, \pm e_n$ (studied, e.g. in
(\cite{Gessel2},\cite{Grabiner1})) also satisfy the same
backward and forward difference equations
(\ref{1.0.11}).

%\end{theorem}

In a subsequent paper, we show the following limit
theorem: Let the spacings between the $n$ walkers and
the number of steps $k$ grow larger, with an appropriate
rescaling; then one finds, for fixed $n$, setting
$\vr:=\sqrt{\frac{n}{k}}$ and letting $k \rightarrow
\iy$ or, what is the same, letting $\vr \rightarrow 0$,
 \bean
   {\cal A}_1
  \Bigr|_{x= \frac{\tx}{\varepsilon\sqrt{t}},
          y= \frac{\ty}{\varepsilon\sqrt{t}}}
%    \\&&\\
%
  &=&
 {\frac{t^{(n+1)/2}}{\varepsilon}} {}
   \sum_1^n \Bigl(\frac{\pl}{\pl \ty_i}
              -\frac{\tx_i-\ty_i}{t}\Bigr)
  \\
  &&~~~~~~
     -2\hspace{.07cm}t^{(n+2)/2} \Bigl(
     \frac{\pl}{\pl t} - \frac{1}{2}
   \sum_1^n \frac{\pl^2}{\pl \ty_i^2}
          \Bigr)%h
          +O\left(\vr\right)
          \\
  {\cal A}_2
  \Bigr|_{x= \frac{\tx}{\varepsilon\sqrt{t}},
          y= \frac{\ty}{\varepsilon\sqrt{t}}}
 &=&~\mbox{same expression, but with}~ \tx\leftrightarrow \ty
  \eean
This expansion in $\vr$ is valid, when acting on an
appropriate function space. The term $O(1)$ in the
$\vr$-expansion of ${\cal A}_1$ contains precisely the
forward diffusion equation for Brownian motion in the
Weyl chamber $\{\ty_1<\ldots<\ty_n\}\subset \BR^n$. The
term $O(1)$ in the $\vr$-expansion of ${\cal A}_2$
contains the corresponding backward equation.  This
Brownian motion in the Weyl chamber is tantamount to the
motion of $n$ non-intersecting Brownian motions, i.e.,
who are killed as soon as they collide. It is related to
Dyson's Brownian motion (\cite{Dyson}).
%\newpage

\bigbreak

  The difference equations (\ref{L-relation}) above for the three cases
  (\ref{case 1}),(\ref{case 2}) and (\ref{case 3})
  are based on replacing the operators
  ${\pl}/{\pl t_n}$ and multiplication by $nt_n$ in the
  Murnaghan-Nakayama
   rule,
  %
%  In this section we study the action of the
%  Virasoro operators on Schur polynomials $ {\bf s}_{\lb}$.
%   It useful to remind the reader of the Murnaghan-Nakayama
%   rule, stated later, from which it follows that
%
\bea
nt_n~\gs_{\lb}(t)&=&\sum_{{\mu}\atop{\mu\backslash\lb\in
B(n)}}(-1)^{\acht(\mu\backslash\lb)}\gs_{\mu}(t) \no\\
\no\\
\frac{\pl}{\pl
t_n}\gs_{\lb}(t)&=&\sum_{{\mu}\atop{\lb\backslash\mu\in
B(n)}}(-1)^{\acht(\lb\backslash\mu)}\gs_{\mu}(t),
 \label{1.0.14}
\eea
by the action of the Virasoro algebra:
 \bea
  V_{-n} {\bf s}_{\lb}&=&\sum_{{\mu}\atop{\mu\backslash\lb\in
B(n)}} d^{(-n)}_{\lb \mu} {\bf s}_{\mu}
 \label{Virasoro10} \no\\
V_{n} {\bf
s}_{\lb}&=&\sum_{{\mu}\atop{\lb\backslash\mu\in B(n)}}
d^{(-n)}_{\mu\lb } {\bf s}_{\mu}
   \label{Virasoro20}   \eea
with $(n\geq 1)$
  \bea
 d^{(-n)}_{\lb \mu}
  &=& \sum_{i\geq 1}~~
    \sum_{
    \left\{
    \begin{array}{l}
           \nu ~~~~~\mbox{\footnotesize such that} \\
           \lb \backslash \nu \in B(i)\\
           \mu \backslash \nu \in B(n+i) \\
            \lb \backslash \nu  \subset \mu \backslash \nu
            \end{array}
            \right\}}
 (-1)^{\mbox{\footnotesize ht}(\lb \backslash \nu) +
       \mbox{\footnotesize ht}(\mu \backslash \nu)  }\no \\
   &&\no  \\&&
   + \frac{1}{2} \sum_{i=1}^n
   ~~
   \sum_{
   \left\{
   \begin{array}{l}
           \nu ~~~~~\mbox{\footnotesize such that} \\
           \nu \backslash \lb \in B(i)\\
           \mu \backslash \nu \in B(n-i) \\
            %\lb \backslash \nu  \subset \mu \backslash \nu
            \end{array}
         \right\}   }
 (-1)^{\mbox{\footnotesize ht}(\nu \backslash \lb) +
       \mbox{\footnotesize ht}(\mu \backslash \nu)  } .
       \label{coefficients}\eea

This will be shown in section 5. To explain the
notation, $h\in B(i)$ denotes a border-strip (i.e., a
connected skew-shape $\lb \ba \mu$ containing $i$ boxes,
with no $2\times 2$ square) and the height $\acht~h$ of
a border strip $h$ is defined as
 \be
 \acht~h:=\#\{\mbox{rows
in~}h\}-1 .\ee
 In view of the infinite sum in the Virasoro algebra
  (\ref{VirasoroIntro}), one would expect $V_{n} {\bf
s}_{\lb}$ to be expressible as
  an {\em infinite sum of Schur polynomials}.
  This is not so: acting with Virasoro $V_{-n}$ (resp.
  $V_n$)
  leads to the same precise sum as acting with $nt_n$ (resp.
  $\pl/ \pl t_n$), except the coefficients
  (\ref{coefficients}) are
  different from the ones in (\ref{1.0.14}). This is to
  say {\em the two operators have the same band structure or
  locality}!

%\newpage

\section{Non-intersecting walks and partitions}

Consider $m$ %arbitrary non-intersecting
 %lock-step
 walkers on $\BZ$ departing from position $x_1,\ldots,x_m$,
 and ending up at $y_1,\ldots, y_m$, such that
 at each instant, only one walker moves either one step to the left,
or one step to the right, with all possible moves
equally likely. The main statement of this section is to
connect the transition probabilities of these walks with
pairs of skew Young tableaux and other combinatorial
formulas. They will play a crucial role in section 4.
Such connections have been known in various situations
in the combinatorial literature; see R. Stanley
\cite{Stanley2} (p. 313), P. Forrester \cite{Forrester},
D. Grabiner \& P. Magyar \cite{Grabiner1, Grabiner2}, J.
Baik \cite{Baik}:

\begin{theorem}

\bean \lefteqn{
 P\left( \begin{array}{l}
 \mbox{that $m$ walkers in $\BZ$,}\\
 \mbox{go from $x_1,\ldots,x_m$ to}\\
 \mbox{$y_1,\ldots,y_m$ in $T$ steps,}\\
 \mbox{and do not intersect}
 \end{array}
   \right)
   } \\
   &=&
   \frac{1}{(2m)^{T}}
   \left(  {T}
    \atop{T_L~~ T_R }
    \right)
   \sum_{
  {{\lb~\mbox{\tiny with}~\lb \supset \mu,\nu}\atop
 {   { {|\lb\ba \mu|=T_L}\atop{|\lb\ba \nu|=T_R}}
 \atop {\lb_1^{\top} \leq m}
   }    }
   }
     f^{\lb\ba \mu}  f^{\lb\ba \nu}
 \\&&\\
 &=& \frac{1}{(2m)^T}
  \sum_{w\in W} (-1)^{\sigma(w)}
    \left.  \left( \sum_1^m (u_i+u_i^{-1})\right)^T
   \right|_{u_1^{y_1-w(x_1)}\ldots u_m^{y_m-w(x_m)}}
   \\ &&\\ &&\\
 &=&\frac{1}{(2m)^T}\sum_{w\in W} (-1)^{\sigma (w)}\\
&&\hspace{-1cm}
  \renewcommand{\arraystretch}{0.5}
\begin{array}[t]{c}
\sum\\
{\scriptstyle k_i\geq 0}\\
{\scriptstyle \sum_1^m k_i=T/2}
\end{array}
\renewcommand{\arraystretch}{1}
\left( \begin{array}{c}
T\\
k_1+\frac{y_1-x_{w(1)}}{2} \,,\,\,
k_1-\frac{y_1-x_{w(1)}}{2} \,, \ldots \,,
k_m+\frac{y_m-x_{w(m)}}{2} \,,\,\,
k_m+\frac{y_m-x_{w(m)}}{2}
\end{array}
\right) \eean
%\varsupsetneq
%
 where $\mu$, $\nu$ are fixed partitions defined by the points
 $x_i$ and $y_i$,
 \bean \mu_k &=& k-1-x_k,~~~
        \nu_k = k-1-y_k  \\
 T_L&=&\frac{1}{2} (T+\sum_1^m(x_i-y_i))=
   \frac{1}{2} (T-|\mu|+|\nu|)\\
      T_R&=&\frac{1}{2} (T-\sum_1^m(x_i-y_i))=
   \frac{1}{2} (T+|\mu|-|\nu|)\\
   T&=&T_L+T_R,~~ \sum_1^m (x_i-y_i)=T_L-T_R.
 \eean

 \end{theorem}

\proof will follow from Propositions 2.3 and 2.5, as
given in the subsequent subsections.

\subsection{Non-intersecting walks and skew-tableaux}

\begin{proposition} There is  a $1-1$ correspondence between
$$
\left\{
\begin{array}{l}
\mbox{ways that $m$ non-intersecting} \\
 \mbox{walkers move
 from $x_1 <\ldots < x_m$\!\!\!}
 \\
 \mbox{to $y_1<\ldots <y_m$, where at}\\
\mbox{each instant $1,\dots,T_L$ one} \\
\mbox{walker moves one step to}\\
\mbox{the left, and at each instant}\\
\mbox{$T_L+1,\dots,T_L+T_R$ one walker} \\
 \mbox{moves one step to the right}
 %\\
%\mbox{}
\end{array}
\right\}\Longleftrightarrow\left\{
\begin{array}{l}
\mbox{all couples $(P,Q)$ of standard}\\
\mbox{skew-tableaux of arbitrary}\\
\mbox{shape $\lb\ba\mu$ and $\lb\ba\nu$,}\\
\mbox{given fixed partitions $\mu$, $\nu$,}\\
\mbox{with $\vert\lb\ba\mu\vert = T_L$,
           $\vert\lb\ba\nu\vert = T_R$,}\\
\mbox{filled with numbers}\\
\mbox{$1,\ldots\!,T_L$ and $1,\ldots\!,T_R$}\\
\mbox{with $\lb^{\top}_1\leq m$}
\end{array}
\right\},
$$
 and so, for this walk,
$$
\# \left\{
\begin{array}{l}
\mbox{ways that $m$ non-intersecting}\\
\mbox{walkers go $T_L$ steps to the left}\\
\mbox{and then $T_R$ steps to right,}\\
\mbox{from $x_1 <\ldots < x_m$ }\\
 \mbox{to $y_1<\ldots <y_m$}
\end{array}
\right\}
 =
 \renewcommand{\arraystretch}{0.5}
\begin{array}[t]{c}
{\displaystyle \sum} \\
{\scriptstyle \lb \vdash \frac{1}{2} (T+|\mu|+|\nu|)      }\\
{\scriptstyle \lb \supset \mu,\nu       }\\
{\scriptstyle \lb^{\top}_1\leq m}
\end{array}
\renewcommand{\arraystretch}{1}
f^{\lb \ba \mu}f^{\lb \ba \nu},
 $$
where $\mu$, $\nu$ are fixed partitions defined by the
points
 $x_i$ and $y_i$, and $\lb \supset \mu, \nu$, such that
 \bean \mu_k &:=& k-1-x_k\\
        \nu_k &:=& k-1-y_k  \\
       |\lb|&=& \frac{1}{2} (T+|\mu|+|\nu|),~~
       T=T_L+T_R.
%        \\
% T_L&=&\frac{1}{2} (T+\sum_1^m(x_i-y_i))=
%   \frac{1}{2} (T-|\mu|+|\nu|)\\
%      T_R&=&\frac{1}{2} (T-\sum_1^m(x_i-y_i))=
%   \frac{1}{2} (T+|\mu|-|\nu|)
 \eean

\end{proposition}

\proof Consider two Young diagrams $\mu,\nu$ and $\lb$
such that

$$
\mu\subset\lb,\quad\nu\subset\lb
$$
and two standard skew-tableaux

$$
(P,Q)=\left\{\begin{array}{l}
\mbox{$P$ standard skew-Young tableaux of shape $=\lb\ba\mu$}\\
\mbox{filled with numbers $1,\ldots,|\lb\ba\mu|$}\\
\\
\mbox{$Q$ standard skew-Young tableaux of shape $=\lb\ba\nu$}\\
\mbox{filled with numbers $1,\ldots,|\lb\ba\nu|$.}
\end{array}
\right.
$$

%\newpage
 $$
\raisebox{12mm}{\begin{tabular}{l}
\raisebox{3mm}{1$^{\rm st}$ walker}\\
\raisebox{3mm}{2$^{\rm nd}$ walker}\\
\raisebox{3mm}{3$^{\rm rd}$ walker}\\
\raisebox{3mm}{4$^{\rm th}$ walker}
\end{tabular}}\overbrace{\squaresize .7cm \thickness .01cm
\Thickness .07cm \Young{ \shade&\shade &\shade &\shade
&c_{11}&c_{12}\cr \shade&\shade &c_{21}&c_{22}\cr
\shade&c_{31}&c_{32}&c_{33}\cr \shade&c_{41}\cr
}}^{\mbox{instants of left move}}\quad
\raisebox{12mm}{\begin{tabular}{l}
\raisebox{3mm}{ }\\
\raisebox{3mm}{ }\\
\raisebox{3mm}{ }\\
\raisebox{3mm}{ }
\end{tabular}} \overbrace{\squaresize .7cm \thickness .01cm
\Thickness .07cm \Young{ \shade&\shade &\shade
&c'_{11}&c'_{12}&c'_{13}\cr \shade&
c'_{21}&c'_{22}&c'_{23}\cr
\shade&c'_{31}&c'_{32}&c'_{33}\cr c'_{41}&c'_{42}\cr }
}^{\mbox{instants of right move}}
 $$
$$\hspace*{17mm}P~\mbox{of shape $\lb\ba\mu$}\hspace*{30mm}
                Q~\mbox{of shape $\lb \ba \nu$}
$$

$$\mbox{figure 1}$$

\bigbreak

%\newpage
\noindent To the standard   skew-Young tableaux $P$ of
shape $\lb\ba \mu$, we associate
%a walk in the following way. For $\lb=(\lb_1\geq
%\ldots\geq\lb_m
%>0)$, the
 $m$ walkers starting at
$$
x_1=-\mu_1+0<%x_2=-\mu_2+1<
 ... <x_k=-\mu_k+k-1<... <x_m=-\mu_m+m-1,
$$
so that $x_k-x_{k-1}=\mu_{k-1}-\mu_k+1$ and requiring
the $k^{\rm th}$ walker, starting
%$(1\leq k\leq m)$ starts
 at $x_k=-\mu_k+k-1$, to move to the left
only, at instants
$$
c_{ki}=\mbox{content of box $(k,i)\in P$},
$$
and thus he has made, in the end, $\lb_k-\mu_k$ steps to
the left. So, at each instant exactly one walker is
moving and this during a time-span $T_L=|\lb\ba\mu|$,
until the $m$ walkers reach the position

$$
-\lb_1+0<... <-\lb_k+k-1 <... <-\lb_m+m-1.
$$
The fact that the skew-tableau $P$ is standard implies
that the walkers have never intersected, as one sees by
imagining $\mu$ filled in a standard fashion with the
numbers $0,-1,\ldots,-|\mu|+1$, thus yielding for each
walker $k$, a path from $x_k=-\mu_k+k-1$ to
$-\lb_k+k-1$, not intersecting the paths of the
neighboring walkers.

In the same way, to $Q$, we associate $m$ walkers
starting at

$$
y_1=-\nu_1+0<
%y_2=-\nu_2+1<
 ...<y_k=-\nu_k+k-1<...<y_m=-\nu_m+m-1
$$
moving left at instants $c'_{ki}$, each making, in the
end $\lb_k-\nu_k$ steps to the left, until the $m$
walkers reach, after time $T_R=|\lb\ba\nu|$ the position

$$
-\lb_1+0<...<-\lb_k+k-1<...<-\lb_m+m-1,
$$
 the same position as before,
 without having ever intersected.

 Now we assemble the two walks. It yields a walk with $m$ non-intersecting walkers,
going from $(x_1<\ldots<x_m)$ to
 $(y_1<\ldots<y_m)$, moving first $T_L$ steps to the left
  and then
 $T_R$ steps to the right, obtained by {\em reversing} the second walk
 (associated with $Q$) such that the final position of the first
 walk (moving left) is the starting position of the second walk
 (moving right). Therefore, we have the total number of
 steps
 \be T=T_L+T_R=2|\lb|-|\mu| - |\nu| , \ee
 from which Proposition 2.2 follows.\qed

%\newpage

%Since $x_k=-\mu_k+k-1$ and $y_k=-\nu_k+k-1$, we have

%$$
%\sum_1^mx_k=-|\mu|-\frac{m(m-1)}{2}\mbox{~~and~~}\sum_1^my_k=-|\nu|-\frac{m(m-1)
%}{2}.
%$$
%From the construction,

%\bean
%T_L&:=&\#\{\mbox{moves to the left}\}=|\lb|-|\mu|\\
%T_R&:=&\#\{\mbox{moves to the right}\}=|\lb|-|\nu| \eean so that

%Hence, using the equations
%above
%
%$$
%T_{\left\{
%\begin{array}{c}
%L\\
%R
%\end{array}\right\}}=|\lb|-|\left\{
%\begin{array}{c}
%\mu\\
%\nu
%\end{array}|\right\}=\frac{1}{2}\left(T\pm\sum^m_1(x_k-y_k)\right).
%$$

\vspace{1cm}

%\newpage

Consider now $m$ %arbitrary non-intersecting
 %lock-step
 walkers on $\BZ$ departing from position $x_1,\ldots,x_m$,
 and ending up at $y_1,\ldots, y_m$, as in Theorem 2.1.

\begin{proposition}

\be %\lefteqn{
 P\left( \begin{array}{l}
 \mbox{that $m$ walkers in $\BZ$,}\\
 \mbox{go from $x_1,\ldots,x_m$ to}\\
 \mbox{$y_1,\ldots,y_m$ in $T$ steps,}\\
 \mbox{and do not intersect}
 \end{array}
   \right)
   %}
   =
   \frac{1}{(2m)^{T}}
   \left(  {T}
    \atop{T_L~~ T_R }
    \right)
   \sum_{
  {{\lb~\mbox{\tiny with}~\lb \supset \mu,\nu}\atop
 {   { {|\lb\ba \mu|=T_L}\atop{|\lb\ba \nu|=T_R}}
 \atop {\lb_1^{\top} \leq m}
   }    }
   }
     f^{\lb\ba \mu}  f^{\lb\ba \nu}
\label{2.1.2}\ee
%\varsupsetneq
%
 where $\mu$, $\nu$ are fixed partitions defined by the
 $x_i$ and $y_i$'s, as in Theorem 2.1.
 %\bean \mu_k &=& k-1-x_k\\
 %       \nu_k &=& k-1-y_k  \\
 %T_L&=&\frac{1}{2} (T+\sum_1^m(x_i-y_i))=
 %  \frac{1}{2} (T-|\mu|+|\nu|)\\
 %     T_R&=&\frac{1}{2} (T-\sum_1^m(x_i-y_i))=
 %  \frac{1}{2} (T+|\mu|-|\nu|),~~~T=T_L+T_R.
 %\eean
\end{proposition}

\vspace{.0cm}

\begin{corollary}
\be
 P\left( \begin{array}{l}
 \mbox{that $m$ walkers in $\BZ$,}\\
 \mbox{go from and return}\\
 \mbox{to $1,\ldots,m$ in $2n$ steps,}\\
 \mbox{and do not intersect}
 \end{array}
   \right)
   =\frac{1}{(2m)^{2n}} \left(  {2n}\atop{n}\right)
   \sum_{{\lb\vdash n}\atop{\lb_1^{\top} \leq m}} \left( f^{\lb}\right)^2
\ee
\end{corollary}

{\medskip\noindent{\it Proof of Proposition 2.3:\/} } At
first, notice that $m$ walkers in $\BZ$, obeying this
rule, is tantamount to a random walk in $\BZ^m$, where
at each point the only moves are
 $$
 \pm e_1,\ldots, \pm e_m.
 $$
 That is to say the walk has at each point $2m$ possibilities and
 thus at time $T$ the walk has \be (2m)^{T}\ee places to go.

Returning to the $\BZ$-picture, we now associate to a
given walk
 a sequence of $T_L$ $L$'s and $T_R$ $R$'s:
 \be
 L~~R~~R~~{R~~L}~~R~~L~~L~~R~\ldots ~R~,
 \label{LR's}\ee
thus recording the nature of the move, left or right, at
the first instant, at the second instant, etc...

If the $k^{th}$ walker is to go from $x_k$ to $y_k$,
then
$$
y_k-x_k=\#\left\{ \begin{array}{c}
  \mbox{right moves }\\
  \mbox{for $k^{th}$ walker}
  \end{array} \right\}
  -
  \#\left\{ \begin{array}{c}
  \mbox{left moves }\\
  \mbox{for $k^{th}$ walker}
   \end{array} \right\}
  $$
  and so, if
  $$
  T_L:=\#\left\{ \begin{array}{c}
  \mbox{left moves }\\
  \mbox{for all $m$ walkers}
  \end{array} \right\}
~~~ \mbox{and} ~~~
  T_R:=\#\left\{ \begin{array}{c}
  \mbox{right moves }\\
  \mbox{for all $m$ walkers}
  \end{array} \right\},
 $$
 we have, since at each instant exactly one walker moves,
  \bean
  T_R+T_L&=& T  \\
  T_R-T_L&=& \sum_1^m (y_k-x_k),
  \eean
  from which
  $$
  T_{\left\{L\atop R\right\}}=\frac{1}{2} \left( T\pm \sum_1^m (x_k-y_k)\right)
  .$$
  Next, we show there is a {\em
canonical} way to map a walk, corresponding to
(\ref{LR's}) into one with left moves only during times
$1,\ldots,T_L$ and then right moves during times
$T_L+1,\ldots,T_L+T_R=T$, thus corresponding to a
sequence
 \be
 \overbrace{L~~L~~L~\ldots ~L}^{T_L}~~\overbrace{R~~R~~R~\ldots
 ~R}^{T_R}~.
 \ee
Indeed, in a typical sequence, as (\ref{LR's}),
 \be
 L~~R~~R~~\underline{R~~L}~~R~~L~~L~~R~\ldots ~R~,
 \ee
 consider the first sequence $R~L$ (underlined) you encounter,
 in reading from left to right. It corresponds to one of the following three
 configurations (in the left column),
%
% \newpage
%
\bean
 \begin{array}{l}
L\\
R \end{array}
 {\displaystyle {\diagdown}\atop{~~~
\big|}}~~{\displaystyle {\big|}\atop{ \big|}}
       ~~{\displaystyle {\big|}\atop{ \big|}}
       ~~{\displaystyle {~~~ \big|}\atop{ \diagup}}
~~~~~~~~~~
  &\Rightarrow&
~~~~~~~~~~
 \begin{array}{l}
R\\
L \end{array}
 {\displaystyle { \big|}\atop{~~ \diagdown}}
 ~~{\displaystyle {\big|}\atop{ \big|}}
       ~~{\displaystyle {\big|}\atop{ \big|}}
       ~
{\displaystyle {~~\diagup}\atop{ \big|}}\\
&&\\&&\\
%@@@@@@@@@@@@@@@@@@@@
 \begin{array}{l}
L\\
R \end{array}
  {\displaystyle {~~~ \big|}\atop{ \diagup}}
       ~~{\displaystyle {\big|}\atop{ \big|}}
       ~~{\displaystyle {\big|}\atop{ \big|}}
       ~~{\displaystyle{\diagdown}\atop{~~~ \big|}}
~~~~~~~~~~
  &\Rightarrow &
~~~~~~~~~~
 \begin{array}{l}
R\\
L \end{array}
 {\displaystyle {~~\diagup}\atop{ \big|}}
 ~~{\displaystyle {\big|}\atop{ \big|}}
       ~~{\displaystyle {\big|}\atop{ \big|}}
       ~{\displaystyle { \big|}\atop{~~ \diagdown}}
 \\
 &&\\
 &&
 \\
%@@@@@@@@@@@@@@@@@@@@@@@@@
  \begin{array}{l}
L\\
R \end{array}
  {\displaystyle {~~~ \big|}\atop{~~~ \big| }}
      ~~~ . ~ ~~{\displaystyle {\big|}\atop{ \big|}}
       ~~{\displaystyle {\big|}\atop{ \big|}}
       ~~{\displaystyle{\diagdown}\atop{\diagup}}
~~~~~~~~~~
  &\Rightarrow&
~~~~~~~~~~
 \begin{array}{l}
R\\
L \end{array}
 {\displaystyle {~~\big|}\atop{ ~~\big|}}
 ~~{\displaystyle {\diagup }\atop{ \diagdown }}
       ~~{\displaystyle {\big|}\atop{ \big|}}
       {\displaystyle { ~~\big|}\atop{~~ \big|}}
       ~~{\displaystyle {\big|}\atop{ \big|}}
 \eean
which then can be transformed into a new
 configuration $L~R$, with same beginning and end,
  thus yielding a new sequence;
 in the third case
 the reflection occurs the first place it can. These
 moves have been considered by Forrester in
 \cite{Forrester}.
 So, by the moves above, the original configuration
  (\ref{LR's}) can be transformed in a new one.
  In the new sequence, pick again the first sequence $RL$,
  reading from left to right, and use again one of the
  moves. So, this leads again to a new sequence, etc...
 \bean
 &L~~R~~R~~\underline{R~~L}~~R~~L~~L~~R~\ldots ~R&\\
 &L~~R~~\underline{R~~L}~~R~~R~~L~~L~~R~\ldots ~R& \\
 &L~~\underline{R~~L}~~R~~R~~R~~L~~L~~R~\ldots ~R&\\
 &L~~L~~R~~R~~R~~\underline{R~~L}~~L~~R~\ldots ~R&\\
 &L~~L~~R~~R~~\underline{R~~L}~~L~~L~~R~\ldots ~R&\\
 &\vdots&
 \eean
%
%
%
%There is a canonical way to transform this walk into a new
%non-intersecting lock-step walk going with the sequence
 \be
 \underbrace{L~~L~~L~\ldots ~L}_{T_L}~~
  \underbrace{R~~R~~R~\ldots ~R}_{T_R}~;
 \ee
  Since this procedure is invertible, it gives a {\em one-to-one} map between all the left-right walks
  corresponding to a given sequence, with $T_L$ $L$'s and $T_R$
  $R$'s
  \be
 L~~R~~R~~{R~~L}~~R~~L~~L~~R~\ldots ~R~;
 \ee
and all the walks corresponding to \be
 \overbrace{L~~L~~L~\ldots ~L}^{T_L}
 ~~\overbrace{R~~R~~R~\ldots ~R}^{T_R}~.
 \ee
 On the one hand, corresponding to this sequence,
 $m$ walkers can walk in
 \be
   \renewcommand{\arraystretch}{0.5}
\begin{array}[t]{c}
{\displaystyle \sum} \\
{\scriptstyle \lb \vdash \frac{1}{2} (T+|\mu|+|\nu|)      }\\
{\scriptstyle \lb^{\top}_1\leq m}
\end{array}
\renewcommand{\arraystretch}{1}
f^{\lb \ba \mu}f^{\lb \ba \nu}
 \label{last}\ee
%ft( f^{\lb}\right)^2
%
ways, as follows from Proposition 2.2. On the other
hand, there are $\left( {T}
    \atop{T_L~~ T_R }
    \right)$
    %$\left( {2n}\atop{n}\right)$
    sequences of $T_L$ $L$'s and $T_R$ $R$'s,
which combined with (\ref{last}) yields the result.\qed

{\medskip\noindent{\it Proof of Corollary 2.4:\/} } In
this situation, we have close packing, and thus
$\mu_k=\nu_k=0$ for all $k$, and so $\mu=\nu=\emptyset $
and $T_L=T_R=T/2$. With these data, (2.1.3) is an
immediate consequence of (2.1.2).   \qed

%\newpage

\subsection{Non-intersecting walks and D. Andr\'e's principle}

We now prove the second and third formulae of Theorem
2.1, namely

\begin{proposition}

\bean \lefteqn{
 P\left( \begin{array}{l}
 \mbox{that $m$ walkers in $\BZ$,}\\
 \mbox{go from $x_1,\ldots,x_m$ to}\\
 \mbox{$y_1,\ldots,y_m$ in $T$ steps,}\\
 \mbox{and do not intersect}
 \end{array}
   \right)
   }
   \\&&\\
 &=&\frac{1}{(2m)^T}\sum_{w\in W} (-1)^{\sigma (w)}
   \left.\left( \sum_1^m (u_i+u_i^{-1})\right)^T
\right|_{u_1^{y_1-x_1} \dots u_m^{y_m-x_m}}
\\
 &=&\frac{1}{(2m)^T}\sum_{w\in W} (-1)^{\sigma (w)}\\
&&\hspace{-1cm}
  \renewcommand{\arraystretch}{0.5}
\begin{array}[t]{c}
\sum\\
{\scriptstyle k_i\geq 0}\\
{\scriptstyle \sum_1^m k_i=T/2}
\end{array}
\renewcommand{\arraystretch}{1}
\left( \begin{array}{c}
T\\
k_1  \!\!  +  \!\!  \frac{y_1-x_{w(1)}}{2}, \,\,\,
 k_1 \!\!  -  \!\!  \frac{y_1-x_{w(1)}}{2}, \, \ldots \, ,
  k_m  \!\!  +   \!\!  \frac{y_n-x_{w(m)}}{2}, \,\,\,
k_m  \!\!  +  \!\!  \frac{y_n-x_{w(m)}}{2}
\end{array}
\right). \eean

\end{proposition}

%\newpage

 \proof Consider again the random walk in $\BZ^m$,
 where at each point the only moves are
 $$
 \pm e_1,\ldots, \pm e_m.
 $$
 Then
\begin{eqnarray*}
&&\#\left\{
\begin{array}{l}
\mbox{ways to walk in $T$ steps} \\
\mbox{from }x=(x_1,\ldots,x_m)\\
\mbox{to } y=(y_1, \ldots, y_m) \mbox{ in } \BZ^m
\end{array}
\right\}\\
&=&\left.\left( \sum_1^m (u_i+u_i^{-1})\right)^T
\right|_{u_1^{y_1-x_1} \dots
u_m^{y_m-x_m}} \\
&=&\left.\renewcommand{\arraystretch}{0.5}
\begin{array}[t]{c}
\sum\\
{\scriptstyle \sum_1^m (a_i+b_i)=T}\\
{\scriptstyle a_i, b_i \geq 0}
\end{array}
\renewcommand{\arraystretch}{1}
\left(\begin{array}{c} T \\a_1~~ b_1\ldots
a_m~~b_m\end{array}\right)u_1^{a_1} u_1^{-b_1} \ldots
u_m^{a_m}
u_m^{-b_m} \right|_{u_1^{y_1-x_1} \ldots u_m^{y_m-x_m}}\\
&=&\renewcommand{\arraystretch}{0.5}
\begin{array}[t]{c}
\sum\\
{\scriptstyle \sum_1^m (a_i+b_i)=T}\\
{\scriptstyle a_i - b_i =y_i-x_i}
\end{array}
\renewcommand{\arraystretch}{1}
\left(\begin{array}{c} T \\a_1 ~~ b_1\ldots
a_m ~~ b_m\end{array}\right)\\
&=&\renewcommand{\arraystretch}{0.5}
\begin{array}[t]{c}
\sum\\
{\scriptstyle \sum_1^m k_i=T/2}\\
{\scriptstyle k_i\geq 0}
\end{array}
\renewcommand{\arraystretch}{1}
\left(\begin{array}{c} T\\
k_1+\frac{y_1-x_1}{2} ~~ k_1-\frac{y_1-x_1}{2} \ldots
k_m+\frac{y_m-x_m}{2} ~~ k_m-\frac{y_m-x_m}{2}
\end{array}\right)
\end{eqnarray*}
%
%
%\newpage
Hence, by the D. Andr\'e reflection principle,
\begin{eqnarray*}
&&\#\left\{
\begin{array}{l}
\mbox{ways to walk in $\BZ^n$ in $T$ steps from $x$ to $y$} \\
%\mbox{} \\
\mbox{within the region } \{z_1 < \ldots < z_n\}\subset
\BZ^n
\end{array}
\right\}\\
&=&\sum_{w\in W} (-1)^{\sigma (w)}
 \# \left\{
\begin{array}{l}
\mbox{ways to walk in }\BZ^n \mbox{ in }  \\
T \mbox{ steps from } w(x)\longmapsto y
\end{array}
\right\}\\
&=&\sum_{w\in W} (-1)^{\sigma (w)}\\
&&\hspace{-1cm}
\renewcommand{\arraystretch}{0.5}
\begin{array}[t]{c}
\sum\\
{\scriptstyle k_i\geq 0}\\
{\scriptstyle \sum_i k_i=T/2}
\end{array}
\renewcommand{\arraystretch}{1}
\left( \begin{array}{c}
T\\
k_1+\frac{y_1-x_{w(1)}}{2}, \,\,\,
k_1-\frac{y_1-x_{w(1)}}{2}, \, \ldots \,,
k_n+\frac{y_n-x_{w(n)}}{2}, \,\,\,
k_n-\frac{y_n-x_{w(n)}}{2}
\end{array}
\right),
\end{eqnarray*}
ending the proof of Proposition 2.5. \qed

%\newpage

\section{A matrix integral}

The main statement of this section is to prove

\begin{proposition} The following matrix integral admits a
``Fourier" expansion
\bea
  \lefteqn{
   \int_{U(n)}e^{\sum_1^{\iy}\Tr(t_j^{(0)}M^j-s_j^{(0)}\bar
M^j)} e^{\sum_1^{\iy}\Tr(t_j M^j-s_j\bar M^j)}dM
  }
  \no\\
  &=&
  \sum_{\lb,\mu~\mbox{\tiny such
 that}\atop{\lb_1^{\top},\mu_1^{\top}\leq n}}
 a_{\lb\mu}(t^{(0)},s^{(0)})s_{\lb}(t)s_{\mu}(-s)
 \label{generalmatrixintegral}\eea
 with Fourier coefficients, taking on many different
 forms:
 \bea
  a_{\lb\mu}(t^{(0)},s^{(0)})
   &=&\det\left(\oint_{S^1}u^{\lb_{\ell}-\ell
-\mu_k+k}e^{\sum^{\iy}_1(t^{(0)}_ju^j-s^{(0)}_ju^{-j})}\frac{du}{2\pi~iu}\right)_{1\leq
\ell,k\leq n}  \no\\
   &=&
   \sum_{{\nu ~\mbox{\tiny with}~\nu \supset \lb,\mu
    \atop{\nu_1^{\top} \leq
n}}\atop{}}
 {\bf s}_{\nu\ba\lb}(t^{(0)})
  {\bf s}_{\nu\ba\mu}(-s^{(0)})
 \no\\
&=&\int_{U(n)} {\bf s}_{\lb}(M){\bf s}_{\mu}(\bar M)
e^{\sum^{\iy}_1\Tr(t^{(0)}_jM^j-s^{(0)}_j\bar M^j)}dM.
 \label{Prop3.1}
   \eea

%\newpage

    \end{proposition}

The proof of this Proposition will be given later in
this section. Specializing the $t^{(0)}$'s and
$s^{(0)}$'s will lead to several examples, discussed in
the next section (section 4).

 Consider the
$(t,s)$-dependent semi-infinite matrix \be
m_{\iy}(t,s)=\left( \mu_{ij}(t,s)\right)_{0\leq i,j <\iy
},\ee
 evolving according to the equations (here $\Lb$ is the
 shift matrix, with zeroes everywhere, except for $1$'s just above
 the diagonal, i.e., $(\Lambda v)_n=v_{n+1}$)

 \be
 \frac{\pl m_{\iy}}{\pl t_k}=\Lb^k m_{\iy}~~\mbox{and}~~
\frac{\pl m_{\iy}}{\pl s_k}=- m_{\iy} (\Lb^{\top})^k,
\label{ode}\ee
 with given initial condition $ m_{\iy}(0,0)$.
 According to \cite{AvM1}, the unique (formal) solution to this problem is given by
\be m_{\iy}(t,s)=e^{\sum_1^{\iy} t_i \Lb^i}m_{\iy}(0,0)
 e^{-\sum_1^{\iy} s_i \Lb^{\top i}}%=S_1^{-1}(t,s)S_2(t,s),
,\label{infinitemomentmatrix}\ee
 where
$$
 e^{\sum_{1}^{\iy}t_i \Lb^i}=\sum_0^{\iy} \Lb^i {\bf s}_i(t)=
\Bigl ({\bf s}_{j-i}(t)\Bigr)_{\tiny\begin{array}{l}
\scriptsize1\leq i < \iy\\ \scriptsize1\leq j<\iy
\end{array}},
 $$
 is a matrix of Schur polynomials\footnote{See footnote 2
 .%The elementary Schur polynomials ${\bf s}_i$ are defined by
% $e^{\sum^{\infty}_1 t_i z^i }=\sum^{\infty}_0 {\bf
%s}_k(t)z^k$ and ${\bf s}_k(t)=0$ for $k<0$.
} ${\bf
s}_i(t)$, of which a truncated version is given by the
following $n\times\iy$ submatrix:
\begin{eqnarray}
E_n(t)&=&\left(\begin{array}{ccccc|cc} 1&{\bf
s}_1(t)&{\bf s}_2(t)&\pp&{\bf s}_{n-1}(t)
 &{\bf s}_n(t)&\pp\\
0&1&{\bf s}_1(t)&\pp&{\bf s}_{n-2}(t)&{\bf s}_{n-1}(t)&\pp\\
\vdots&\vdots& & & \vdots&\vdots& \\
0&0&0&\pp&{\bf s}_1(t)&{\bf s}_2(t)&\pp\\
0&0&0&\pp&1&{\bf s}_1(t)&\pp
\end{array}\right)
\nonumber\\ &=&\Bigl({\bf s}_{j-i}(t)\Bigr)_{ 1\leq
i\leq n \atop 1\leq j<\iy }\label{3.0.6}
\end{eqnarray}
Then the $n \times n$ upper-left corner
$m_n(t,s):=\left( \mu_{ij}(t,s)\right)_{0\leq i,j\leq
n-1}$ of $m_{\iy}(t,s)$ is given by
%So, for the
%semi-infinite initial condition $m_{\iy}(0,0)$, setting
%$ m_{n}(t,s):=\left( \mu_{ij}(t,s)\right)_{0\leq i,j\leq
%n-1}$,we have
%the $\tau$-functions of the
%2-Toda problem are given by
%
 \be
m_{n}(t,s)
 =  E_n(t)~m_{\iy}(0,0)~E_n^{\top}(-s)
.  \label{finitemomentmatrix}\ee
%Incidentally, the wave vectors $\Psi_1$ and
%$\Psi_2^*$ define monic polynomials $p^{(1)}(x)$ and
%$p^{(2)}(y)$, \be
%\begin{tabular}{lll}
%$\Psi_1:=e^{\sum t_kz^k}p^{(1)}(z)$&\mbox{and}&
% $\Psi_2^*:=e^{-\sum
%s_kz^{-k}}h^{-1}p^{(2)}(z^{-1}) $   \\
%$\,\,\,~~~=e^{\sum t_kz^k}S_1\chi(z)$&
%&$\,\,\,~~~=e^{-\sum
%s_kz^{-k}}(S_2^{-1})^{\top}\chi(z^{-1}),$
%\end{tabular}
%\ee which are {\em bi-orthogonal} with regard to the
%original matrix $m_{\iy}$; that is, for all $t,s$:
%\be\la p_n^{(1)},p_m^{(2)}\ra =\dt_{n,m}h_n~~\mbox{for
%the inner-product defined by}~~
% \la x^i,y^j  \ra:= \mu_{ij},
% \ee
%with $h_n$ as in (2.1.8).

%\newpage
\vspace{.5cm}

%\newpage

\begin{lemma}  Given the semi-infinite initial
condition $m_{\iy}(0,0)$, and given integers
 $$
 1\leq a_1<\ldots <a_n~~\mbox{and}~~
 1\leq b_1<\ldots <b_n,
 $$
 the
following determinants have a ``Fourier" expansion in
Schur and skew-Schur polynomials\footnote{The sum below
is taken over all Young diagrams $\lb$ and $\nu$, with
the first columns $\lb_1^{\top}$ and $\nu_1^{\top}\leq n
$.}

\be
  \det (\mu_{k,{\ell}}(t,s))_{1\leq k,\ell\leq
n}
 = \sum_{\lb,~ \nu \atop
  \lb^{\top}_1,~  \nu^{\top}_1 \leq n}
    \det(m_{}^{\lb,\nu}(0,0)) s_{\lb}(t) s_{\nu}(-s),
 ~~ \mbox{for $n>0$},
 \label{minor1}\ee
%where the sum is taken over all Young diagrams $\lb$ and
%$\mu$, with first columns $\lb^{\top}_1$ and
%$\mu^{\top}_1\leq n $ and where
%  \be
%   m_{}^{\lb,\mu}:=\left( \mu_{\lb_{i}-i+n,\mu_{j}-j+n}
%    \right)_{1\leq i,j\leq n}.
% \ee

%function $\tau(t,s)$ has the following expansion in
%Schur polynomials\footnote{For a given Young diagram
%$\lb_1 \geq ...\geq\lb_n$, define $s_{\lb}(t)=\det
%(p_{\lb_i-i+j}(t))_{1\leq i,j \leq n}$. },
 \be
%\tau_n(t,s)
\det (\mu_{a_k,b_{\ell}}(t,s))_{1\leq k,\ell\leq n}
% =\det \left(
%E_n(t)~m~E_n^{\top}(-s)\right)
 = \sum_{
 {{\lb,~ \nu}\atop {\lb \supset \al,~~ \nu\supset \beta}}
 \atop
  \lb^{\top}_1,~  \nu^{\top}_1 \leq n}
    \det(m_{}^{\lb,\nu}(0,0))
    \gs_{\lb\backslash \alpha}(t)
     \gs_{\nu\backslash \beta}(-s)
%     s_{\lb}(t) s_{\nu}(-s),
 ~~ \mbox{for $n>0$}, \label{minor2}\ee
with Fourier coefficients, involving the matrices
  \be
   m_{}^{\lb,\nu}(0,0):=\left( \mu_{\lb_{i}-i+n,\nu_{j}-j+n}
    (0,0)\right)_{1\leq i,j\leq n}~\mbox{for
     $ \lb_1^{\top},\nu_1^{\top}\leq n$}
    ,
 \ee
  and where $\al$ and $\beta$ are partitions defined by
  \be
  a_j= \alpha_{n-j+1}+j  ~~\mbox{and}~~
 b_j=   \beta_{n-j+1}+j.
  \ee
 \end{lemma}

\proof Note that every strictly increasing sequence
$1\leq k_1<...<k_n<\iy$ of integers can be mapped into a
Young diagram $\lb_1\geq \lb_2\geq...\geq \lb_n\geq 0$,
by setting $k_j=j+\lb_{n+1-j}$; and similarly with the
increasing sequences $\ell_j,a_j,b_j$. So, this leads to
partitions $\lb, \nu,\alpha , \beta$:
  \bean
 k_j&=& \lb_{n-j+1}+j\\
 \ell_j&=& \nu_{n-j+1}+j\\
 a_j&=& \alpha_{n-j+1}+j\\
 b_j&=& \beta_{n-j+1}+j.
 \eean
 It is also useful to
relabel the index $i$ with $1\leq i\leq n$, by setting
$i':=n-i+1$, also with $1\leq i'\leq n$ and so on with
the other indices.

%\newpage

 Applying the Cauchy-Binet formula twice, and setting
  $$
  A_{a_1a_2\ldots a_n}:= \mbox{matrix formed with the
  rows $a_1\ldots a_n$ of} ~A
  ,$$
  the expression (\ref{finitemomentmatrix}) leads to:
  \bean
 \lefteqn{\det \left( m_{n}(t,s)_{a_i,b_j}\right)_{1\leq i,j \leq n}
     }\\ %&&\\
 &=& \det \left( E(t)_{a_1a_2\ldots a_n}~m_{\iy}(0,0)
  ~\left(E_n(-s)_{b_1b_2\ldots b_n}\right)^{\top}\right)
       \\
 &=& %\hspace{-.3cm}
    \sum_{1\leq k_1<...<k_n<\iy}\det
    \left(p_{k_j-a_i}(t)\right)_{1\leq i,j \leq n}
    \det
    \left(\left(m_{\iy}(0,0)
     E(-s)_{b_1b_2\ldots b_n}\right)^{\top}
      )_{k_i,b_r}
    \right)_{1\leq i,r\leq n}
 \\
 %\eean
 %\bean
 &=& % \hspace{-.3cm}
    \sum_{1\leq k_1<...<k_n<\iy}\det
    \left(p_{k_j-a_i}(t)\right)_{1\leq i,j \leq n}
    \det
    \left( (\mu_{k_i-1,j-1})_{ 1\leq i\leq n\atop
    1\leq j<\iy }
      (p_{i-b_r}(-s))_{ 1\leq i<\iy \atop 1\leq r\leq n
     })\right) \\
 &=& %\hspace{-.3cm}
   \sum_{1\leq k_1<...<k_n<\iy}\det
    \left(p_{k_j-a_i}(t)\right)_{1\leq i,j \leq n}\\
 &&~~~~~~  \sum_{1\leq \ell_1<...<\ell_n<\iy}\det
    \left( \mu_{k_i-1,\ell_j-1}\right)_{1\leq i,j \leq n}
    \det
    \left( p_{\ell_j-b_r}(-s)\right)_{1\leq j,r \leq n}
    \\
% &=&  @@ %\hspace{-.3cm}
%   \sum_{\lb \atop   \lb^{\top}_1\leq n}\det
%    \left(p_{(\lb_{n+1-j}+j)-(\al_{n+1-i}+i)}(t)
%    \right)_{1\leq i,j \leq n}
%    \\
% &&~~~~~~  \sum_{\nu \atop   \nu^{\top}_1\leq n}\det
%    \left( \mu_{\lb_{i}-i+n,\nu_{j}-j+n}
%    \right)_{1\leq i,j \leq n}
%    \det
%    \left( p_{(\nu_{n+1-j}+j)-(\beta_{n+1-r}+r)
%     }(-s)
%    \right)_{1\leq j,r \leq n}
  %   \\
 \eean
 \bean
 &=&
    \sum_{\lb,~ \nu \atop \lb^{\top}_1,~ \nu^{\top}_1 \leq n}
    \det\left( \mu_{\lb_{i'}-i'+n,\nu_{j'}-j'+n}
    \right)_{1\leq i',j' \leq n}
     \\
     && \hspace{1cm}
     \det
    \left(p_{(\lb_{j'}-j')-(\al_{i'}-i')}(t)
    \right)_{1\leq i,j \leq n}
     \det
    \left( p_{(\nu_{ j'}-j')-(\beta_{r'}-r')
     }(-s)
    \right)_{1\leq j,r \leq n}\\
 &=&
    \sum_{\lb,~ \nu \atop \lb^{\top}_1,~ \nu^{\top}_1 \leq n}
    \det\left( \mu_{\lb_{i'}-i'+n,\nu_{j'}-j'+n}
    \right)_{1\leq i',j' \leq n}
     \gs_{\lb\backslash \alpha}(t)
     \gs_{\nu\backslash \beta}(-s),
\eean
 establishing Lemma 3.2. \qed

%\begin{proposition} Given the
%semi-infinite initial condition
%  $m_{\iy}(0,0)=(\mu_{ij})_{i,j\geq 0} $,
% and given
%\be E_n(t):=
% ({\bf s}_{j-i}(t))_{{1\leq i\leq n}\atop {1\leq j\leq \iy}}
% ,\label{Ematrix}\ee
%   the determinant of the matrix (incidentally, a $2$-Toda $\tau$-function)
% \be
% m_{n}(t,s)
% :=  E_n(t)~m_{\iy}(0,0)~E_n^{\top}(-s)
% \label{m-matrix}
%\ee
%  admits the following expansion in Schur polynomials
%\footnote{For a
%given Young diagram $\lb_1 \geq ...\geq\lb_n$, define
%$s_{\lb}(t)=\det (p_{\lb_i-i+j}(t))_{1\leq i,j \leq n}$. },
% \be
%\tau_n(t,s):= \det m_n(t,s)
% =\det \left(
%E_n(t)~m~E_n^{\top}(-s)\right)
% = \sum_{\lb,~ \mu \atop
%\hat \lb_1,~\hat \mu_1 \leq n}
%    \det(m_{}^{\lb,\mu}) s_{\lb}(t) s_{\mu}(-s),
% ~~ \mbox{for $n>0$}, \ee
%where the sum is taken over all Young diagrams $\lb$ and
%$\mu$, with first columns $\lb^{\top}_1$ and
%$\mu^{\top}_1\leq n $ and where
%  \be
%   m_{}^{\lb,\mu}:=\left( \mu_{\lb_{i}-i+n,\mu_{j}-j+n}
%    \right)_{1\leq i,j\leq n}.
% \ee

%\end{proposition}

%\vspace{1cm}

{\medskip\noindent{\it Proof of Proposition 3.1:\/} }
%\proof
 In the integral
$$
\int_{U(n)}e^{\sum_1^{\iy}\Tr(t_jM^j-s_i\bar M^j)}dM,
$$
the shifts $t_i\mapsto t_i+t_i^{(0)}$, $s_i\mapsto
s_i+s_i^{(0)}$ lead to a Toeplitz matrix of the form
  \bean
  \lefteqn{
   \int_{U(n)}e^{\sum_1^{\iy}\Tr(t_j^{(0)}M^j-s_j^{(0)}\bar
M^j)} e^{\sum_1^{\iy}\Tr(t_j M^j-s_j\bar M^j)}dM
  }
  \\
  &=&
  \int_{U(n)}e^{%\sum_1^{\iy}
   \Tr V(M)}
%  (t_j^{(0)}M^j-s_j^{(0)}\bar M^j)}
 e^{\sum_1^{\iy}\Tr(t_j M^j-s_j\bar M^j)}dM
 \\
  &=&
   \det
   \left( \oint u^{\ell-k}e^{V(u)}e^{\sum_1^{\iy}
    (t_ju^j-s_ju^{-j})}  \frac{du}{2\pi i~u}
    \right)_{1\leq \ell,k\leq n}
   \eean
 with
 \be
   V(u):=\sum_1^{\iy}(t_j^{(0)}u^j-s_j^{(0)}u^{-j})
   \ee
   The following matrix of integrals has the form
(\ref{infinitemomentmatrix}), using footnote 2 and
  (\ref{3.0.6}),
 \bean
 \lefteqn{
  \left( \oint u^{\ell-k}e^{V(u)}e^{\sum_1^{\iy}
    (t_ju^j-s_ju^{-j})}  \frac{du}{2\pi i~u}
    \right)_{1\leq \ell,k\leq n}
       }\\
 %\eean  \bean
      &=&
     \left( \sum_{\al,\beta=0}^{\iy}{\bf s}_{\al}(t)
          {\bf s}_{\beta}(-s)
          \oint u^{\ell-k+\al-\beta}e^{V(u)}
            \frac{du}{2\pi i~u}
       \right)_{ 1\leq \ell,k\leq n}
       \\
      &=&
      E_n(t) \left(
        \oint_{S^1} u^{\ell-k}e^{V(u)}
          \frac{du}{2\pi i~u}
       \right)_{ 1\leq \ell,k\leq \iy}
        E_n(-s)^{\top}
  \eean
  or, alternatively, the matrix above satisfies
   the differential equations
  (\ref{ode}), with initial condition
  $$
  \left(
         \oint_{S^1} u^{\ell-k}e^{V(u)}
          \frac{du}{2\pi i~u}
       \right)_{ 1\leq \ell,k\leq n}
   .$$
Now, using Lemma 3.2, its determinant admits a Fourier
expansion in Schur polynomials
\bea \lefteqn{\int_{U(n)}e^{\sum_1^{\iy}
 \Tr(t_j^{(0)}M^j-s_j^{(0)}\bar M^j)}
 e^{\sum_1^{\iy}\Tr(t_j M^j-s_j\bar M^j)}dM
   }\nonumber\\
& & \nonumber\\
&=&\sum_{\lb,\mu~\mbox{\tiny such
 that}\atop{\lb_1^{\top},\mu_1^{\top}\leq n}}
 a_{\lb\mu}(t^{(0)},s^{(0)})s_{\lb}(t)s_{\mu}(-s),
  \eea
  where the Fourier coefficients can be expressed in two
  different ways, first as a determinant,
  using (\ref{minor1}), and secondly as a Fourier
  series, to be explained below,
  \bean
  a_{\lb\mu}(t^{(0)},s^{(0)})
   &=&\det\left(\oint_{S^1}u^{\lb_{\ell}-\ell
-\mu_k+k}e^{\sum^{\iy}_1(t^{(0)}_ju^j-s^{(0)}_ju^{-j})}\frac{du}{2\pi~iu}\right)_{1\leq
\ell,k\leq n}  \\
   &=&
   \sum_{\nu~\mbox{\tiny with}\atop{\nu_1^{\top} \leq
n}}
 {\bf s}_{\nu\ba\lb}(t^{(0)})
  {\bf s}_{\nu\ba\mu}(-s^{(0)})
  .  \eean
To prove the second expression above, we apply
(\ref{minor2}) of Lemma 3.2. Indeed, switching points of
view (i.e., $t \rightarrow t^{(0)}, ~s \rightarrow
s^{(0)}$),
 the initial condition for the differential
 equation (\ref{ode}) is given here by
 $$
 \left. \left(
         \oint_{S^1} u^{\ell-k}
         e^{\sum^{\iy}_1(t^{(0)}_ju^j-s^{(0)}_ju^{-j})}
          \frac{du}{2\pi i~u}
       \right)_{ 1\leq \ell,k\leq n}
       ~~\right|_{t^{(0)}=s^{(0)}=0}= I_n
   ,$$
   which implies that in the sum (\ref{minor2}), the
   Fourier coefficients all vanish, except when
   $\lb=\nu$, for which they equal $1$.
   Moreover, the
sequences $1\leq a_1<a_2<\ldots < a_n$ and $1\leq
b_1<b_2<\ldots < b_n$ are given by
  \bean
 && 1\stackrel{*}{\leq} \lb_{n}-n+c<
          \lb_{n-1}-(n-1)+c<\ldots<\lb_{1}-1+c
  \\
 && 1\stackrel{*}{\leq} \mu_{n}-n+c<
          \mu_{n-1}-(n-1)+c<\ldots<\mu_{1}-1+c
   \eean
   with $c$ the smallest integer for which the
   inequality
   $\stackrel{*}{\leq}$ is satisfied; e.g.,
   \bean
    a_j&=&\lb_{n-j+1}-(n-j+1)+c\\
       &=& \lb_{n-j+1} +j-(n+1-c).
    \eean
 The orthonormality of the Schur
 polynomials\footnote{where the inner-product is defined by
 \be
  \la f(t),g(t) \ra:=
    \left.f(\tilde \pl_t)g(t)\right|_{t=0},
    ~~\mbox{with}
     \tilde \pl_t:=(\frac{\pl}{\pl t_1},\frac{1}{2}\frac{\pl}{\pl t_2},
     ,\frac{1}{3}\frac{\pl}{\pl t_3},
     \ldots )\label{inproduct}\ee}
  $\la {\bf s}_{\lb},{\bf s}_{\mu}\ra=\delta_{\lb\mu}$
  allows for an
  alternative way of expressing
  the ``Fourier" coefficients
$a_{\lb\mu}(t^{(0)},s^{(0)})$ in (\ref{Prop3.1}):
\bean
 \lefteqn{a_{\lb\mu}(t^{(0)},s^{(0)})}\\
 &=& \left\la {\bf s}_{\lb}(t){\bf s}_{\mu}(-s),
  \sum_{\al,\beta~\mbox{\tiny such
 that}\atop{\al_1^{\top},\beta_1^{\top}\leq n}}
 a_{\al\beta}(t^{(0)},s^{(0)})s_{\al}(t)s_{\beta}(-s)
 \right\ra\\
 &=&
 {\bf s}_{\lb}(\tilde\partial_t)
  {\bf s}_{\mu}(-\tilde\partial_s)
   \sum_{\al,\beta~\mbox{\tiny such
 that}\atop{\al_1^{\top},\beta_1^{\top}\leq n}}
 a_{\al\beta}(t^{(0)},s^{(0)})s_{\al}(t)s_{\beta}(-s)
 \Big|_{t=0,~s=0}\\
& & \\
 &=&{\bf s}_{\lb}(\tilde\partial_t)
  {\bf s}_{\mu}(-\tilde\partial_s)
    \int_{U(n)}e^{\sum_1^{\iy}
 \Tr(t_j^{(0)}M^j-s_j^{(0)}\bar M^j)}
 e^{\sum_1^{\iy}\Tr(t_j M^j-s_j\bar M^j)}dM
 \Big|_{t=0,~s=0}
 \\
 &&\\
&=&\int_{U(n)}{\bf s}_{\lb}(M){\bf s}_{\mu}(\bar M)e^{\sum^{\iy}_1\Tr(t^{(0)}_jM^j-s^{(0)}_j\bar M^j)}dM,\\
  \eean
 upon using the identity\footnote{Thus
                           $\gs_{\lb}(M)$ is viewed as a
                           symmetric function of the
                           eigenvalues $x_1,\ldots,x_n$
                           of the unitary matrix $M$}
\bean
  {\bf s}_{\lb}(\tilde\partial_t)e^{\Tr\sum_1^{\iy}
t_jM^j}\Big|_{t=0}
 &=&
  {\bf s}_{\lb}\left(\Tr M,\frac{1}{2}\Tr
M^2,\frac{1}{3}\Tr M^3,...\right)
 =: {\bf s}_{\lb}(M),\\
% &=:& {\bf s}_{\lb}(M),
\eean
 thus ending the proof of Proposition 3.1.\qed

%\newpage

\section{Matrix integrals and Random Walks}

   In this section, we consider some interesting special cases,
   based on special values of $t$ and $s$, for which the
   skew Schur polynomials take on the following form:
\bean
 {\bf s}_{\lb\ba\al}(t)\Big|_{t_i=u\delta_{i1}}
 &=&
  u^{|\lb\ba \al|}\gs_{\lb} (1,0,\ldots)\\
 &=&
  \frac{u^{|\lb\ba\al |}}{|\lb\ba\al |!}
   \qquad\#\left\{\begin{array}{l}
\mbox{standard skew-tableaux}\\
\mbox{of shape $\lb\ba\al$, filled}\\
\mbox{with numbers $1,...,|\lb \ba\al |$}
\end{array}\right\}\\
& & \\
&=&\frac{u^{|\lb\ba\al |}}{|\lb\ba\al |!}f^{\lb\ba\al}\\
\eean
%
%\vspace*{2cm}
and
 \bean
 \gs_{\lb\ba\al}(t)\Big|_{it_i=qu^i}
 &=&u^{|\lb\ba\al
|}~\gs_{\lb\ba\al}
\left(q,\frac{q}{2},\frac{q}{3},...\right)
\\
&=&%&=&
 u^{|\lb\ba\al |}
\qquad\#\left\{\begin{array}{l}
\mbox{semi-standard skew-tableaux}\\
\mbox{of shape $\lb\ba\al$, filled}\\
\mbox{with numbers $1,...,q$}
\end{array}\right\}.\\
%& &\,\,\,\downarrow\,\,\,\uparrow \\
%& & \\
\eean We now study the three integrals, appearing in the
introduction:

%\newpage

\subsection{Integral 1:
  $\displaystyle{\int_{U(n)}}s_{\lb}(M)s_{\mu}(\bar M)e^{z\Tr(M+\bar M)}dM
$}

%\noindent{\bf Locus 1}: $\LR_1=\{$all
%$t_k^{(0)}=s_k^{(0)}=0$, except $t_1^{(0)}=z$,
%$s_1^{(0)}=-z\}$.

%\noindent{\bf Example 1}:
{\em A generating function for

\be
b^{(k)}_{xy}=
 \#\left\{\begin{array}{l} \mbox{ways that
$n$ non-intersecting walkers
  in $\mathbb Z$ move during}\\
\mbox{$k$ instants from $x_1<x_2< ... <x_n$ to $y_1<y_2<
... <y_n$,}\\
 \mbox{where at each instant exactly one walker moves}\\
 \mbox{either one step to the left, or
one step to the right}
%   \\ \mbox{leading to $k$ effective moves  .}
\end{array}\right\}
\ee
%
%
%
%
%$$b^{(k)}_{xy}=
%  \#\left\{\begin{array}{l}
%\mbox{ways that $n$ non-intersecting}\\
%\mbox{walkers in $\mathbb Z$ move in $k$ steps}\\
%\mbox{from $x_1<x_2< ... <x_n$}\\
%\mbox{to $y_1<y_2< ... <y_n$}
%\end{array}\right\}
%$$
 is given by the matrix integral
$$
\sum_{k\geq 0}\frac{z^k}{k!}b^{(k)}_{xy}=
 \int_{U(n)}s_{\lb}(M)s_{\mu}(\bar M)e^{z\Tr(M+\bar M)}dM
=:a_{\lb\mu}(z).$$
Furthermore, a generating function for the
$a_{\lb\mu}$'s is given by
$$
\sum_{{\lb,\mu\mbox{\tiny ~such
that}}\atop{\lb_1^{\top},\mu_1^{\top}\leq
n}}a_{\lb\mu}(z)s_{\lb}(t)s_{\mu}(-s) =
 \int_{U(n)}e^{z\Tr(M+\bar M)}e^{\sum_1^{\iy}\Tr(t_iM^i-s_i\bar M^{i})}dM\\
,$$
%
%$$
%\int_{U(n)}s_{\lb}(M)s_{\mu}(\bar M)e^{z\Tr(M+\bar M)}dM
%= \sum_{k\geq 0}\frac{z^k}{k!}\#\left\{\begin{array}{l}
%\mbox{ways that $n$ non-intersecting}\\
%\mbox{walkers in $\mathbb Z$ move in $k$ steps}\\
%\mbox{from $x_1<x_2< ... <x_n$}\\
%\mbox{to $y_1<y_2< ... <y_n$}
%\end{array}\right\}
%$$
 where
 \be
 \mu_{n-k+1}:=x_k-k+1~,~~~\lb_{n-k+1}:=y_k-k+1.
  ~~\mbox{for} ~k=1,\ldots,n.  \label{4.1.2}\ee
%   \bean
%\begin{array}{ll}
%   k_1=\frac{1}{2}(k-|\lb|+|\mu|),
% & k_2=\frac{1}{2}(k+|\lb|-|\mu|)\\
% \mu_{n-k+1}:=x_k-k,&\lb_{n-k+1}:=y_k-k.
%\end{array}
%  \eean

}

\proof Consider the locus
 $$\LR_1=\{\mbox{all
$t_k^{(0)}=s_k^{(0)}=0$, except $t_1^{(0)}=z$,
$s_1^{(0)}=-z\}$}.$$
 Then, since
$$
e^{\sum_1^{\iy}(t_i^{(0)}u^i-s_i^{(0)}u^{-i})}\Big|_{\LR_1}=e^{z(u+u^{-1})},
$$
we have, using (\ref{generalmatrixintegral}),
\be \int_{U(n)}e^{z\Tr(M+\bar M)}
  e^{\sum_1^{\iy}\Tr(t_iM^i-s_i\bar M^{i})}dM
   =
   \sum_{{\lb,\mu\mbox{\tiny ~such
that}}\atop{\lb_1^{\top},\mu_1^{\top}\leq
n}}a_{\lb\mu}(z)s_{\lb}(t)s_{\mu}(-s), \ee
with
\bean
a_{\lb\mu}(z)&=&\int_{U(n)}s_{\lb}(M)s_{\mu}(\bar M)e^{z\Tr(M+\bar M)}dM\\
& & \\
&=&\det\left(\oint_{S^1}u^{\lb_{\ell}-\ell-\mu_k+k}e^{z(u+u^{-1})}\frac{du}{2\pi~iu}\right)_{1\leq\ell,k\leq n}\\
& & \\
 &=&\sum_{{\nu\mbox{\tiny
~with}}\atop{\nu_1^{\top}\leq n}}
   s_{\nu\ba\lb}(t^{(0)})s_{\nu\ba\mu}(-s^{(0)})\Big|_{\LR_1}\\
& & \\
&=&\sum_{{\nu\mbox{\tiny
~with~}\nu\supset\lb,\mu}\atop{\nu_1^{\top}\leq
n}}\frac{z^{|\nu\ba\lb|}}{|\nu\ba\lb|!}~f^{\nu\ba\lb}\frac{z^{|\nu\ba\mu|}}{|\nu\ba\mu|!}~f^{\nu\ba\mu}\\
& & \\
&=&\sum^{\iy}_{k=0}\frac{z^k}{k!}\frac{k!}{k_1!~k_2!}\sum_{
  {{\nu~\mbox{\tiny with}~\nu \supset \lb,\mu}\atop
 {   { {|\nu\ba \lb|=k_1}\atop{|\nu\ba \mu|=k_2}}
 \atop {\nu_1^{\top} \leq n}
   }    }
   }f^{\nu\ba\lb}f^{\nu\ba\mu}
%   \\& & \\
 \eean \bean&=&\sum_{k\geq 0}\frac{z^k}{k!}\#\left\{\begin{array}{l}
\mbox{ways that $n$ non-intersecting}\\
\mbox{walkers in $\mathbb Z$ move in $k$ steps}\\
\mbox{from $x_1<x_2< ... <x_n$}\\
\mbox{to $y_1<y_2< ... <y_n$}
\end{array}\right\}
\eean
 as a consequence of (\ref{2.1.2}), where
   \bean
\begin{array}{ll}
   k_1=\frac{1}{2}(k-|\lb|+|\mu|),
 & k_2=\frac{1}{2}(k+|\lb|-|\mu|).
\end{array}
  \eean
An alternative way of proving the final formula is to
invoke the D. Andr\'e reflection principle. Indeed
 \bean
 a_{\lb\mu}(z)&=& \det(m^{\lb,\mu})\\
%&=&
%  \det \left( J_{(\lb_i-i)-(\mu_j-j)}
%   (z)\right)_{1\leq i,j\leq n}
%\\
%&=&
%   \det \left( J_{y_i-x_j} (z)
%     \right)_{1\leq i,j \leq n},~~~~~\mbox{using (4.0.1)}
%\\
&=&
   \det \left( \left. e^{z(u+u^{-1})}\right|_{u^{y_i-x_j}}
     \right)_{1\leq i,j \leq n}
     \mbox{using (\ref{Prop3.1}) and (\ref{4.1.2})}\\
&=&\sum_{w\in W} (-1)^{\sigma(w)}
  \prod_1^n \left.
   e^{z(u_i+u_i^{-1})}\right|_{u_i^{y_i-x_{w(i)}}}
   \\
&=& \sum_{k=0}^{\iy}\# \left\{
\begin{array}{l}
  \mbox{walks of $k$ steps}\\
  \mbox{from $x \mapsto y$ in $\BZ^n$}\\
  \mbox{within $\{ u_1 <\ldots
  <u_n\}$}
  \end{array} \right\}  \frac{z^k}{k!},~\mbox{using
    Theorem 2.1}
  \\
&=& \sum_{k\geq 0} %\sum_{k=0}^{\iy}
  \#\left\{
   \begin{array}{l}
     \mbox{ways that $n$ non-intersecting walkers}\\
     %\mbox{(at each step only one walker walks)}\\
     \mbox{%by $n$ walkers
     in $\BZ$ move in $k$ steps%, one at the time
     }\\
     \mbox{from $x_1<x_2<...<x_n$ }\\
     \mbox{to $y_1<y_2<...<y_n$%in $k$ steps
     }\\
     \end{array}
     \right\}~\frac{z^k}{k!}\\&=&
     \sum_{k\geq 0}b^{(k)}_{x,y} ~~\frac{z^k}{k!} ,
\eean

%\newpage

\subsection{Integral 2:
$\displaystyle{\int_{U(n)}}s_{\lb}(M)s_{\mu}(\bar
M)\det(I+M)^q
  e^{z \Tr\bar M}dM$}

%\noindent{\bf Locus 2}:

%\noindent{\bf Example 2}:
%
{\em A generating function for
\be
 b_{xy}^{(k)}=
  \# \left\{
\begin{array}{l}
  \mbox{ways that n non-intersecting walkers move during $q+k$}\\
  \mbox{instants from $x_1< ... <x_n$ to $y_1< ... <y_n$,
   where at}\\
  \mbox{the instants $1$ to $q$, walkers
  \underline{may} move one step to the}\\
  \mbox{right, or \underline{stay put},
   and at the instants $q+1,...,q+k$
   }\\
  \mbox{\underline{exactly one} walker moves one step to the left, with}\\
  \mbox{total $\#\{\mbox{effective moves}\}=2k+\sum_1^n(y_i-x_i)
  $.}
  \end{array} \right\}
  \label{4.0.4}\ee
%
%  @@@@@@@@@@@
%$$b^{(k)}_{xy}=
%  \# \left\{
%\begin{array}{l}
%  \mbox{ways that n non-intersecting walkers move}\\
%  \mbox{from $x_1< ... <x_n$ to $y_1< ... <y_n$,
%  where at}\\
%  \mbox{each instant $1,...,q=T_R$ walkers
%   \underline{may move}}\\
%  \mbox{\underline{one step to the right or stay put}
%  and at each}\\
%  \mbox{instant $q+1,...,q+k=T_R+T_L$ \underline{exactly one}}\\
%  \mbox{walker moves one step to the left.}
%  \end{array} \right\}
%$$
 is given by the matrix integral
$$
\sum_{k\geq 0}\frac{z^k}{k!}b^{(k)}_{xy}=
 \int_{U(n)}s_{\lb}(M)s_{\mu}(\bar M)\det(I+M)^q
  e^{z \Tr\bar M}dM
=:a_{\lb\mu}(z).$$
Furthermore, a generating function for the
$a_{\lb\mu}$'s is given by
$$
\sum_{{\lb,\mu\mbox{\tiny ~such
that}}\atop{\lb_1^{\top},\mu_1^{\top}\leq
n}}a_{\lb\mu}(z)s_{\lb}(t)s_{\mu}(-s) =
 \int_{U(n)}\det(I+M)^q
  e^{z \Tr\bar M}
   e^{\sum_1^{\iy}\Tr(t_iM^i-s_i\bar M^{i})}dM\\
,$$
  with
  \bean
\begin{array}{ll}
%T_L=q,&T_R=k\\
\mu_{n-j+1}:=x_j-j+1~,~~&~\lb_{n-j+1}:=y_j-j+1,~~\mbox{for
$j=1,\ldots,n$.}
\end{array}
\eean

}

 \proof Consider the locus
  $$
  \mbox{$\LR_2=\{$all
   $it^{(0)}_i=-q(-1)^i$, $is^{(0)}_i=-z\delta_{i1}\}$.}
 $$
 Then, using
$$
e^{\sum_1^{\iy}t_i u^i}\Big|_{it_i=qx^i}=e^{q\sum_{i\geq
1}\frac{(xu)^i}{i}}=(1-xu)^{-q},
$$
we have by (\ref{generalmatrixintegral}) that
\bean \lefteqn{\int_{U(n)}\det(I+M)^q
  e^{z \Tr\bar M}
   e^{\sum_1^{\iy}\Tr(t_iM^i-s_i\bar M^{i})}dM
   }
   \\
& & \\
&=&\sum_{{\lb,\mu\mbox{\tiny ~such
that}}\atop{\lb_1^{\top},\mu_1^{\top}\leq
n}}a_{\lb\mu}(z)s_{\lb}(t)s_{\mu}(-s). \eean
Using
\bean
  {\bf s}_{(\lb \ba \al)^{\top}}(t)
    &=&
        (-1)^{|\lb\ba \al|}
            {\bf s}_{\lb \ba \al}(-t)  \\
   {\bf s}_{\lb \ba \al}  (-t_1,t_2,-t_3,\ldots)
    &=&
        (-1)^{|\lb\ba \al|}
            {\bf s}_{\lb \ba \al} (t),
  \eean
the Fourier coefficients $a_{\lb\mu}(z)$ of
(\ref{Prop3.1}) have the following interpretation:
\bean
a_{\al\beta}(z)&=&\det\left(\int_{S^1}u^{\al_{\ell}-\ell-\beta_k
+k}e^{\sum_1^{\iy}(t^{(0)}_j u^j-s^{(0)}_j
u^{-j})}\frac{du}{2\pi~iu}\right)_{1\leq \ell,k\leq
n}\left|_{
\renewcommand{\arraystretch}{0.5}
\begin{array}[t]{l}
\\
{\scriptstyle it_i=-q(-1)^i}\\
{\scriptstyle s_i=-z\delta_{i1}}
\end{array}
\renewcommand{\arraystretch}{1}
}\right.\\
& & \\
&=&\det\left(\int_{S^1}u^{\al_{\ell}-\ell-\beta_k
+k}(1+u)^{+q}e^{zu^{-1}}\frac{du}{2\pi~iu}
\right)_{1\leq \ell,k\leq n}
  \\
 \eean\bean
 & &
  \\
 &=&\sum_{{\lb\supset\al,\beta}\atop{\lb_1^{\top}\leq
n}}s_{\lb\ba\al}(t)\Big|_{it_i=-q(-1)^i}s_{\lb\ba\beta}(-s)
   \Big|_{s_i=-z\delta_{i1}}
   \\
 & & \\
 &=&
 \sum_{{\lb\supset\al,\beta}\atop{\lb_1^{\top}\leq
n}}s_{(\lb\ba\al)^{\top}}(-t)\Big|_{it_i=-q(-1)^i}
  s_{(\lb\ba\beta)^{\top}}(s)\Big|_{s=-z\delta_{i1}}
   (-1)^{|\lb\ba \al|+|\lb\ba \beta|}\\
 & &
 \\
 &=&
 \sum_{{\lb\supset\al,\beta}\atop{\lb_1^{\top}\leq
n}}s_{(\lb\ba\al)^{\top}}(t)\Big|_{it_i=q(-1)^i}(-1)^{|\lb\ba\al
|}
\frac{z^{|\lb\ba\beta|}}{|\lb\ba\beta|!}f^{\lb\ba\beta}
 \\
 & &
 \\
&=&%(-1)^{|\beta |-|\al|}
  \sum^{\iy}_{k=0}\frac{z^k}{k!}\sum_{
  {{\lb\vdash k+|\beta|}\atop
 {   { {\lb\supset\al,\beta}\atop{\lb_1^{\top}\leq n}}}    }
}s_{(\lb\ba\al)^{\top}}\left(q,\frac{q}{2},\frac{q}{3},...\right)f^{\lb\ba\beta}\\
%& &\hspace{3cm}\mbox{using
%  $|\lb |=|\lb\ba\beta |+|\beta |  $}\\
%
& & \\
&=&  \sum^{\iy}_{k=0}\frac{z^k}{k!}\# \left\{
\begin{array}{l}
  \mbox{all $(P,Q)$, with $P$ semi-standard and}\\
   \mbox{$Q$ standard skew-tableaux of arbitrary}\\
  \mbox{shape $(\lb\ba\al)^{\top}$ and
   $(\lb\ba\beta)^{\top}$, filled}\\
   \mbox{with numbers
               $1,...,q$ and $1,...,k$,
               }\\
  %\mbox{$1,...,q$,}\\
  \mbox{with $\lb  \vdash k+|\beta|$,
  $\lb_1^{\top}\leq n$ and fixed $\al$, $\beta$}
  \end{array} \right\}\\
  &&\\
&=&%(-1)^{|\beta |-|\al |}
   \sum^{\iy}_{k=0}\frac{z^k}{k!}b_{xy}^{(k)},
   \eean

   \vspace{-1cm}\be\label{4.0.5}\ee
%
%   \# \left\{
%\begin{array}{l}
%  \mbox{ways that n non-intersecting walkers move}\\
%  \mbox{from $x_1< ... <x_n$ to $y_1< ... <y_n$,
%  where at}\\
%  \mbox{each instant $1,...,q=T_R$ walkers
%   \underline{may move}}\\
%  \mbox{\underline{one step to the right or stay put}
%  and at each}\\
%  \mbox{instant $q+1,...,q+k=T_R+T_L$ \underline{exactly one}}\\
%  \mbox{walker moves one step to the left.}
%  \end{array} \right\}
%  \\
%%
%\eean
  upon setting
   $$
   x_j:= \al_{n-j+1}+j-1 ~~\mbox{and}~~ y_j:=
   \beta_{n-j+1}+j-1.
   $$
  Then\footnote{The number of ``effective steps"
  counts the actual steps taken by all walkers; i.e.,
  two walkers walk simultaneously count for two steps,
  a walker not walking contributes nothing!}
      \bean
   \mbox{Total} ~\# ~\{\mbox{effective steps}\}=T_R+T_L=
         &=& 2|\lb|-|\al|-|\beta| \\
         &=& 2(|\lb|-|\beta|)-|\al|+|\beta|  \\
         &=& 2k+\sum_1^n (y_i-x_i),
         \eean
proving statement (\ref{4.0.4}).

The last equality in (\ref{4.0.5}) requires a little
explanation. We refer here to figure 1 in section 2.1,
but somewhat modified, with $\mu \rightarrow \al$ and
$\nu \rightarrow \beta$. Namely, the left diagram
corresponds to {\em right} steps and the right diagram
to {\em left} steps. The last row corresponds to the
first walker $x_1=\al_{n}$ and the first to the last
walker $x_{n}=\al_1+n-1$. With this slight change, a
semi-standard skew-tableau of shape $(\lb \ba
\al)^{\top}$, with $\lb_1^{\top}\leq n$, filled with
numbers $1,\ldots, q$, is tantamount to a skew-partition
$\lb \ba \al$, with $\lb_1^{\top}\leq n$, filled with
numbers among $1,\ldots,q$, which are strictly
increasing from left to right and non-decreasing from
top to bottom; not all numbers need to be used. If some
number $1\leq r\leq q$ is not used in the skew Young
diagram, then no walker moves at that instant. So, the
pair $(P,Q)$, as in figure 1, corresponds to a walk
departing from $x_k=\al_{n-k+1}+k-1$ and ending up at
$y_k=\beta_{n-k+1}+k-1$. The most-left walker moves to
the right only at the moments indicated by the integer
appearing in the first row, the second walker moves
right at the moments indicated in the second row, etc...
Thus, corresponding to $P$, two walkers can
simultaneously move right!

 %\newpage

\subsection{Integral 3:
$\displaystyle{\int_{U(n)}}s_{\lb}(M)s_{\mu}(\bar M)
\det (I+z M)^p
 \det (I+z\bar M)^q dM$}

%\noindent{\bf Example 3}:
 %
{\em A generating function for
\be
b_{xy}^{(k)}= k!
  \# \left\{
\begin{array}{l}
  \mbox{ways that n non-intersecting walkers move during $p+q$}\\
  \mbox{instants from $x_1< ... <x_n$ to $y_1< ... <y_n$,
   where at}\\
  \mbox{the instants $1$ to $p$, walkers
  \underline{may} move one step to the}\\
  \mbox{right, or \underline{stay put},
   and at the instants $p+1,...,p+q$
   }\\
  \mbox{walkers
  \underline{may} move one step to the left or \underline{stay put},}\\
  \mbox{with total $\#\{\mbox{effective moves}\}=k
  $.}
  \end{array} \right\}
   \label{4.0.6}\ee
is given by the matrix integral
 $$
\sum_{k\geq 0}\frac{z^k}{k!}b^{(k)}_{xy} =
 \int_{U(n)}s_{\lb}(M)s_{\mu}(\bar M) \det (I+z M)^p
 \det (I+z\bar M)^q dM
=:a_{\lb\mu}(z).$$
Furthermore, a generating function for the
$a_{\lb\mu}$'s is given by
 \bean
 \hspace{-1cm} \lefteqn{
   \hspace{-1cm}\sum_{{\lb,\mu\mbox{\tiny ~such
that}}\atop{\lb_1^{\top},\mu_1^{\top}\leq
n}}a_{\lb\mu}(z)s_{\lb}(t)s_{\mu}(-s)
 ~~~~~~~~}
  \\
 ~~~~~~~~~~~~~~~& =&
 \int_{U(n)}\det (I+z M)^p \det (I+z\bar M)^q
   e^{\sum_1^{\iy}\Tr(t_iM^i-s_i\bar M^{i})}dM,
 \eean
  with
  \bean
\begin{array}{ll}
%T_L=q,&T_R=k\\
\mu_{n-j+1}:=x_j-j+1~,~~&~\lb_{n-j+1}:=y_j-j+1,
\end{array}
\eean
for $j=1,\ldots,n$.

}

\proof Consider the locus:

$$ \mbox{ $\LR_3=\{$all
$it_i=-p(-z)^i$, $is_i=q(-z)^i\}$.}
$$
Then, using
$$
e^{\sum_1^{\iy}(t_iu^i-s_iu^{-i})}\Big|_{\LR_3}
 =(1+zu)^{p}(1+zu^{-1})^{q},
$$
we have by (\ref{Prop3.1}),
\bean
a_{\al\beta}(z)&=&\det\left(\int_{S^1}
 u^{\al_{\ell}-\ell-\beta_k+k}(1+zu)^p
(1+zu^{-1})^q\frac{du}{2\pi~iu}\right)_{1\leq \ell,k\leq n}\\
 & & \\
&=&\sum_{{\lb\mbox{\tiny ~with}}\atop{\lb_1^{\top}\leq
n}}
 s_{\lb\ba\al}(t)s_{\lb\ba\beta}(-s)\Big|_{\LR_3}\\
& & \\
 &=&\sum_{{\lb\mbox{\tiny ~with}}\atop{\lb_1^{\top}\leq n}}
  s_{(\lb\ba\al)^{\top}}(-t)(-1)^{|\lb \ba \al|}
    s_{(\lb\ba\beta)^{\top}}(s)(-1)^{|\lb \ba \beta|}\Big|_{\LR_3}\\
& & \\
% \eean
% \bean
 &=&\sum_{{\lb\mbox{\tiny
~with}}\atop{\lb_1^{\top}\leq
n}}s_{(\lb\ba\al)^{\top}}\left(p,\frac{p}{2},
\frac{p}{3}, ...\right)z^{|\lb\ba\al
|}s_{(\lb\ba\beta)^{\top}}\left(q,\frac{q}{2},
\frac{q}{3},
...\right)z^{|\lb\ba\beta |}\\
& & \\
% \eean
% \bean
  &=&\sum^{\iy}_{k=0}z^k\sum_{{{\lb\mbox{\tiny
~with} }\atop{{{\lb_1^{\top}\leq n}\atop{|\lb|=
  \frac{1}{2} (k+|\al|+|\beta|)}}}} }
s_{(\lb\ba\al)^{\top}}\left(p,\frac{p}{2},
...\right)s_{(\lb\ba\beta)^{\top}} \left(q,\frac{q}{2},
...\right)\\
 \eean
 \bean
 &=& \sum^{\iy}_{k=0}{z^k}\# \left\{
\begin{array}{l}
  \mbox{all $(P,Q)$, with $P$ and $Q$ semi-standard}\\
   \mbox{skew-tableaux of arbitrary shapes $(\lb\ba\al)^{\top}$}\\
  \mbox{and
   $(\lb\ba\beta)^{\top}$, filled with numbers
               $1,...,p$}\\
   \mbox{and $1,...,q$, with $\lb  \vdash \frac{1}{2}(k+|\al|+|\beta|)
  $,}\\
  \mbox{$\lb_1^{\top}\leq n$ and fixed $\al$, $\beta$.}
  \end{array} \right\}\\
  &&\\
&=&%(-1)^{|\beta |-|\al |}
   \sum^{\iy}_{k=0}{z^k}b_{xy}^{(k)}
%   \# \left\{
%\begin{array}{l}
%  \mbox{ways that in non-intersecting}\\
%  \mbox{walkers move from $x_1< ... <x_m$}\\
%  \mbox{to $y_1< ... <y_m$, where at}\\
%  \mbox{each instant $1,...,T_L$ at least}\\
%  \mbox{one walker moves one step to the}\\
%  \mbox{left, and at each instant}\\
%  \mbox{$T_L+1,...,T_L+T_R$ at least one walker}\\
%  \mbox{moves one step to the right.}
%  \end{array} \right\}\\
 \eean
as given in (\ref{4.0.6}).

%\newpage

\section{The action of Virasoro on Schur polynomials}

A border-strip $h\in B(i)$ is a connected skew-shape
$\lb \ba \mu$ containing $i$ boxes, with no $2\times 2$
square. The height of a border strip $h$ is defined as
 \be
 \acht~h:=\#\{\mbox{rows
in~}h\}-1 .\ee
 Consider the Virasoro operator in the variables $t_1,t_2,\ldots
 $ for $k \in \BZ$,
 \be
V_k(t)=\frac{1}{2}\sum_{i+j=k}\frac{\pl^2}{\pl t_i\pl
t_j}+\sum_{-i+j=k}it_i\frac{\pl}{\pl
t_j}+\frac{1}{2}\sum_{-i-j=k}(it_i)(jt_j)
\label{5.0.2}\ee
  In this section we study the action of the
  Virasoro operators on Schur polynomials $ {\bf s}_{\lb}$.
  From the
   {\bf Murnaghan-Nakayama rule }
   \cite{Stanley2}, stated in Proposition 5.2, it follows
   that (Corollary 5.3)
\bea
nt_n~\gs_{\lb}(t)&=&\sum_{{\mu}\atop{\mu\backslash\lb\in
B(n)}}(-1)^{\acht(\mu\backslash\lb)}\gs_{\mu}(t) \no\\
\no\\
\frac{\pl}{\pl
t_n}\gs_{\lb}(t)&=&\sum_{{\mu}\atop{\lb\backslash\mu\in
B(n)}}(-1)^{\acht(\lb\backslash\mu)}\gs_{\mu}(t).
 \label{3.0.14}
\eea
In view of this, one would expect $V_{n} {\bf s}_{\lb}$
to be expressible as
  an {\em infinite sum of Schur polynomials}.
  This is not so: acting with Virasoro
  leads to the same precise sum, except for different
  coefficients:

\begin{theorem} The Virasoro operator acts as follows on Schur
polynomials
 \bea
  V_{-n} {\bf s}_{\lb}&=&\sum_{{\mu}\atop{\mu\backslash\lb\in
B(n)}} d^{(-n)}_{\lb \mu} {\bf s}_{\mu}
 \label{Virasoro10} \\
V_{n} {\bf
s}_{\lb}&=&\sum_{{\mu}\atop{\lb\backslash\mu\in B(n)}}
d^{(-n)}_{\mu\lb } {\bf s}_{\mu}
   \label{Virasoro20}   \eea
with $(n\geq 1)$
  \bea
 d^{(-n)}_{\lb \mu}
  &=& \sum_{i\geq 1}~~
    \sum_{
    \left\{
    \begin{array}{l}
           \nu ~~~~~\mbox{\footnotesize such that} \\
           \lb \backslash \nu \in B(i)\\
           \mu \backslash \nu \in B(n+i) \\
            \lb \backslash \nu  \subset \mu \backslash \nu
            \end{array}
            \right\}}
 (-1)^{\mbox{\footnotesize ht}(\lb \backslash \nu) +
       \mbox{\footnotesize ht}(\mu \backslash \nu)  }\no \\
   &&\no  \\&&
   + \frac{1}{2} \sum_{i=1}^n
   ~~
   \sum_{
   \left\{
   \begin{array}{l}
           \nu ~~~~~\mbox{\footnotesize such that} \\
           \nu \backslash \lb \in B(i)\\
           \mu \backslash \nu \in B(n-i) \\
            %\lb \backslash \nu  \subset \mu \backslash \nu
            \end{array}
         \right\}   }
 (-1)^{\mbox{\footnotesize ht}(\nu \backslash \lb) +
       \mbox{\footnotesize ht}(\mu \backslash \nu)  } . \eea

\end{theorem}

%\newpage

%@@@@@@@@@@@@

%A border-strip is a connected skew-shape with no $2\times 2$
%square.

\vspace{1cm}
 We need a few combinatorial preliminaries. Given
 \be
 \al:= (\al_1,\al_2,\ldots ),~~\al_i \in \BZ_{\geq 0},
 \label{3.0.5}\ee
 define
 \bean
T&:=&\mbox{border-strip tableau of shape $\lb\backslash
\mu$
 and type $\al,~\sum \al_i=|\lb\backslash \mu|$,}\\
   &&\mbox{is an assignment of integers ($\geq 0$) to the boxes of
  $\lb\backslash \mu$, such that}\\
& & \mbox{(1) the positive integers are weakly
increasing from $\left\{
\begin{array}{l}
\mbox{left to right}\\
\mbox{top to bottom}
\end{array}
\right.$}\\
& & \mbox{(2) the integer $i$ appears $\al_i$ times}\\
& & \mbox{(3) $\{$squares containing $i\}$ is a border
strip $B_i$.} \eean
Define
 \bea
 \acht~B_i&=&\#\{\mbox{rows
in~}B_i\}-1
   \no\\
  \acht~T&=&\sum_i \acht~B_i \no\\
   &&\no \\
   \chi^{\lb\backslash \mu}(\al)&=&\sum_{
    \left\{
    \begin{array}{c}
           \mbox{\footnotesize all border-strip tableaux $T$} \\
           \mbox{\footnotesize of shape $\lb\backslash \mu$
             and type $\al$}
%           \\  \mbox{\footnotesize }
            \end{array}
            \right\}}(-1)^{\acht~T}
    \eea
Given $\al$ as in (\ref{3.0.5}), define
$$
p_{\al}=p_{\al_1}p_{\al_2}\ldots
=\left(\sum_ix_i^{\al_1}\right)
\left(\sum_ix_i^{\al_2}\right)\ldots=\al_1
t_{\al_1}\al_2t_{\al_2}\ldots
$$
%
%Remember the inner-product (\ref{inproduct}) in footnote
%6 between symmetric functions, where $kt_k=\sum_{i\geq
%1} x_i^k$.
% \be
%  \langle f(t), g(t) \rangle
% :=\left. f\left(\frac{\pl}{\pl t}, \frac{\pl}{\pl t_2}, \ldots\right)
%  g(t)\right|_{t=0}.
% \label{inproduct}\ee

The Schur polynomials $\gs (t)$ form an orthonormal
basis in the space of symmetric functions, for the
customary inner-product (\ref{inproduct}) in footnote 6,
between symmetric functions, where $kt_k=\sum_{i\geq 1}
x_i^k$:
  \be
   \left\langle \gs_{\lb},
 \gs_{\mu}\right\rangle = \delta_{\lb \mu}.
  \label{orthonormality}\ee
The identity
 \be
\langle it_i \gs_{\lb}, \gs_{\mu}\rangle = \langle
\gs_{\lb}, \frac{\pl}{\pl t_i} \gs_{\mu}\rangle
 \label{transpose}\ee
shows that, with regard to the inner-product
(\ref{inproduct}),
 \be
  \left( \frac{\pl}{\pl t_i}\right)^{\top}
    = it_i,  \label{duality}
  \ee
 and so the matrices representing $it_i$ and
  $ \pl / \pl t_i$ in the orthonormal basis  are
 transpose of each other.
 Also, for $k\geq 0$,
   \bea
    V_k^{\top}&=& \left(\frac{1}{2}\sum_{i+j=k}\frac{\pl^2}{\pl
t_i\pl t_j}+\sum_{-i+j=k}it_i\frac{\pl}{\pl
t_j}+\frac{1}{2}\sum_{-i-j=k}it_i~jt_j \right)^{\top} \no\\
  &=&
   \frac{1}{2}\sum_{i+j=k}
    i t_i jt_j +\sum_{-i+j=k}jt_j\frac{\pl}{\pl
t_i}+\frac{1}{2}\sum_{-i-j=k}\frac{\pl^2}{\pl t_i\pl t_j} \no \\
 &=& V_{-k} ; \label{Virasoroduality}
   \eea
%
%\bea
%  V_k^{\top}&=&V_{-k}; \label{Virasoroduality}
%  \eea
 hence also the matrices associated with $V_k$ and $V_{-k}$ are
 transpose of each other.

%\newpage

We now have:
\begin{proposition}({\bf Murnaghan-Nakayama rule })\cite{Stanley2} We have

 \be
 p_{\al}s_{\lb}=\sum_{\mu}\chi^{\mu\backslash\lb}(\al) s_{\mu}
 \ee
\end{proposition}

\begin{corollary} The following holds:
\bea
it_i~\gs_{\lb}(t)&=&\sum_{{\mu}\atop{\mu\backslash\lb\in
B(i)}}(-1)^{\acht(\mu\backslash\lb)}\gs_{\mu}(t) \no\\
\no\\
\frac{\pl}{\pl
t_i}\gs_{\lb}(t)&=&\sum_{{\mu}\atop{\lb\backslash\mu\in
B(i)}}(-1)^{\acht(\lb\backslash\mu)}\gs_{\mu}(t).
 \label{3.0.14}
\eea

\end{corollary}

\example For a partition $\lb=(\lb_1\geq \ldots \geq
\lb_n\geq 0)$, Corollary 5.3 implies:
  \bean
%V_0(t){\bf s}_{\lb}(t)&=&|\lb|{\bf s}_{\lb}(t)\\
%\\
\frac{\pl}{\pl t_1}{\bf s}_{\lb}(t)&=&\sum_{
         { \nu~\mbox{\tiny such that}}\atop
         { |\lb\backslash\nu|=1  }
          }{\bf s}_{\nu}(t)
          =
          \sum_{{
          1\leq i~\mbox{\tiny such that }}\atop
          {\lb_{i+1}<\lb_{i}
          }
          }
          ~ {\bf s}_{\lb-e_i}(t)
          =\sum_1^n {\bf s}_{\lb-e_i}(t),\\
%\\
%& & \mbox{since $\lb\backslash\mu\in B(1)$ is equivalent to
%$|\lb\backslash\mu|=1$}\\
%\\
t_1{\bf s}_{\lb}(t)&=&\sum_{
         { \nu~\mbox{\tiny such that}}\atop
         { |  \nu\backslash\lb  |=1  }
          }{\bf s}_{\nu}(t)
          =
          \sum_{{
          1\leq i~\mbox{\tiny such that }}\atop
          {\lb_i<\lb_{i-1}
          }
          }
          ~ {\bf s}_{\lb+e_i}(t)
          =\sum_1^{n+1} {\bf s}_{\lb+e_i}(t)\\
\\
\frac{\pl}{\pl t_2}{\bf
s}_{\lb}(t)&=&\sum_{\lb\backslash\nu\in
B(2)}(-1)^{\acht(\lb\backslash\nu)} {\bf s}_{\nu}(t)
 =
 \sum_1^{n} {\bf s}_{\lb-2e_i}(t)
   -  \sum_1^{n-1} {\bf s}_{\lb-e_i-e_{i+1}}(t)
   ,\\
\\
& &\mbox{with $B(2)=\{~\squaresize .4cm \thickness .01cm
\Thickness .07cm \Young{ & \cr
 }~,~
 \squaresize .4cm \thickness .01cm \Thickness .07cm \Young{
 \cr
 \cr
 }~~\}$},
 \eean
 where $e_i=(0,\ldots,0,{\stackrel{i}
 {\stackrel{\downarrow}{1}}}
 ,0,\ldots)^{\top}$. Note that the right most sum in each of the
 expressions means "wherever it makes sense!"; i.e., you
 sum over $i$, wherever $\lb-e_i$, ~$\lb+e_i$,
 $\lb-2e_i$ and $\lb-e_i-e_{i+1}$ are again partitions.

\vspace{1cm}

 \medskip\noindent{\it Proof of Corollary 5.3:\/}   The first identity follows by
applying the Murnaghan-Nakayama rule to $p_i=\sum_{k\geq
1} x_k^i = it_i$. In that case,
\begin{eqnarray*}
\chi^{\mu\backslash \lb} (i)&=& \sum_{ \left\{
\begin{array}{l}
\mbox{\footnotesize border-strip tableaux}~T\\
\mbox{\footnotesize of shape $\mu\backslash \lb $ and
type $i$}
\end{array}
\right\}} (-1)^{{\rm ht} T} \\ &&\\
 &=& (-1)^{{\rm
ht}(\mu\backslash \lb)},~\mbox{with} ~\mu\backslash
\lb\in B(i),
\end{eqnarray*}
since the only border-strip tableau $T$ of shape
$\mu\backslash \lb$ and type $i$ is a border-strip
$\mu\backslash \lb\in B(i)$.

%\newpage

To prove the second relation (\ref{3.0.14}), let
$$
\frac{\pl}{\pl t_{i}} \gs _{\mu} (t) = \sum_{\nu}
c_{\mu\nu} \gs_{\nu} (t).
$$
Then using the duality (\ref{duality}), the result
follows immediately, by taking the transpose of the
first relation (\ref{3.0.14}).\qed

%Then, using orthonormality () of the Schur polynomials with regard
%to $\langle\quad,\quad\rangle$, one computes, using (),
%
%
%\begin{eqnarray*}
%c_{\mu \lb}
%&=& \sum_{\nu} c_{\mu\nu}\langle \gs_{\lb},\gs_{\nu}\rangle\\
%&=& \langle \gs_{\lb}, \frac{\pl}{\pl t_i} \gs_{\mu}\rangle,
%~~\mbox{using ( )}\\
%%
%&=&\langle it_i \gs_{\lb}, \gs_{\mu}\rangle\\
%&=&\left\langle \sum_{\nu\backslash \lb \in B(i)} (-1)^{{\rm
%ht}(\nu\backslash
%\lb)} \gs_{\nu}, \gs_{\mu}\right\rangle, ~~\mbox{using ( )}\\
%&=&\sum_{\nu\backslash \lb \in B(i)} (-1)^{{\rm ht}(\nu\backslash
%\lb)} \langle \gs_{\nu}, \gs_{\mu}\rangle\\ &&\\
%&=&(-1)^{\acht(\mu\backslash \lb)}\langle \gs_{\mu},
% \gs_{\mu}\rangle\\ &&\\
%&=&(-1)^{\acht(\mu\backslash \lb)} \neq 0, \mbox{ iff }
%{\mu\backslash \lb}\in B(i),
%\end{eqnarray*}
%establishing ( ).

%\qed

%\newpage

%\begin{theorem} For $n\geq 0$,
%$$
%V_{-n}\gs_{\lb}=\sum_{\mu\backslash\lb\in
%B(n)}d^{(-n)}_{\lb\mu}\gs_{\mu}.
%$$
%\end{theorem}

\medskip\noindent{\it Proof of Theorem 5.1:\/} At first, it
suffices to prove (\ref{Virasoro10}) for $n\geq 1$; the
second identity (\ref{Virasoro20}) follows immediately
from the duality (\ref{Virasoroduality}).

% Indeed,

\bigbreak

\noindent\underline{\it Step 1}: From (\ref{5.0.2}),
(\ref{orthonormality}), (\ref{transpose})
%Virasoroduality})
  and from Corollary 5.3, it follows that
{\footnotesize  \bean
d_{\lb\mu}^{(-n)}&=&\la V_{-n}\gs_{\lb},\gs_{\mu}\ra \\
\\
&=&\sum_{i\geq 1}\left\la\frac{\pl}{\pl
t_i}\gs_{\lb},\frac{\pl}{\pl t_{n+i}} \gs_{\mu}\right\ra
+ \frac{1}{2} \sum_{i=1}^{n-1}\left\la
it_i\gs_{\lb},\frac{\pl}{\pl t_{n-i}}
 \gs_{\mu}\right\ra
 \\  \\
%\eean \bean
&= &\sum_{i\geq 1}\left\la
\sum_{\lb\backslash\nu\in
B(i)}(-1)^{\acht(\lb\backslash\nu)}\gs_{\nu},\sum_{\mu\backslash\nu'\in
B(i+n)}(-1)^{\acht(\mu\backslash\nu')}\gs_{\nu'}\right\ra \\
\\
& &\quad +\frac{1}{2}\sum_{i=1}^{n-1}
 \left\la \sum_{\nu\backslash\lb\in
B(i)}(-1)^{\acht(\nu\backslash\lb)}\gs_{\nu},\sum_{\mu\backslash\nu'\in
B(n-i)}(-1)^{\acht(\mu\backslash\nu')}\gs_{\nu'}\right\ra \\
\\
  \eean
  \bean
  &=&\sum_{i\geq 1}\sum_{\left\{\begin{array}{l}
           \mbox{\footnotesize $\nu$ such that} \\
           \mbox{\footnotesize $\lb\backslash\nu \in B(i)$} \\
           \mbox{\footnotesize $\mu\backslash\nu\in B(i+n)$}
\end{array}\right\}}(-1)^{\acht(\lb\backslash\nu)+\acht(\mu\backslash\nu)}\\
\\
& &%\quad
+\frac{1}{2}\sum_{i=1}^{n-1}\sum_{\left\{\begin{array}{l}
           \mbox{\footnotesize $\nu$ such that} \\
           \mbox{\footnotesize $\nu\backslash\lb \in B(i)$} \\
           \mbox{\footnotesize $\mu\backslash\nu\in B(n-i)$}
            \end{array}\right\}}
(-1)^{\acht(\nu\backslash\lb)+\acht(\mu\backslash\nu)}
%       \mbox{\footnotesize ht}(\mu \backslash \nu)  }
\eean
}%
Hence
{\footnotesize \bea V_{-n}\gs_{\lb}
             &=&\sum_{\mu}\left\la V_{-n} \gs_{\lb},\gs_{\mu}
               \right\ra \gs_{\mu} \no\\
             &=& \sum_{\mu}\gs_{\mu}\sum_{i\geq 1}
    \sum_{
    \left\{
    \begin{array}{c}
           h_i\in B(i) \\
                                h_{n+i}\in
B(n+i)\\
\mbox{\footnotesize such that}\\
           \lb \backslash h_i=\mu\backslash h_{n+i}\\
           =\mbox{\footnotesize some partition~} \nu
            \end{array}
            \right\}}
 (-1)^{\mbox{\footnotesize ht}(h_i)+\acht(h_{n+i})}\no \\
   && \no \\&&
   + \frac{1}{2} \sum_{\mu}\gs_{\mu}\sum_{i=1}^{n-1}\sum_{
   \left\{
   \begin{array}{c}
           h_i\in B(i) \\
                                h_{n-i}\in
B(n-i)\\
\mbox{\footnotesize such that}\\
           \mu \backslash h_{n-i}=\nu \\
           \nu \backslash h_i=\lb\\
 \mbox{\footnotesize for some partition~} \nu
            \end{array}
         \right\}   }
 (-1)^{\mbox{\footnotesize ht}(h_i) +
       \mbox{\footnotesize ht}(h_{n-i})  }\no\\&&
\label{VirasoroSchur}\eea}

\noindent\underline{\it Step 2}: Having shown that
$$
 V_{-n}\gs_{\lb}=\sum_{{\mu} }
  d_{\lb\mu}^{(-n)} s_{\mu},
  $$
  we now prove that
 $$
 V_{-n}\gs_{\lb}=\sum_{{\mu} \atop {\mu \supset \lb}}
  d_{\lb\mu}^{(-n)} s_{\mu}
.   $$

 %\newpage

The first sum in (\ref{VirasoroSchur}) contains a
summation over $h_i\in B(i)$, $h_{n+i} \in B(n+i)$, such
that $\lb\backslash h_i=\mu\backslash h_{n+i}= $ some
partition $\nu$.  Therefore $\mu= \lb \ba h_i+h_{n+i}$.
We distinguish two cases:
\begin{enumerate}
\item[(i)] $\lb\subset \mu , ~~\mbox{in which case}\quad
h_i\subset h_{i+n}$ \item[(ii)] $\lb\not\subset \mu
\Leftrightarrow h_i \not\subset h_{i+n}$
\end{enumerate}

Assume we are in case (ii).  Then
\begin{eqnarray*}
{\rm either~~ (a)} && h_i \cap h_{i+n}=\emptyset\\
{\rm or ~~(b)}&& h_i \cap h_{i+n}=b\neq \emptyset
\end{eqnarray*}
so in case (ii,b), we have
$$
h_i\backslash b=\alpha \neq \phi, \quad
h_{i+n}\backslash b=\beta\neq \phi$$ and thus $\alpha$
and $\beta$ are border-strips such that
$$
h_i=\alpha+b  ~~\mbox{ and }~~ h_{i+n}=\beta+b.
$$
Therefore, since $\lb \ba h_i=\mu\ba h_{n+i}$, we have
found border strips $\alpha$ and $\beta$, such that
   $$
 \lambda\backslash \alpha=\mu\backslash \beta
  \quad \quad \alpha\in
B(i-|b|), \, \beta\in B(n+i-|b|)
$$
and \be \lambda\backslash (\alpha +b)=\mu\backslash
(\beta + b)=
  \mbox{some partition}~ \nu, %\quad
\quad \alpha +b \in B(i),  ~\, \beta +b\in B(n+i).
 \label{1}\ee
  Hence, in
the first sum of (\ref{VirasoroSchur}), the coefficient
of $\gs_{\mu}$ contains
$$
(-1)^{{\rm ht}(\alpha)+ {\rm ht} (\beta)}+(-1)^{{\rm
ht}(\alpha +b)+ {\rm ht} (\beta +b)} =0,
$$
and therefore does not contribute. The point is that,
if, say, $\beta$ is above $\al$, then the the first
statement in (\ref{1}) forces $\beta$ and $b$ to have
boxes lying on the same row, and
 $\al$ and $b$ on the same column; hence
  \bean
   {\rm ht}(\beta +b)&=&{\rm ht}(\beta)+ {\rm ht}(b)\\
   {\rm ht}(\al +b)&=&{\rm ht}(\al)+ {\rm ht}(b)+1.
 ,\eean
and thus
$$({\rm ht}(\al +b)+ {\rm ht}(\beta +b))
  - ({\rm ht}(\al )+ {\rm ht}(\beta ))
  = 2 {\rm ht}( b)+1,
  $$
  which is odd.
%\newpage

\noindent As an illustration, we give the following
example. Let
$$\lb = \squaresize .4cm \thickness .01cm \Thickness .07cm \Young{
&&&\cr
 &&\cr
 &&*\cr
 &*&*\cr
 &*            \cr
 &*              \cr
 &\shade         \cr
 &\shade\cr }\qquad\qquad\qquad \mu = \squaresize .4cm \thickness .01cm
\Thickness .07cm \Young{
 &&&&\bshade&\bshade\cr
 &&&\bshade&\bshade\cr
 &&*&\bshade\cr
 &*&*\cr
 &*     \cr
 &*       \cr
       \cr
  \cr }
 $$

 \be h_7=
 \squaresize .4cm \thickness .01cm \Thickness .07cm
\Young{ \blank&*\cr *&*\cr *\cr *\cr \shade\cr \shade\cr
 }  ~~\mbox{,}~~~
 h_{10}= \squaresize .4cm \thickness .01cm \Thickness .07cm \Young{
\blank&\blank&\blank&\bshade&\bshade\cr
\blank&\blank&\bshade&\bshade\cr \blank&*&\bshade\cr
*&*\cr *\cr *\cr}
 ~~~\mbox{and}~~~
 h_{7}\cap h_{10}=b=
 \squaresize .4cm \thickness .01cm \Thickness .07cm
\Young{ \blank&*\cr *&*\cr *\cr *\cr %\shade\cr \shade\cr
 }
\label{b}\ee
 $$ \al=h_7-b=h_2=\squaresize .4cm \thickness .01cm \Thickness .07cm
\Young{ \shade\cr \shade\cr
 }~~~\mbox{and}~~~
  \beta=h_{10}-b=h_5=\squaresize .4cm \thickness .01cm \Thickness .07cm \Young{
\blank&\bshade&\bshade\cr \bshade&\bshade\cr
\bshade\cr},
$$
 $$
 \mbox{figure 2}$$
Clearly
 $$
%  \lb -h_{2}=\mu-h_5~~~~\mbox{and}~~~~
\lb-h_7=\mu-h_{10},
 $$
  where $h_7$ and $h_{10}$ are border strips, having a border strip
   $b$ in common ($b$ is given by the boxes containing stars *).
 In
the first sum of (\ref{VirasoroSchur}) above, the
couples of border-strips $(h_2,h_5)$ and $(h_7,h_{10})$
contribute to the sum, as follows:
 \bean
 (-1)^{\acht (\al)+\acht (\beta)}+
 (-1)^{\acht (\al +b)+\acht (\beta +b)}
 &=&
 (-1)^{\acht (h_{2})+\acht (h_5)}+(-1)^{\acht (h_7)+\acht
 (h_{10})}\\&=&(-1)^3+(-1)^{5+5}\\&=&0,
 \eean
 and thus
 \be
 (\acht(h_7)+\acht(h_{10}))- (\acht(h_2)+\acht(h_{5}))
 =2\acht (b)+1 =7 \label{example}
 .\ee
In case (ii,a), we consider $h_i$ and $h_{n+i}~ (h_i\cap
h_{n+i}=\emptyset)$ with the unique connector $b$, which
make $h_i+h_{n+i}+b$ a border strip and again, by the
same reason given before,
  \be (-1)^{{\rm ht}(h_i)+ {\rm
ht} (h_{n+i})}+(-1)^{{\rm ht}(h_i +b)+ {\rm ht} (h_{n+i}
+b)} =0, \label{example2}
  \ee
In Figure 2, consider $h_2$ and $h_5$ with $h_2\cap
h_5=\emptyset$.
The border strip $b$, as in Figure 2, is the unique
connector of $h_2$ and $h_5$. So (\ref{example2}) is
illustrated again by formula (\ref{example}).

We conclude that only case (i) remains, i.e., where
$\lambda \subset \mu$ in which case $h_i \subset
h_{i+n}$.
 Therefore
 \be
 \mu=(\lambda\backslash h_i)+ h_{n+i}=\lambda+ (h_{n+i}\backslash
 h_i), ~~\mbox{with }~ |\mu\backslash \lambda|=|h_{n+i}\backslash h_i|=n%  \in B(n)
 .\label{a}\ee
 In the second summation of (\ref{VirasoroSchur}), it is obvious that
 \be
 \lambda \subset \mu, ~~h_i\cap h_{n-i}=\emptyset ~~\mbox{ and }~~
 |\mu\backslash \lambda|=n% = h_i+h_{n-i} \in B(n)
 .
 \label{b}\ee
So, altogether, we have shown that
 $$
 V_{-n}\gs_{\lb}=\sum_{{{\mu} \atop {\mu \supset \lb}}
                      \atop {|\mu\ba \lb|=n}                         }
  d_{\lb\mu}^{(-n)} s_{\mu}
.   $$

%\newpage

\vspace{.8cm}

\noindent\underline{\it Step 3}:  Finally, we show that

$$
 V_{-n}\gs_{\lb}=\sum_{{{\mu} \atop {\mu \supset \lb}}
                      \atop { \mu\ba \lb \in B(n)}                         }
  d_{\lb\mu}^{(-n)} s_{\mu}
.   $$
 Indeed, if $\mu \backslash \lb  \notin B(n)$, then,
  by (\ref{VirasoroSchur}), we have two
  possible contributions to the coefficient of
  $\gs_{\mu}$; we first deal with the case
  as in (\ref{a}),
  and
 $$
 \lb \subset (\lb \backslash h_i) + h_{i+n}=\mu,~~ \mbox{with}~~
 \mu-\lb \notin B(n)
 $$
   Then $\mu \ba \lb=h_{i+n} \backslash h_i =\sigma$
     is not a border strip, while $h_i+\sigma$ is a
     border strip and hence connected.
  Thus $\sigma$ must come in exactly two border strips; say
  $$h_{i+n} \backslash h_i =\al +\beta,  $$
 and thus we are precisely in the situation of (\ref{b}), with
 {\em either}
 $h_{i'}=\al, ~h_{n-i'}=\beta$, {\em or} $h_{i'}=\beta, ~h_{n-i'}=\al$,
$$
 \lb \subset \mu, ~h_{i'}\cap h_{n-i'}=\emptyset~~\mbox{and} ~~
    \mu \ba \lb=  h_{i'}+ h_{n-i'}.
    $$
    That puts us in the situation of the second
    contribution to $\gs_{\mu}$ in
    (\ref{VirasoroSchur}).

 In Figure 3 below, $\lb$ is the Young diagram to the left of the
bold face line and $\mu$ denotes everything; $\mu
\backslash \lb$ consists of the two border strips $\al$
and $\beta$ below, while $h_{n+i}$ is the large border
strip bordering $\lb$ on the right.
 \bean  &&%~~~~~~~~~
           ~~~~~~~~~ \lb ~~~~~~~~~~~~~~~~~~~\mu\\
        && ~~~~~~~~~ \downarrow
                         ~~~~~~~~~~~~~~~~~\downarrow
 \eean

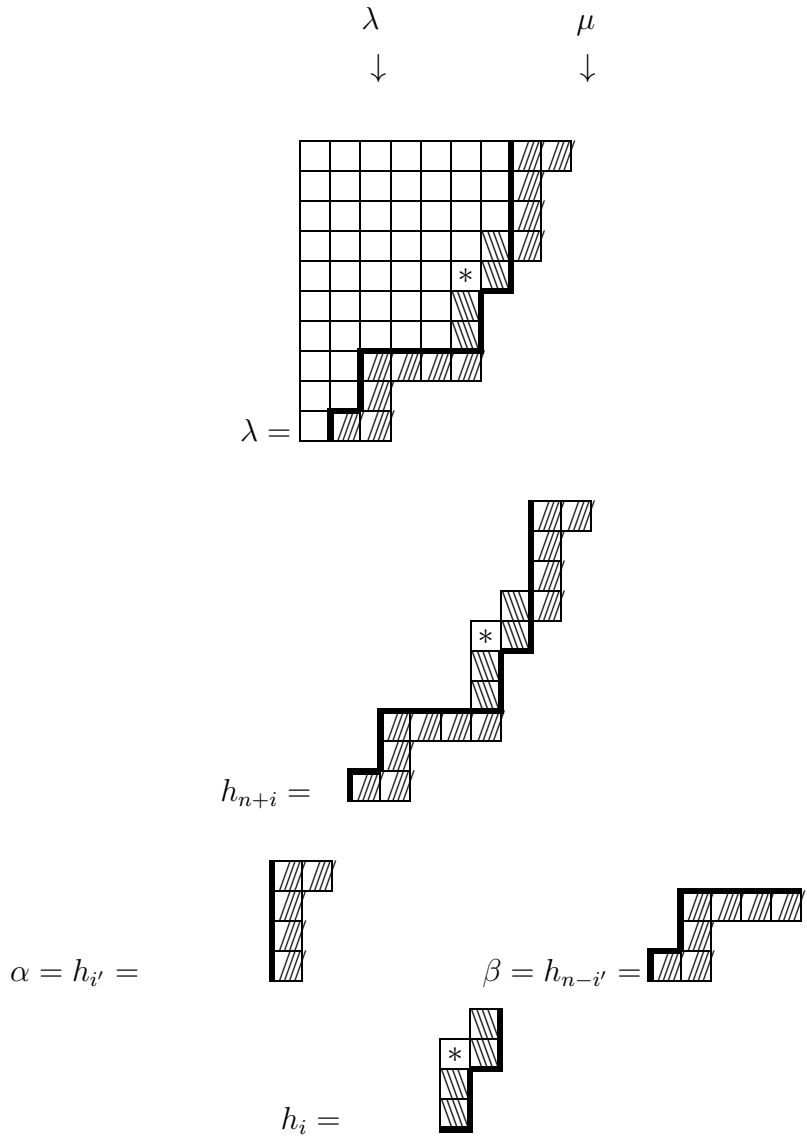
$$
 \lb=\squaresize .4cm \thickness .01cm \Thickness .07cm
  \Young{
 && & & & &  &l\shade&\shade\cr
  &&  & & & & &l\shade\cr
  &&    & & & &  &l\shade\cr
  &&  & & & &\bshade&l\shade\cr
 && & & &*&rd\bshade \cr
 && & & &r\bshade \cr
 && &  &  &rd\bshade  \cr
   & &lu\shade&u\shade&u\shade&\shade  \cr
&  &l\shade\cr
    &ul\shade&\shade \cr }
 $$

$$
 h_{n+i}=\squaresize .4cm \thickness .01cm \Thickness .07cm
  \Young{
 \blank&\blank&\blank &\blank &\blank &\blank &\blank  &l\shade&\shade\cr
  \blank&\blank&\blank  &\blank &\blank &\blank &\blank &l\shade\cr
  \blank&\blank&\blank    &\blank &\blank &\blank &\blank
    &l\shade\cr
  \blank&\blank&\blank  &\blank &\blank &\blank
   &\bshade&l\shade\cr
 \blank&\blank&\blank &\blank &\blank &*&rd\bshade \cr
 \blank& \blank&\blank &\blank &\blank &r\bshade \cr
 \blank& \blank &\blank &\blank  &\blank  &rd\bshade  \cr
   \blank&\blank &lu\shade&u\shade&u\shade&\shade  \cr
\blank& \blank &l\shade\cr
    \blank &ul\shade&\shade \cr }
 $$

 $$
\alpha=h_{i'}=\!\!\!\!\!
 \squaresize .4cm \thickness .01cm \Thickness .07cm \Young{
 \blank&\blank&\blank&\blank&\blank&l\shade&\shade\cr
 \blank&\blank&\blank&\blank&\blank &l\shade\cr
 \blank&\blank&\blank&\blank&\blank &l\shade\cr
 \blank&\blank&\blank&\blank&\blank &l \shade\cr
%\blank&\blank&\blank&r\bshade \cr
 %\blank&\blank&\blank&rd\bshade  \cr
%\blank&lu\shade&u\shade&u\shade \cr
% \blank&l\shade\cr
% ul\shade&\shade
} %%
\hspace{2cm} \beta= h_{n-i'}=
  \squaresize .4cm \thickness .01cm \Thickness .07cm
    \Young{
%\blank&\blank&\blank&\blank&\blank&l\shade&\shade\cr
%\blank&\blank&\blank&\blank&\blank &l\shade\cr
%\blank&\blank&\blank&\blank&\blank &l\shade\cr
%\blank&\blank&\blank&\blank&\bshade&l
%\blank&\blank&\blank&r\bshade \cr
 %\blank&\blank&\blank&rd\bshade  \cr
\blank&lu\shade&u\shade&u\shade &u\shade \cr
\blank&l\shade\cr ul\shade&\shade \cr }
 $$
$$
 h_i
  = \squaresize .4cm \thickness .01cm \Thickness .07cm
  \Young{
%  \blank&\blank&\blank&\blank&\blank&l\shade&\shade\cr
%  \blank&\blank&\blank&\blank&\blank &l\shade\cr
%  \blank&\blank&\blank&\blank&\blank &l\shade\cr
  \blank&\blank&\blank&\blank&r\bshade&\blank\cr
  \blank&\blank&\blank&*&rd\bshade   \cr
 \blank&\blank&\blank&r\bshade \cr
  \blank&\blank&\blank&rd\bshade  \cr
% \blank&lu\shade&u\shade&u\shade \cr
%  \blank&l\shade\cr
%  ul\shade&\shade
 }
 $$
%
% \vspace{-1.2cm}
%
 $$
 \mbox{figure 3}$$
 In the first contribution to $\gs_{\mu}$ in
    (\ref{VirasoroSchur}),
 $$
 \lb\backslash h_i=\mbox{~partition~}=\mu\backslash h_{n+i}
 $$
 contributes \be(-1)^{\acht (h_i)+\acht (h_{i+n})}.
  \label{x}\ee
 In the second contribution to $\gs_{\mu}$ in
    (\ref{VirasoroSchur}).
 $$
 \mu\backslash h_{n-i'}=\nu,\quad\nu\backslash h_{i'}=\lb
 $$
 contributes, upon first setting
  $\al =h_{i'}\quad\beta =h_{n-i'}$, and then setting
  $\al=h_{n-i'},~ \beta = h_{i'},$ yielding
 \be
 2\left(\frac{1}{2}\right)
  (-1)^{\acht
(h_{i'})+\acht (h_{n-i'})}.\label{y}\ee
 The sum of (\ref{x}) and (\ref{y}) vanishes, the point being that
   \bea
   \acht (h_{i})+\acht (h_{i+n})
  &=&  \acht (h_{i})+\acht (h_{\al+h_i+\beta}) \no\\
  &=&  \acht (h_{i})+\acht (h_{\al})
                    +\acht (h_{\beta})
                    +\acht (h_{i}) +1\no\\
   &=&
   \acht (h_{i'}) + \acht (h_{n-i'})
     + 2~\acht (h_{i})+1
   ;\eea
 this is a consequence of the fact that the upper strip $\al$ and
 $h_i$ must have a cell on the same row, where they meet, since
 $\mu \ba(\al +\beta)=\lb$ is a partition.
 For the same reason,
 $\beta$ and $h_i$ cannot have a cell on the same row.
Thus the two contributions to $\gs_{\mu}$ in
    (\ref{VirasoroSchur}) must cancel out in unique pairs,
    if $\mu \ba \lb \notin B(n)$.
     This ends the proof of Theorem 5.1.
    \qed

 %\newpage

 \begin{corollary} For a partition $\lb$
  with at most $n$ rows,
 %$\lb_1^{\top}\leq n$,
  the following holds:

  \bea
  V_{0}{\bf s}_{\lb}
 &=&~~ |\lb| ~{\bf s}_{\lb}   \label{V0} \\&&\no\\
  V_{-1}{\bf s}_{\lb}
 &=&\sum_{{
          1\leq i~\mbox{\tiny such that }}\atop
          {\lb_i<\lb_{i-1}
          }
          }
          (\lb_i-i+1)~ {\bf s}_{\lb+e_i}
          =
          \sum_1^{n+1}
          (\lb_i-i+1)~ {\bf s}_{\lb+e_i}~~~~
           \label{V-1}\\
  V_{1}{\bf s}_{\lb}
 &=&\sum_{{
          1\leq i~\mbox{\tiny such that }}\atop
          {\lb_{i+1}<\lb_{i}
          }
          }
          (\lb_i~-~i)~~ {\bf s}_{\lb-e_i}
          ~~=~
          \sum_{1}^{n}
          (\lb_i-i)~ {\bf s}_{\lb-e_i}
          ~~, \label{V1}
  \eea
where $e_j$ stands for a box added to the right of the
$j$th row of $\lb$. The second sum in (\ref{V-1}) and
(\ref{V1}) refers to a sum over all $i$, with the
understanding that $\gs_{\lb\pm e_i}=0$, if
 $\lb\pm e_i$ is not a Young diagram.
 \end{corollary}

\remark Compare with

 \bea
 t_1{\bf s}_{\lb}&=&
          \sum_1^{n+1} {\bf s}_{\lb+e_i}
          \label{t-action}\\
  \frac{\pl}{\pl t_1}{\bf s}_{\lb}&=&
  \sum_1^n {\bf s}_{\lb-e_i}.
   \label{pl-t-action}
   \\
 \frac{\pl}{\pl t_2}{\bf s}_{\lb}(t) &=&
 \sum_1^{n} {\bf s}_{\lb-2e_i}(t)
   -  \sum_1^{n-1} {\bf s}_{\lb-e_i-e_{i+1}}(t)
   \label{pl-t2-action}
   \eea

 \proof To see identity (\ref{V0}), observe
that for monomials
 \bea
   V_0(\prod_j t_j^{r_j})
  &=& \sum it_i \frac{\pl}{\pl t_i}(\prod_j t_j^{r_j})\no\\
   &=& (\sum jr_j) \prod_j t_j^{r_j}
   \eea
   and so, since Schur polynomials have the form
   $${\bf s}_{\lb}= \sum_{\sum jr_j=|\lb|} \al_r \prod_{j} t_j^{r_j}
   ~,$$
    the result follows.

 We now turn to the identity (\ref{V-1}). From the
identity in Theorem 5.1, it follows that the second sum
in (\ref{VirasoroSchur}) vanishes and so
 \bean
V_{-1}\gs_{\lb}&=& \sum_{\mu\backslash \lb\in B(1)}
 d_{\lb\mu}\gs_{\mu}\no\\
\no\\
&=&\sum_{{ j\geq 1~\mbox{\tiny %\footnotesize
                             such that}}\atop
          { \lb_{j-1}>\lb_j
            }
            }d_{\lb,\lb+e_j}\gs_{\lb+e_j}\no\\
&&\no\\
  \eean
  \bea
  &=&\sum_{{ j\geq 1~\mbox{\tiny%\footnotesize
                              such that}}\atop
           {\lb_{j-1}>\lb_j
            }
            }\gs_{\lb+e_j}\sum_{i\geq 1}
            \sum_{\begin{array}{c}
                        h_i\in B(i) \\
                   h_{i+1}\in B(i+1)\\
           \mbox{\tiny such that}\\
           \lb \backslash h_{i}=(\lb+e_j)\backslash h_{i+1}\\
 =\mbox{\footnotesize some partition~} \nu
            \end{array}}
 (-1)^{\mbox{\footnotesize ht}(h_i) +
       \mbox{\footnotesize ht}(h_{i+1})  }\no\\&& \label{sign}
\eea
In the sum above
$$
h_{i+1}=h_i+\left\{\begin{array}{c}
\mbox{one box, appearing}\\
\mbox{at the right most place}\\
\mbox{in the $j$th row of $h_{i+1}$}
\end{array}\right\}
$$
Consider now a $\lb$ such that $\lb \backslash
h_{i}=(\lb+e_j)\backslash h_{i+1}=\mbox{some
partition}$. Given the rightmost box $b$ in the $j$th
row such that $\lb_{j-1}> \lb_j$; the extra box $e_j$ is
precisely added to the right of this box. Two kinds of
$h_i$, satisfying these requirements can be formed.

$\bullet$ All possible border strips $h$ of $\lb$,
containing that box $b$ and contained in the rows $j,
j+1, \ldots$ of $\lb$, such that $\lb \ba h=\nu$ remains
a partition; see figure 3 for an example. Then $h\in
B(i)$ for some $i\geq 1$. Writing $h_i:=h$ and letting
$h_{i+1}:= h_i+e_j\in B(i)$, we have that $\lb \ba
h_{i}=(\lb+e_j)\ba h_{i+1}$ and $\acht(h_{i+1})=
\acht(h_{i})$.

$\bullet$ All border strips $h$ of $\lb$, containing
that box $b$ and contained in the rows $1,\ldots, j$ of
$\lb$, such that $\lb \ba h=\nu$ remains a partition;
see figure 4 for an example. Then $h\in B(i)$ for some
$i\geq 1$. Writing $h_i:=h$ and letting $h_{i+1}:=
h_i+e_j\in B(i)$, we have that $\lb \ba
h_{i}=(\lb+e_j)\ba h_{i+1}$ and
$\acht(h_{i+1})=\acht(h_{i})+1$.

 \example Consider
 $$
 \lb=
\squaresize .4cm \thickness .01cm \Thickness .07cm
\Young{
 &&&&&& \cr
 &&&&&  \cr
 &&&&&  \cr
 &&&\cr
             \cr
               \cr
                }
 % \right.
 =(7,6,6,4,1,1)
  $$

  $$ \lb+e_j=
\squaresize .4cm \thickness .01cm \Thickness .07cm
\Young{
 &&&&&& \cr
 &&&&&  \cr
 &&&&&  \cr
 &&&&\shade\cr
             \cr
               \cr
                }
 % \right.
 =(7,6,6,5,1,1)
  $$

%\newpage

$\bullet$ On the one hand, the different $h_i$'s,
corresponding to a $+$ sign in (\ref{sign}), are
generated as follows:
 \newline  $
\squaresize .4cm \thickness .01cm \Thickness .07cm
\Young{
 &&&&&& \cr
 &&&&&  \cr
 &&&&&  \cr
 &&&+&\shade\hspace{-.2cm}+\cr
             \cr
               \cr
                }
 % \right.
 $,
  $
\squaresize .4cm \thickness .01cm \Thickness .07cm
\Young{
 &&&&&& \cr
 &&&&&  \cr
 &&&&&  \cr
 &&+&+&\shade\hspace{-.2cm}+\cr
             \cr
               \cr
                }
 % \right.
 $
 ,
 $
\squaresize .4cm \thickness .01cm \Thickness .07cm
\Young{
 &&&&&& \cr
 &&&&&  \cr
 &&&&&  \cr
 &+&+&+&\shade\hspace{-.2cm}+\cr
             \cr
               \cr
                }
 % \right.
 $,
 $
\squaresize .4cm \thickness .01cm \Thickness .07cm
\Young{
 &&&&&& \cr
 &&&&&  \cr
 &&&&&  \cr
 +&+&+&+&\shade\hspace{-.2cm}+\cr
      +       \cr
       +        \cr
              }
 % \right.
 $

 $$\mbox{figure}~ 4$$
\newline They correspond respectively to
$$
\begin{array}{c|c}
\lb \backslash h_i=(\lb+e_j)\backslash h_{i+1}
  &
  (-1)^{\mbox{\footnotesize ht}(h_i) +
       \mbox{\footnotesize ht}(h_{i+1})} \\
\hline
  \lb\backslash
 \squaresize .4cm \thickness .01cm \Thickness .07cm \Young{
 +%&\shade\hspace{-.2cm}+
  \cr
              }=(\lb+e_j)\backslash
 \squaresize .4cm \thickness .01cm \Thickness .07cm \Young{
 + &\shade\hspace{-.2cm}+
  \cr
              }   &
              =1
\\
\lb\backslash
 \squaresize .4cm \thickness .01cm \Thickness .07cm \Young{
 +&+%&\shade\hspace{-.2cm}+
  \cr
              }=(\lb+e_j)\backslash
 \squaresize .4cm \thickness .01cm \Thickness .07cm \Young{
 +&+ &\shade\hspace{-.2cm}+
  \cr
              } & =1\\
\lb\backslash
 \squaresize .4cm \thickness .01cm \Thickness .07cm \Young{
 +&+&+%&\shade\hspace{-.2cm}+
  \cr
              }=(\lb+e_j)\backslash
 \squaresize .4cm \thickness .01cm \Thickness .07cm \Young{
 +&+&+ &\shade\hspace{-.2cm}+
  \cr
              } &=1\\
\lb\backslash
 \squaresize .4cm \thickness .01cm \Thickness .07cm \Young{
 +&+&+&+%&\shade\hspace{-.2cm}+
  \cr
      +       \cr
       +        \cr
              }=(\lb+e_j)\backslash
 \squaresize .4cm \thickness .01cm \Thickness .07cm \Young{
 +&+&+&+ &\shade\hspace{-.2cm}+
  \cr
      +       \cr
       +        \cr
              } &=1
\end{array} $$

\vspace{.7cm}

$\bullet$ On the other hand, the different $h_i$'s,
corresponding to a $-$ sign in (\ref{sign}), are
generated by:
  $$
\squaresize .4cm \thickness .01cm \Thickness .07cm
\Young{
 &&&&&& \cr
 &&&&&  \cr
 &&&&-& - \cr
 &&&&\shade\hspace{-.2cm}-\cr
             \cr
               \cr
                }
 % \right.
 ~~,~~
\squaresize .4cm \thickness .01cm \Thickness .07cm
\Young{
 &&&&&& \cr
 &&&&&-  \cr
 &&&&-& - \cr
 &&&&\shade\hspace{-.2cm}- \cr
             \cr
               \cr
                }
 % \right.
 ~~,~~
\squaresize .4cm \thickness .01cm \Thickness .07cm
\Young{
 &&&&&-&- \cr
 &&&&&-  \cr
 &&&&-& - \cr
 &&&&\shade\hspace{-.2cm}-\cr
             \cr
               \cr
                }
 % \right.
 $$
 $$\mbox{figure}~ 5$$
They correspond respectively to
$$
\begin{array}{c|c}
\lb \backslash h_i=(\lb+e_j)\backslash h_{i+1}
  &
  (-1)^{\mbox{\footnotesize ht}(h_i) +
       \mbox{\footnotesize ht}(h_{i+1})} \\
\hline
\lb \backslash \squaresize .4cm \thickness .01cm
\Thickness .07cm \Young{
 -&- \cr
                }
  =
 (\lb+e_j) \backslash  \squaresize .4cm \thickness .01cm \Thickness .07cm \Young{
  -&- \cr
 \shade\hspace{-.2cm}-\cr
                } & =-1   \\
 \lb \backslash \squaresize .4cm \thickness .01cm \Thickness .07cm
\Young{
  \blank&-\cr
 -&- \cr
                }
  =
 (\lb+e_j) \backslash  \squaresize .4cm \thickness .01cm \Thickness .07cm \Young{
 \blank&-\cr
  -&- \cr
 \shade\hspace{-.2cm}-\cr
                } & =-1   \\
 \lb \backslash \squaresize .4cm \thickness .01cm \Thickness .07cm
\Young{
  \blank&-&-\cr
  \blank&-\cr
 -&- \cr
                }
  =
 (\lb+e_j) \backslash  \squaresize .4cm \thickness .01cm \Thickness .07cm \Young{
 \blank&-&-\cr
  \blank&-\cr
  -&- \cr
 \shade\hspace{-.2cm}-\cr
                } & =-1   \\
 \end{array}
$$

The point is that in the first case each additional $h$
corresponds to adding to the previously added $h$ a new
box to the immediate left of the previous $h$ (in the
$j$th row) and as many boxes as there are in $\lb$ below
the new box added, thus yielding $\lb_j$ cases in all,
each contributing $+1$, since $\acht ~h_i=\acht
~h_{i+1}$. In the second case, we add to the previously
added $h$ one box above the left most box of $h$ and
then all the boxes in $\lb$ to the right of that new
box, yielding $j-1$ cases, each contributing $-1$, since
$\acht ~h_i =\acht~h_{i+1}-1$. So, the total sum equals
$\lb_j-(j-1)=\# \{1\} -\# \{(-1)\}=4-3=1$. For this
example, one computes
 \bean
\lefteqn{V_{-1}{\bf s}_{(7,6,6,4,1,1)}
 }\\ \\
 &=&7~ {\bf s}_{(8,6,6,4,1,1)}
  +5~{\bf s}_{(7,7,6,4,1,1)}
  +  {\bf s}_{(7,6,6,5,1,1)}
  -3~{\bf s}_{(7,6,6,4,2,1)}
  -6~ {\bf s}_{(7,6,6,4,1,1,1)}
  ,\eean
 % [a = 7, b = 5, c = 1, d =  - 3, f = - 6]
corresponding to the different ways to add a box $e_j$
to the partition $\lb$. This ends the proof of
(\ref{V-1}), whereas identity (\ref{V1}) follows from
(\ref{V-1}) by taking the transpose. \qed

 %\newpage

%\newpage

\section{Virasoro actions, Fourier series and difference equations}

Remember the notation (\ref{BVirasoroNotation}) for
\begin{eqnarray}
\BV_{-1}&=&V_{-1}(t)-V_1(s)+n\left(t_1+\frac{\pl}{\pl
s_{1}}\right)\nonumber\\
 \BV_{0}&=&V_0(t)-V_0(s)
 \label{Virasoro}\\
%{\cal D}_{1}I_n
 \BV_{1}&=&-V_{-1}(s)+V_1(t)
  +n\left(s_1+\frac{\pl}{\pl t_1}
\right).\nonumber
\end{eqnarray}
with
 the $V_k$-operators as in (\ref{VirasoroNotation}); e.g.,
 \be
\begin{array}{ll}
V_0(t):=\displaystyle{\sum_{i\geq 1}it_i\frac{\pl}{\pl t_i}}& \\
\\
V_{-1}(t):=\displaystyle{\sum_{i\geq
1}(i+1)t_{i+1}\frac{\pl}{\pl
t_i}},&V_1(t):=\displaystyle{\sum_{i\geq
2}(i-1)t_{i-1}\frac{\pl}{\pl t_i}}.
\end{array} \label{notation}
\ee
In \cite{AvM2}, we have shown the following:
\begin{proposition}
 For all integers $n\geq 0$, the integrals
 \bea
 \tau_n(t,s)&:=& \int_{U(n)}e^{\sum_1^{\iy}\Tr
(t_iM^i-s_i\bar M^i)}dM
 \no\\
&=& %\frac{1}{n!}
 \int_{(S^1)^{n}}|\Dt_n(z)|^{2}
 \prod_{k=1}^n
\left(e^{\sum_1^{\iy}(t_i z_k^i-s_iz_k^{-i})}
 \frac{dz_k}{2\pi i z_k}\right)
 \label{t-s-integral}\eea
satisfy the following three Virasoro constraints:
$$
\BV_i\tau_n=0~~~~\mbox{for}~i=-1,0,1.
$$

\end{proposition}

\vspace{1cm}

Two strictly increasing sets of integers
$x=(x_1<x_2<\ldots<x_n)$ and $y=(y_1<y_2<\ldots<y_n)$ in
$\BZ_{\geq 0}$ are equivalent to two partitions
$\lb=(\lb_1\geq\ldots\lb_n)$ and
 $\mu=(\mu_1\geq\ldots\mu_n)$, with $\lb_1^{\top}\leq n$
  and $\mu_1^{\top}\leq n$, by setting $x_i=\lb_{n-i+1}+(i-1)$ and
   $y_i=\mu_{n-i+1}+(i-1)$, as in (\ref{def(x-lb)}).
%    Written out,
%  \bea
% x&:=&  (0+\lb_n,1+\lb_{n-1},\ldots , n-1+\lb_1)\no\\
% y&:=&  (0+\mu_n,1+\mu_{n-1},\ldots ,n-1+ \mu_1)
%\label{def(x-lb)}\eea

%

\newpage

\subsection{Integral 1}

%\noindent{\bf Example 6.1}:
 {\em \noindent {\bf (i)} For
the fixed partitions $\lb$ and $\mu$, with at most $n$
rows, the integral
 \be
   a_{\lb\mu}(z):=\int_{U(n)} {\bf s}_{\lb}(M){\bf s}_{\mu}
    (\bar M)
  e^{z\Tr (M+\bar M)}dM=
  \sum_0^{\iy}   \frac{z^k}{k!}b^{(k)}_{xy}
  \label{correspondence}
 \ee
satisfies the difference equations, with $ {\cal
L}^{(1)}_{\pm \atop 0}:={\cal L}^{(1)}_{\pm \atop
0}(\lb,\mu,z)$,
\bea
  {\cal L}_-^{(1)}(a_{\lb\mu})
   &:=&
   \sum_{  1\leq i\leq n   }\Bigl(a_{\lb-e_i,\mu}(\lb_i-i+n)
   -a_{\lb,\mu+e_i}  (\mu_i-i+n+1)\Bigr)
 \no\\
 &&
  ~~+z\left(n~a_{\lb\mu}-\sum_1^n a_{\lb,\mu+2e_i}+
                        \sum_1^n a_{\lb,\mu+e_i+e_{i+1}}
                                               \right)=0
 \no\\
  &&
 \label{a-difference2}
  \\ \no  \\
% \eea
%
%
% \bea
{\cal L}_0^{(1)}(a_{\lb\mu}) &:=&
\left(|\lb|-|\mu|\right)~a_{\lb,\mu}
  +z\sum_{ 1\leq i\leq n  }
  \left(    a_{\lb+e_i,\mu}
  -
            a_{\lb,\mu+e_i}
    \right) =0
  \label{a-difference1}
  \\
  \no\\
  {\cal L}_+^{(1)}(a_{\lb\mu})&=&0,
  \label{a-difference1'}\eea
  where $$ {\cal L}_+^{(1)}(\lb,\mu):=
           {\cal L}_-^{(1)}(\mu,\lb)%\Bigr|_{\lb \leftrightarrow\mu}
           .$$

\noindent{\bf (ii)} Moreover,
  \be
 b_{xy}^{(k)}:= \#\left\{\begin{array}{l}
     \mbox{ways that $n$ non-intersecting walkers}\\
     %\mbox{(at each step only one walker walks)}\\
     \mbox{%by $n$ walkers
     in $\BZ$ move from $x_1<x_2<...<x_n$
     }\\
     \mbox{to $y_1<y_2<...<y_n$ in $k$ steps
     }\\
     \end{array}\right\}
     \label{b's}\ee
 satisfies the difference equations, alluded to in
 table 2 of the introduction, with $ {
L}^{(1)}_{\pm \atop 0}:={ L}^{(1)}_{\pm \atop
0}(x,y,\Lb_k)$,

 \bea
   L^{(1)}_{-} \left(b^{(k)}_{xy}\right)
    &:=&
 -
 \sum_{1\leq i \leq n} \left(
  (y_i+1) b^{(k)}_{x,y+e_i}-x_i b^{(k)}_{x-e_i,y}
  \right)
  \no\\
  &&
  +k\Bigl( nb^{(k-1)}_{x,y}-
  \sum_{1\leq i \leq n}   b^{(k-1)}_{x,y+2e_i}
  +
  \sum_{1\leq i\leq n\atop{y_{i+1}-y_i=1}}
%  \sum_{1\leq i \leq n}\delta_{y_{i+1},y_i+1}
 b^{(k-1)}_{x,y+e_i+e_{i+1}}\Bigr)=0 \no\\&&
  \label{b-difference2}
%\\ &&\no\\&&\no\\ \no
 \eea
  \bea
   L^{(1)}_{0}\left(b^{(k)}_{xy}\right)
 &:=& -b_{xy}^{(k)}\sum_{1\leq i \leq n}  (y_i-x_i)
 +k%\frac{k}{\sum_1^n (y_i-x_i)}
 \left(
  \sum_{1\leq i \leq n}b_{x+e_i,y}^{(k-1)}
  -
  \sum_{1\leq i \leq n}b_{x,y+e_i}^{(k-1)}
 \right)=0\no\\&&
 \label{b-difference1}
 \\&&\no\\&&\no\\
 %  \eean
 %  \bea
 L^{(1)}_{+} \left(b^{(k)}_{xy}\right) &=& 0,~~~
 \mbox{where}~~  L^{(1)}_{+}(x,y):=
  L^{(1)}_{-}(y,x).\label{b-difference0}
 \eea

 }

\remark In the formulae above, the integer
$b^{(k)}_{x,y}=0$, when the strict inequalities
$x_1<x_2<...<x_n$ and $y_1<y_2<...<y_n$ are not
satisfied.

 \proof Applying the shifts
  $t_1\longmapsto t_1+z,~~s_1\longmapsto s_1-z$ to the
   Virasoro constraints of Proposition 6.1,
   and to the matrix integral (\ref{t-s-integral}) lead
   to the following equations, for $k=-1,0,1$,

\bean  0&=&\BV_k
 \int_{U(n)}e^{z\Tr(M+\bar M)}
 e^{\sum_1^{\iy}\Tr(t_iM^i-s_i\bar M^{i})}dM
   \\
   &=&
    \BV_k
    \sum_{{\lb,\mu\mbox{\tiny ~such
that}}\atop{\lb_1^{\top},\mu_1^{\top}\leq
n}}a_{\lb\mu}(z)s_{\lb}(t)s_{\mu}(-s)
  ,\eean
where
\begin{eqnarray}
\BV_{-1}
 &=&
    V_{-1}(t)-V_1(-s)
 +n\left(t_1+\frac{\pl}{\pl s_1}\right)
  +z\left(n+\frac{\pl} {\pl s_2}\right)
   \no\\\no\\
   \BV_{0}
 &=&
  V_0(t)-V_0(-s)+z\left(\frac{\pl}{\pl
t_1}+\frac{\pl}{\pl s_1}\right)\label{Virasoro1}
  \\\no\\
   \BV_{1}
    &=&
   -V_{-1}(s)+V_1(-t)
 +n\left(s_1+\frac{\pl}{\pl t_1}\right)-
  z\left(n-\frac{\pl}
{\pl t_2}\right).\no
\end{eqnarray}
%
%
%
%\newpage
%
Using the action (\ref{t-action}),(\ref{pl-t2-action}),
(\ref{V-1}), (\ref{V1}) of $t_1$, $\pl / \pl t_2$,
$V_{-1}$ and $V_{1}$ on Schur polynomials, one computes
   \bean
%\mathbb V_{-1}(t,s)&=&V_{-1}(t)-V_1(-s)
% +n\left(t_1-\frac{\pl}{\pl(-s_1)}\right)+z\left(n-\frac{\pl}
%{\pl(-s_2)}\right)\\
0&=& \mathbb
V_{-1}(t,s)\sum_{\lb,\mu\atop{\lb^{\top}_1,\mu^{\top}_1\leq
n}}a_{\lb\mu} (z) \gs_{\lb}(t)\gs_{\mu}(-s)%\tau_n
 \\
 &=&\sum_{\lb',\mu'\atop{\lb^{'\top}_1,\mu^{'\top}_1\leq
n}}a_{\lb'\mu'}\left\{\begin{array}{l}
(V_{-1}(t)\gs_{\lb'}(t))\gs_{\mu'}(-s)
 -\gs_{\lb'}(t)(V_1(-s)\gs_{\mu'}(-s))\\  \\
+n(t_1\gs_{\lb'}(t))\gs_{\mu'}(-s)-n\gs_{\lb'}(t)
 \left(\frac{\pl}{\pl(-s_1)}\gs_{\mu'}(-s)\right)\\ \\
+z~n~\gs_{\lb'}(t)\gs_{\mu'}(-s)-z~\gs_{\lb'}(t)\frac{\pl}{\pl(-s_2)}\gs_{\mu'}(-s)
\end{array}
\right\}\\
\\&&\\
&=&\sum_{\lb',\mu'\atop{\lb^{'\top}_1,\mu^{'\top}_1\leq
n}}a_{\lb'\mu'}\left\{\begin{array}{l} \displaystyle{
\sum_1^{n+1}}(\lb'_i-i+1)
 \gs_{\lb'+e_i}(t)\gs_{\mu'}(-s)  \\
 -\gs_{\lb'}(t)
  \displaystyle{\sum_1^{n}}
   (\mu'_i-i)\gs_{\mu'-e_i}(-s)
   \\
   +n\left(\displaystyle{\sum_{1}^{n+1}}
   \gs_{\lb'+e_i}(t)\gs_{\mu'}(-s)-
\gs_{\lb'}(t) \displaystyle{\sum_{1}^n}
\gs_{\mu'-e_i}(-s)\right)\\
+z n~\gs_{\lb'}(t)\gs_{\mu'}(-s)\\
 -z\gs_{\lb'}(t)
\left(\displaystyle{\sum_1^{n}} {\bf s}_{\mu'-2e_i}(-s)
   -  \displaystyle{\sum_1^{n-1}}
      {\bf s}_{\mu'-e_i-e_{i+1}}(-s)  \right)
%   @@@@@@@@@\displaystyle{\sum_{\nu\atop{\mu'
%\backslash \nu\in
%B(2)}}}
%@@\gs_{\nu}(-s)(-1)^{\acht(\mu'\backslash\nu)}
%
\end{array}
\right\}(*)\\
& &\\
&=&\sum_{\lb^{\top}_1\leq n\atop{\mu^{\top}_1\leq
n}}\gs_{\lb}(t)\gs_{\mu}(-s){\cal
L}_-^{(1)}(a_{\lb\mu}(z))~~~~~~~~~~(**)
 . \eean
The coefficients ${\cal L}_-^{(1)}(a_{\lb\mu}(z))$ of
$\gs_{\lb}(t)\gs_{\mu}(-s)$ in the formula above are
linear expressions in the $a_{\lb\mu}(z)$, given by
 the difference equation (\ref{a-difference2}). Also
 notice that,
 from the expression (*), it would seem one would have
 to include
  in the sum (**) partitions
  $\lb$ of the form $\lb=\lb'+e_{n+1}$,
  with $\lb_1^{\prime\top}=n$, so that
  $\lb_1^{\top}=n+1$. However, the coefficient
  of such terms vanish by visual inspection!

The Virasoro duality
$$\BV_{1}=-\BV_{-1}\Bigl|_{t\longleftrightarrow -s}$$
implies at once expression (\ref{a-difference1'}).

Finally, \bean   0&=& \mathbb
V_0\sum_{\lb,\mu\atop{\lb^{\top}_1,\mu^{\top}_1\leq
n}}a_{\lb\mu} (z) \gs_{\lb}(t)\gs_{\mu}(-s)\\
&=&\sum_{\lb,\mu\atop{\lb^{\top}_1,\mu^{\top}_1\leq
n}}a_{\lb\mu}\left\{
\begin{array}{l}
(V_0(t)\gs_{\lb}(t))\gs_{\mu}(-s)-\gs_{\lb}(t)(V_0(-s)\gs_{\mu}(-s))\\
\\
+z\left(\left(\frac{\pl}{\pl
t_1}\gs_{\lb}(t)\right)\gs_{\mu}(-s)-\gs_{\lb}(t)\left(
\frac{\pl}{\pl (-s_1)}\gs_{\mu}(-s)\right)\right)
\end{array}\right\}\\
\\
&=&\sum_{\lb,\mu\atop{\lb^{\top}_1,\mu^{\top}_1\leq
n}}a_{\lb\mu}\left\{
\begin{array}{l}
(|\lb|-|\mu|)\gs_{\lb}(t)\gs_{\mu}(-s)\\
\\
+z\left(\displaystyle{\sum_{1}^n}
\gs_{\lb-e_i}(t)\right)\gs_{\mu}(-s)-z~\gs_{\lb}(t)
 \left(\displaystyle{\sum_1^{n}}
 \gs_{\mu-e_i}(-s)\right)
\end{array}\right\}\\
& &\\
&=&\sum_{\lb^{\top}_1\leq n\atop{\mu^{\top}_1\leq
n}}\gs_{\lb}(t)\gs_{\mu}(-s)
 {\cal L}_0^{(1)}(a_{\lb\mu}(z))
 .
\eean
In this expression, the coefficient ${\cal
L}_0^{(1)}(a_{\lb\mu}(z))$ of
 $
\gs_{\lb}(t)\gs_{\mu}(-s)$
 is precisely equation (\ref{a-difference1}).

Then setting
$$
a_{\lb\mu}(z)=\sum_{k\geq 0}b_{xy}^{(k)}\frac{z^k}{k!}
$$
in (\ref{a-difference2}), (\ref{a-difference1}) and
(\ref{a-difference1'}), one finds, using the map
$x_{n+1-i}=\lb_i+n-i$, ~$y_{n+1-i}=\mu_i+n-i$,
 \be
  0= {\cal L}_{\pm \atop 0}^{(1)}(a_{\lb\mu}(z))
  ={\cal L}_{\pm \atop 0}^{(1)}
   \left(\sum_{k\geq 0}b_{xy}^{(k)}\frac{z^k}{k!}\right)
  =\sum_{k\geq 0}  L^{(1)}_{\pm \atop 0}\left(b_{xy}^{(k)}\right)
     \frac{z^k}{k!}, \label{L-argument}
  \ee
leading to (\ref{b-difference2}), (\ref{b-difference1})
and (\ref{b-difference0}).
\qed

%\newpage

\subsection{Integral 2}

%  {\bf Example 6.2}:

 {\em \noindent {\bf (i)} For the fixed partitions $\lb$ and
$\mu$, with at most $n$ rows,
 the integral
 \be
   a_{\lb\mu}(z):=\int_{U(n)} {\bf s}_{\lb}(M)
   {\bf s}_{\mu}(\bar M)
  \det(I+M)^q
  e^{z \Tr\bar M}dM=\sum^{\iy}_0
b_{xy}^{(k)}\frac{z^k}{k!}
  \label{correspondence}\ee
satisfies the difference equations, with $ {\cal
L}^{(2)}_{\pm }:={\cal L}^{(2)}_{\pm}(\lb,\mu,z)$,
\bea
  {\cal L}_-^{(2)} (a_{\lb\mu})
   &:=&
  \sum_1^{n}a_{\lb-e_i,\mu}(\lb_i-i+n)
 -\sum_1^n a_{\lb,\mu+e_i}(\mu_i-i+n+z+1)\no\\
\no\\
& +&\!\!\!a_{\lb\mu}(|\lb|-|\mu|+nq)
-z \Bigl(
  \sum_1^n a_{\lb,\mu+2e_i}-\sum_1^{n-1}
  a_{\lb,\mu+e_i+e_{i+1}}\Bigr)=0\no\\
 \label{example2-a1}
 \\
% \lefteqn{\hspace{-1.1cm} \mbox{and}}\no\\
 %
 %
  {\cal L}_+^{(2)} (a_{\lb\mu})
   &:=&\sum_1^{n}a_{\lb,\mu-e_i}
 (\mu_i-i+n)-\sum_1^n a_{\lb+e_i,\mu}(\lb_i-i+n+q+1)\no\\
   \no\\
& &+~a_{\lb\mu}(|\mu|-|\lb|+nz)
+z
  \sum_1^n a_{\lb,\mu+e_i}=0
%\sum_{{{\mu'
%}\atop{{{\mu'\ba\mu\in B(2)}\atop{\mu^{'\top}_1\leq
%n}}}} }a_{\lb,\mu'}(-1)^{\acht(\mu'\ba\mu)}
. \no\\
\label{example2-a2}\eea
\noindent {\bf (ii)} The coefficients $b_{xy}^{(k)}$
satisfy, with $ { L}^{(2)}_{\pm }:=
   L^{(2)}_{\pm}(x,y,\Lb_k)$,
 \bea
 L^{(2)}_{-}\left(b^{(k)}_{xy}\right)
   &:=&
  \sum_1^n\left(x_i b_{x-e_i,y}^{(k)}
-(y_i+1)b^{(k)}_{x,y+e_i}
+(x_i-y_i+q)b_{xy}^{(k)}\right)
   \no\\
& &~
 +k\left(-\sum_1^n(
   b_{x,y+e_i}^{(k-1)}
  - b_{x,y+2e_i}^{(k-1)}) +
   \sum_1^{n-1} b_{x,y+e_i+e_{i+1}}^{(k-1)}  \right)
  =0 \no\\
  \label{example2-b1}
  \\
  &&
  \no\\
  &&\no\\
    L^{(2)}_{+}  \left(b^{(k)}_{xy}\right)
   &:=&
%
%\lefteqn{
\sum_1^n\left(y_i b_{x,y-e_i}^{(k)}
 -(x_i+1+q)b^{(k)}_{x+e_i,y}+(y_i-x_i)b_{xy}^{(k)}\right)
% }
  \no \\
   & &~
 +k
  \sum_1^{n}\left(b_{xy}^{(k-1)}+
      b_{x,y+e_i}^{(k-1)}  \right)
  =0.
  \label{example2-b2}
  \eea
}

 \proof Applying the shifts
  $it_i\longmapsto it_i-q(-1)^i$, $is_i\longmapsto is_i-z\delta_{i1}$ to the
   Virasoro constraints and the
   matrix integral of Proposition 6.1,
   leads to the equations
\bean  0&=&\BV_k
 \int_{U(n)}\det(I+M)^q
  e^{z \Tr\bar M}
   e^{\sum_1^{\iy}\Tr(t_iM^i-s_i\bar M^{i})}dM\\
   \\
   &=&
    \BV_k
    \sum_{{\lb,\mu\mbox{\tiny ~such
that}}\atop{\lb_1^{\top},\mu_1^{\top}\leq
n}}a_{\lb\mu}(z)s_{\lb}(t)s_{\mu}(-s)
  ,\eean
where
%
% {\bf Virasoro}:
\begin{eqnarray}
\BV_{-1} &=&
    V_{-1}(t)-V_1(-s)
 +n\left(t_1+q+\frac{\pl}{\pl s_1}\right)
  +q\sum_{i\geq 1}(-1)^i\frac{\pl} {\pl t_i }+z\frac{\pl} {\pl s_2 } \no\\\no\\
   \BV_{0}
        &=&
  V_0(t)-V_0(-s)-q\sum_{i\geq 1}(-1)^i\frac{\pl}{\pl t_i}
      +z\frac{\pl}{\pl s_1}
         \label{Virasoro1}
  \\\no\\
   \BV_{1}
        &=&
   -V_{-1}(s)+V_1(-t)
 +n\left(s_1-z+\frac{\pl}{\pl t_1}\right)
   +q\sum_{i\geq 2}(-1)^i\frac{\pl}{\pl t_i}.\no
\end{eqnarray}
The only linear combinations in the span of
$\BV_{-1},\BV_{0},\BV_{1}$, involving finite sums of
$V_k(t),V_k(s), \pl /\pl t_k, \pl /\pl s_k, t_k, s_k$
are as follows:
\bean \mathbb V_{-1}+\mathbb V_{0}
 &=&\Bigl(V_{-1}(t)+nt_1\Bigr)
  -\Bigl(V_1(-s)+(n+z)\frac{\pl}{\pl
(-s_1)}\Bigr)\\
&&~~~~+V_0(t)-V_0(-s)+nq-z\frac{\pl}{\pl
(-s_2)}\\
 \\
-\mathbb V_{0}-\mathbb V_{1}
 &=&
  \Bigl(V_{-1}(-s)+n(-s_1)\Bigr)
   - \Bigl(V_1(t)+(n+q)\frac{\pl}{\pl t_1}\Bigr)\\
   &&~~~~ +V_0(-s)-V_0(t)+nz+z\frac{\pl}{\pl
(-s_1)} \\
%&&~~~~+(n+q)\frac{\pl}{\pl t_1}-n(-s_1-z)
 \eean
and thus, using Corollary 5.3 and 5.4, compute
 \bean \lefteqn{0= (\mathbb V_{-1}+\mathbb
V_0)\sum_{{\lb,\mu\mbox{\tiny ~such
that}}\atop{\lb_1^{\top},\mu_1^{\top}\leq
n}}a_{\lb\mu}(z)s_{\lb}(t)s_{\mu}(-s)}
 \\
 &=&\sum_{\lb',\mu'\atop{\lb^{'\top}_1,\mu^{\top}_1\leq
n}}a_{\lb'\mu'}\left\{\begin{array}{l}
(V_{-1}(t)+nt_1)\gs_{\lb'}(t)\gs_{\mu'}(-s)
    \\  \\
    -\gs_{\lb'}(t)\left(V_1(-s)+(n+z)
  \frac{\pl}{\pl(-s_1)}\right)\gs_{\mu'}(-s)\\ \\
+(V_0(t)\gs_{\lb'}(t))\gs_{\mu'}(-s)-\gs_{\lb'}(t)
 (V_0(-s)\gs_{\mu'}(-s))\\ \\
+nq~\gs_{\lb'}(t)\gs_{\mu'}(-s)-z~\gs_{\lb'}(t)\frac{\pl}{\pl(-s_2)}\gs_{\mu'}(-s)
\end{array}
\right\}\\
&&\\
 &=&\sum_{\lb',\mu'\atop{\lb^{'\top}_1,\mu^{'\top}_1\leq
n}}a_{\lb'\mu'}\left\{\begin{array}{l}
\displaystyle{\sum_1^{n+1}}(\lb'_i-i+1+n)\gs_{\lb'+e_i}(t)\gs_{\mu'}(-s)\\
\\
-\gs_{\lb'}(t)\displaystyle{\sum_1^n}(\mu'_i-i+n+z)\gs_{\mu'-e_i}(-s)\\
\\
+(|\lb'|-|\mu'|+nq)\gs_{\lb'}(t)\gs_{\mu'}(-s)\\
\\
-z~\gs_{\lb'}(t) %\\
 \displaystyle{(\sum_{1}^{n}{\bf s}_{\mu'-2e_i}(-s)
        -\sum_{1}^{n-1}{\bf s}_{\mu'-e_i-e_{i+1}}(-s)})
        %(-1)^{\acht(\mu'\ba\nu)}\gs_{\nu}(-s)
\end{array}
\right\}\\
& &\\
&=:&\sum_{\lb^{\top}_1\leq n\atop{\mu^{\top}_1\leq
n}}\gs_{\lb}(t)\gs_{\mu}(-s)
 {\cal L}_-^{(2)}(a_{\lb\mu}(z))
 ,\eean
implying the vanishing of all the coefficients ${\cal
L}_-^{(2)}(a_{\lb\mu}(z))$ of
$\gs_{\lb}(t)\gs_{\mu}(-s)$, leading to the difference
equation (\ref{example2-a1}). The same remark as for
integral 1 holds for this case.

The remaining identities are obtained in a similar
fashion:
 \bean \lefteqn{0= -(\mathbb V_{1}+\mathbb
V_0)\tau_n}
 \\
 &=&\sum_{\lb',\mu'\atop{\lb^{'\top}_1,\mu^{\top}_1\leq
n}}a_{\lb'\mu'}\left\{\begin{array}{l}
 \gs_{\lb'}(t)(V_{-1}(-s)-ns_1)\gs_{\mu'}(-s)
    \\  \\
    -\left(V_1(t)+(n+q)
  \frac{\pl}{\pl t_1}\right)\gs_{\lb'}(t)\gs_{\mu'}(-s)
   \\  \\
+(\gs_{\lb'}(t))V_0(-s)\gs_{\mu'}(-s)-V_0(t)\gs_{\lb'}(t)
 \gs_{\mu'}(-s)\\ \\
+nz~\gs_{\lb'}(t)\gs_{\mu'}(-s)
 +z~\gs_{\lb'}(t)\frac{\pl}{\pl(-s_1)}\gs_{\mu'}(-s)
\end{array}
\right\}\\
&&\\
 &=&\sum_{\lb',\mu'\atop{\lb^{'\top}_1,\mu^{'\top}_1\leq
n}}a_{\lb'\mu'}\left\{\begin{array}{l}
\displaystyle{\sum_1^{n+1}}\gs_{\lb'}(t)
  (\mu'_i-i+1+n)\gs_{\mu'+e_i}(-s)\\
\\
-\displaystyle{\sum_1^n}(\lb'_i-i+n+q)\gs_{\lb'-e_i}(t)
 \gs_{\mu'}(-s)\\
\\
+(|\mu'|-|\lb'|+nz)\gs_{\lb'}(t)\gs_{\mu'}(-s)\\
\\
   +z~\gs_{\lb'}(t) %\\
 \displaystyle{ \sum_{1}^{n}{\bf s}_{\mu'-e_i}(-s)
        )}
        %(-1)^{\acht(\mu'\ba\nu)}\gs_{\nu}(-s)
\end{array}
\right\}\\
& &\\
&=&\sum_{\lb^{\top}_1\leq n+1\atop{\mu^{\top}_1\leq
n}}\gs_{\lb}(t)\gs_{\mu}(-s){\cal
L}_+^{(2)}(a_{\lb\mu}(z))
 \eean
implying ${\cal L}_+^{(2)}(a_{\lb\mu}(z))=0$, and thus
(\ref{example2-a2}).
Identities (\ref{example2-b1}) and (\ref{example2-b2})
follow by the precise same method as for integral 1.\qed

%\newpage

\subsection{Integral 3
 }

%  {\bf Example 6.3}:
{\em \noindent {\bf (i)} For the fixed partitions $\lb$
and $\mu$, with at most $n$ rows,
 the integral
 \be
   a_{\lb\mu}(z):=\int_{U(n)}
   s_{\lb}(M)s_{\mu}(\bar M)
   \det (I+z M)^p \det (I+z\bar M)^q dM
=\sum_0^{\iy} b_{xy}^{(k)}{z^k}
  \label{correspondence}\ee
satisfies the difference equation, with $ {\cal
L}^{(3)}_{}:={\cal L}^{(3)}_{}(\lb,\mu,z)$,
  \bea
   {\cal L}^{(3)} (a_{\lb\mu})
   &=&\sum^n_{1}\left(a_{\lb-e_i,\mu}(\lb_i-i+n)+a_{\lb+e_i,\mu}
(\lb_i-i+n+p+1)\right)\no\\
&
&-\sum_{1}^n\left(a_{\lb,\mu-e_i}(\mu_i-i+n)+a_{\lb,\mu+e_i}
(\mu_i-i+n+q+1)\right)\no\\
& &+(z+z^{-1})(|\lb|-|\mu|)a_{\lb\mu}=0.
  \label{ExampleThreea}\eea
\noindent {\bf (ii)} The coefficients $b_{xy}^{(k)}$
satisfy, with $ { L}^{(3)}_{}:=L^{(3)}_{}(x,y,\Lb_k)$,
\bea
   L^{(3)} \left(b^{(k)}_{xy}\right)
   &= &
    \sum^n_1 \left(\begin{array}{l}
   \left(x_i b^{(k)}_{x-e_i,y}
   +(x_i+p+1)b^{(k)}_{x+e_i,y} \right)\\
   \\
   -\left(y_i b^{(k)}_{x,y-e_i}
       +(y_i+q+1)b^{(k)}_{x,y+e_i}\right)\end{array}\right)\no\\
\no\\
& &\quad
+\sum_1^n(x_i-y_i)\left(b_{xy}^{(k-1)}+b_{xy}^{(k+1)}\right)=0.
\label{ExampleThreeb}
  \eea
}

\proof
 Applying the shifts
  $it_i\longmapsto it_i-p(-z)^i$, $is_i\longmapsto is_i+q(-z)^i
   $ to the
   Virasoro constraints and the matrix integral of Proposition 6.1,
   leads to the equations
    \bean  0&=&\BV_k
 \int_{U(n)}\det (I+z M)^p \det (I+z\bar M)^q
   e^{\sum_1^{\iy}\Tr(t_iM^i-s_i\bar M^{i})}dM\\
   \\
   &=&
    \BV_k
    \sum_{{\lb,\mu\mbox{\tiny ~such
that}}\atop{\lb_1^{\top},\mu_1^{\top}\leq
n}}a_{\lb\mu}(z)s_{\lb}(t)s_{\mu}(-s)
  ,\eean
where
% {\bf Virasoro}:
{\footnotesize
\begin{eqnarray} \BV_{-1}
  &=&
    V_{-1}(t)-V_1(-s)
 +n\left(t_1-\frac{\pl}{\pl (-s_1)}\right)
  \no\\
  &&  ~~~~~~~~~~~~~~~~~~~~~~~~~~~~~~~~~
  -p\sum_{i\geq 1}(-z)^{i+1}\frac{\pl} {\pl t_i }-q
  \sum_{i\geq 2}(-z)^{i-1}\frac{\pl} {\pl s_i }
     \no\\\no\\
   \BV_{0}
        &=&
  V_0(t)-V_0(-s)-p\sum_{i\geq 1}(-z)^i\frac{\pl}{\pl t_i}
   -q\sum_{i\geq 1}(-z)^i\frac{\pl}{\pl s_i}\label{Virasoro1}
  \\\no\\
   \BV_{1}
        &=&
   -V_{-1}(-s)+V_1(t)
 +n\left(s_1+\frac{\pl}{\pl t_1}\right)
  \no\\
  &&  ~~~~~~~~~~~~~~~~~~~~~~~~~~~~~~~~~
  -q\sum_{i\geq 1}(-z)^{i+1}\frac{\pl}
{\pl s_i}
  -p\sum_{i\geq 2}(-z)^{i-1}\frac{\pl}
{\pl t_i}.\no
\end{eqnarray}
}
Here the only linear combination in the span of
$\BV_{-1},\BV_{0},\BV_{1}$, involving finite sums of
$V_k(t),V_k(s), \pl /\pl t_k, \pl /\pl s_k, t_k, s_k$ is
the following expression:
\bean \lefteqn{ \mathbb V_{-1}+(z+z^{-1})\mathbb
V_{0}+\mathbb V_{1}   }
 \\
 &=&\left(V_{-1} (t)+V_{1} (t)+n\left(t_{1} +\frac{\pl}{\pl t_{1}
 }\right)\right)\\
 & &-\left(V_{-1} (-s)+V_{1} (-s)+n\left((-s_{1}) +\frac{\pl}{\pl (-s_{1} )
 }\right)\right)\\
 & &+(z+z^{-1})(V_{0}(t)-V_{0}(-s))+\left(p\frac{\pl}{\pl t_{1}}-q
 \frac{\pl}{\pl (-s_{1})}\right),
\eean
and so

{\footnotesize
\bean \lefteqn{0= \left(\mathbb V_{-1}+(z+z^{-1})\mathbb
V_{0}+\mathbb V_{1} \right)\tau_{n}  }
 \\
 \\
 &=&\sum_{\lb',\mu'\atop{\lb_1^{'\top},\mu_1^{'\top}\leq n}}
   a_{\lb'\mu'}
\left\{
\begin{array}{l}
\left(
\begin{array}{l}
V_{-1}(t)+V_1(t)+(z+z^{-1})V_0(t)\\
+~nt_1+(n+p)\frac{\pl}{\pl t_1}
\end{array}
\right)
\gs_{\lb'}(t)\gs_{\mu'}(-s)\\
\\
-\gs_{\lb'}(t)\left(
\begin{array}{l}
V_{-1}(-s)+V_1(-s)+(z+z^{-1})V_0(-s)\\
+~n(-s_1)+(n+q)\frac{\pl}{\pl (-s_1)}
\end{array}
\right)\gs_{\mu'}(-s)
\end{array}\right\}\\
\\  \\
&=&\sum_{\lb',\mu'\atop{\lb_1^{'\top},\mu_1^{'\top}\leq
n}}a_{\lb'\mu'} \left\{
\begin{array}{l}
\left(
\begin{array}{l}
\displaystyle{\sum^{n+1}_1}(\lb'_i-i+n+1)\gs_{\lb'+e_i}(t)
 \\
 +\displaystyle{\sum^n_1}(\lb'_i-i+n+p)\gs_{\lb'-e_i}(t)
  \\   \\
+~(z+z^{-1})|\lb'|\gs_{\lb'}(t')
\end{array}
\right)
 \gs_{\mu'}(t)\\
\\
-\gs_{\lb'}(t)\left(
\begin{array}{l}
\displaystyle{\sum^{n+1}_1}(\mu'_i-i+n+1)\gs_{\mu'+e_i}(-s)
 \\  \\
 +\displaystyle{\sum^n_1}
    (\mu'_i-i+n+q)\gs_{\mu'-e_i}(-s)\\  \\
+(z+z^{-1})|\mu'|\gs_{\mu'}(-s)
\end{array}
\right)
\end{array}\right\}\\
\\
&=&\sum_{\lb,\mu\atop{\lb_1^{\top},\mu_1^{\top}\leq
n}}\gs_{\lb}(t)\gs_{\mu}(-s){\cal L}^{(3)}
(a_{\lb\mu})\eean
  }
implying (\ref{ExampleThreea}) and similarly
(\ref{ExampleThreeb}).\qed.

\subsection{The action of Virasoro on two-dimensional
Fourier series}

In this section we prove identity (\ref{box1}) of the
introduction:

\begin{corollary}
The action of the Virasoro operator $\tilde\BV_{\Lb_k}$
on``Fourier series", with arbitrary coefficients,
depending on an integer parameter $k$, translates itself
into a linear action of $\tilde L_{\Lb_k}$ on the
coefficients:%$\tilde b^{(k)}_{\lb\mu}$:
 \be
  \tilde\BV_{\Lb_k}
   \sum_
{
  { \lb,~\mu, ~\mbox{\tiny such that}}\atop
  {  \lb_1^{\top},~\mu_1^{\top}\leq n }
          }
    b^{(k)}_{\lb\mu}
    {\bf  s}_{\lb}(t){\bf s}_{\mu}(-s)
   =
   \sum_
{
  { \lb,~\mu, ~\mbox{\tiny such that}}\atop
  {  \lb_1^{\top},~\mu_1^{\top}\leq n }
          }
          \tilde L_{\Lb_k} (\tilde b^{(k)}_{\lb\mu}){\bf  s}_{\lb}(t)
            {\bf s}_{\mu}(-s)
    \ee
where
$$
\begin{array}{c|c}
\tilde\BV_{\Lb_k}(t,s) &  \tilde L_{\Lb_k}(\lb,\mu) \\
\hline
   \\
  \BV_0,~\BV_{\pm 1}\Bigr|_{{it_i\longmapsto it_i+k\Lb^{-1}_k
     \delta_{i1}}\atop
    {is_i\longmapsto is_i-k\Lb^{-1}_k\delta_{i1}}} &
        L^{(1)}_0,~
        L^{(1)}_{\pm}
      %\left\{\begin{array}{ll}
       %  L^{(1)}_+\\  L^{(1)}_-:=L^{(1)}_+\Bigr|_{x
       %  \leftrightarrow y}
        %     \end{array}\right.
     \\  \\

   \pm(\BV_0+\BV_{\pm 1})
    \Bigr|_{{it_i\longmapsto it_i-q(-1)^i}
       \atop {is_i\longmapsto is_i- k\Lb_k^{-1}\dt_{i1}}}

      &
       L^{(2)}_{\pm}
      \\  \\

   \left(\mathbb V_{-1}+(\Lb_k+\Lb_k^{-1})\mathbb
V_{0}+\mathbb V_{1} \right)\Bigr|_{{it_i\longmapsto
it_i-p(-\Lb_k)^{-i}}  \atop  {is_i\longmapsto
is_i+q(-\Lb_k)^{-i}}}
     &

     L^{(3)}

\end{array}
$$
\bigbreak

\noindent The right column of $\tilde L_{\Lb_k}$'s are
the precise equations satisfied by the three types of
random walks considered in this paper.

\end{corollary}

\proof
%
%Denote by $\tilde \BV_z$ the Virasoro operators for each
%of the shifts and by $\tilde \BV_{\Lb}$ the Virasoro
%operators, with $z$ replaced -roughly speaking- by
%$k\Lb_k^{-1}$, as made precise in table 2. Then, u
 Using
the commutation relation,
 $$ z\Lambda^{-1}_k \frac{z^{k}}{k!}
=\frac{z^{k}}{k!}
 k\Lambda^{-1}_k,
%
 %\frac{z^{k-1}}{(k-1)!}
$$
  one computes, on the one hand,
 \bean
  \tilde V_z
    \sum_{{\lb,\mu\mbox{\tiny ~such
that}}\atop{\lb_1^{\top},\mu_1^{\top}\leq
n}}a_{\lb\mu}(z)s_{\lb}(t)s_{\mu}(-s)
 & =&
  \tilde \BV_z
    \sum_{{\lb,\mu\mbox{\tiny ~such
that}}\atop{\lb_1^{\top},\mu_1^{\top}\leq n}}
\left(\sum_0^{\iy} \frac{z^k}{k!}
 \tilde b^{(k)}_{\lb\mu}\right)
 s_{\lb}(t)s_{\mu}(-s)
 \eean \bean\\
 %
 %\newpage
 &=& \tilde \BV_z
    \sum_0^{\iy} \frac{z^k}{k!}
    \sum_{{\lb,\mu\mbox{\tiny ~such
that}}\atop{\lb_1^{\top},\mu_1^{\top}\leq n}}
   \tilde b^{(k)}_{\lb\mu}
 s_{\lb}(t)s_{\mu}(-s)
 \\
 &=&
    \sum_0^{\iy} \frac{z^k}{k!}
   \tilde \BV_{\Lb}
     \sum_{{\lb,\mu\mbox{\tiny ~such
that}}\atop{\lb_1^{\top},\mu_1^{\top}\leq n}}
    \tilde b^{(k)}_{\lb\mu}
 s_{\lb}(t)s_{\mu}(-s)
 \eean
and, on the other hand,
\bean
  \lefteqn{
   \tilde\BV_z
    \sum_{{\lb,\mu\mbox{\tiny ~such
that}}\atop{\lb_1^{\top},\mu_1^{\top}\leq
n}}a_{\lb\mu}(z)s_{\lb}(t)s_{\mu}(-s)
       }\\
 &=&
   \sum_{{\lb,\mu\mbox{\tiny ~such
that}}\atop{\lb_1^{\top},\mu_1^{\top}\leq n}}
 {\cal L}(a_{\lb\mu}(z))s_{\lb}(t)s_{\mu}(-s)
  \\
 &=&
   \sum_{{\lb,\mu\mbox{\tiny ~such
that}}\atop{\lb_1^{\top},\mu_1^{\top}\leq n}}
 {\cal L}\left(\frac{z^k}{k!}b^{(k)}_{xy}\right)
    s_{\lb}(t)s_{\mu}(-s)
 \\
  &=&
   \sum_{k=0}^{\iy} \frac{z^k}{k!}\sum_{{\lb,\mu\mbox{\tiny ~such
that}}\atop{\lb_1^{\top},\mu_1^{\top}\leq n}}
 { L}_{\Lambda}\left(b^{(k)}_{xy}\right)
    s_{\lb}(t)s_{\mu}(-s),
  \eean
using the argument in (\ref{L-argument}).\qed

%\newpage

\section{The discrete backward and forward
equations for a random walk in a Weyl chamber}

 Remembering the definition of the difference
 operators (\ref{definition}) (see
footnote 4),
%
%\subsection{The discrete Brownian motion }
%Define, for $\al\in\BZ$, $\al\neq 0$, the following
%difference operators, acting on functions $f(k,x,y),$
%with $ k\in \BZ_+, ~x,y\in \BZ$:
%
%\bea
%\pl^+_{\al x_i}f &:=&f( k, x+\al e_i,y)-f(k,  x,y)  \no\\
%\pl^-_{\al x_i}f &:=&f(  k,x,y)-f( k, x-\al e_i,y) \no\\
%\Lb^{-1}_k f&:=& f(k-1,x,y) \label{definition}\eea
 consider now the following difference operators:
\bea
  {\cal A}_1&:=&
    \sum_{i=1}^n
   \left( \frac{k}{2n}\Lb^{-1}_k
     \pl^+_{2y_i}
     +
      x_i\pl^-_{x_i}
      +\pl^+_{y_i} y_i
   -
      (x_i-{y_i})
    \right)\no\\
 {\cal A}_2&:=&
    \sum_{i=1}^n
   \left( \frac{k}{2n}\Lb^{-1}_k
     \pl^+_{2x_i}
     +
      y_i\pl^-_{y_i}
      + \pl^+_{x_i} x_i
   -
      (y_i-{x_i})
    \right)\label{A-operator}
\eea

\begin{theorem} The probability
 $$
 P(k,x,y):=P\left( \begin{array}{l}
 \mbox{that $n$ walkers in $\BZ$,}\\
 \mbox{go from $x_1,\ldots,x_n$ to}\\
 \mbox{$y_1,\ldots,y_n$ in $k$ steps,}\\
 \mbox{and do not intersect}
 \end{array}
   \right)=\frac{b^{(k)}_{xy}}{(2n)^{k}}
   $$
satisfies
  \be
 {\cal A}_i P(k,x,y)=0
  \label{5.1}\ee

\end{theorem}

\proof Indeed
\bean
 (2n)^k {\cal A}_1 P(k,x,y)
  &=&
  {k}
  \sum_1^n (b^{(k-1)}_{x,y+2e_i}-b^{(k-1)}_{x,y})\\
  &&+
   \sum_1^n\Bigl(x_i\bigl( b_{xy}^{(k)}-b_{x-e_i,y}^{(k)}
               \bigr) -  %\sum_1^n
               x_i b_{xy}^{(k)}\Bigr)
         \\
   &&+
   \sum_1^n\Bigl(
   (y_i+1) \bigl( b_{x,y+e_i}^{(k)}-b_{x,y}^{(k)}
               \bigr) +  %\sum_1^n
               (y_i+1) b_{xy}^{(k)}
               \Bigr)\\
    \no\\
 &=&
 k\left(
  \sum_{1\leq i \leq n}   b^{(k-1)}_{x,y+2e_i}
  -nb^{(k-1)}_{x,y}
%  +
%  \sum_{1\leq i\leq n\atop{y_{i+1}-y_i=1}}
%  \sum_{1\leq i \leq n}\delta_{y_{i+1},y_i+1}
% b^{(k-1)}_{x,y+e_i+e_{i+1}}
 \right) \no  \\
 &&
  -\sum_{1\leq i\leq n}
  \Bigl(x_i b^{(k)}_{x-e_i,y} -(y_i+1) b^{(k)}_{x,y+e_i}\Bigr)
  \\
  &=&    -{ L}_{-}^{(1)}(x,y)~ b^{(k)}_{xy}=   0,
   \eean
 using (\ref{b-difference2}), insofar none of the
 final positions are adjacent. The second
 equation ${\cal A}_2 P(k,x,y)=0$ follows immediately by the
 duality $x\leftrightarrow y$. \qed

\end{document}